\documentclass{amsart}
\usepackage{float,mathrsfs,amsmath,amssymb,amsthm}
\usepackage[dvips]{graphicx}
\usepackage[all]{xy}
\usepackage{multirow,tabls}

\allowdisplaybreaks[4]

\numberwithin{equation}{section}

\theoremstyle{plain}
\newtheorem{problem}{Problem}
\newtheorem{theorem}{Theorem}
\newtheorem*{acknowledgement}{Acknowledgement}

\DeclareMathOperator{\Km}{Km}

\DeclareMathOperator{\I}{I}
\DeclareMathOperator{\II}{II}
\DeclareMathOperator{\III}{III}
\DeclareMathOperator{\IV}{IV}

\title[Weierstrass equations]{Weierstrass equations for Jacobian
  fibrations on a certain $K3$ surface}
\date{}
\author{Kazuki Utsumi}
\address{Department of Mathematics, Graduate School of Science, Hiroshima University}
\email{kazu-utsumi@hiroshima-u.ac.jp}
\subjclass[2010]{14J27, 14J28,14H52}
\keywords{$K3$ surface, elliptic surface, elliptic curve}

\begin{document}
        
\maketitle

\begin{abstract}
    In this paper we give the Weierstrass equations for Jacobian fibrations on the
    $K3$ surface that is the minimal resolution of the double
    covering of $\mathbb{P}^2$ ramified along generic six lines.
\end{abstract}

\section{Introduction}

\subsection{Problem setting}

Let $X$ be a $K3$ surface defined over an algebraically closed field
$k$ with ${\rm char}(k) \neq 2, 3$. Suppose $f: X \to
\mathbb{P}^1$ is a Jacobian fibration, that is, an elliptic fibration
on $X$ with a section $O:\mathbb{P}^1 \to X$. Let $t$ be an affine
coordinate of $\mathbb{P}^1$. Then $f^*(t)$ defines a non-constant
rational function on $X$, which is called an \textit{elliptic
  parameter} for the Jacobian fibration $f$. We also denote $f^* (t)$ by
$t$ and regard $t$ as a rational function on $X$. The generic fiber of
$f$ defines an elliptic curve $E$ over the rational function field
$k(t)$.

Kuwata and Shioda~\cite{Kuwa-Shio} proposed the following problems.

\begin{problem}
    Given a $K3$ surface $X/ k$ and a Jacobian fibration $f$,
    determine {\rm (i)} the elliptic parameter $t$ for $f$, {\rm
      (ii)} the defining equation of the elliptic curve
    $E/k(t)$, and {\rm (iii)} the Mordell-Weil lattice
    {\rm(}MWL{\rm )} $E(k(t))$.
\end{problem}

\begin{problem}
    Given a $K3$ surface $X/k$, determine all the {\rm(}essentially
    distinct{\rm)} elliptic parameters.
\end{problem}

Problem~2 is a combination of Problem~1 and the following standard problem:

\begin{problem}
    Given a $K3$ surface $X/k$, classify the Jacobian fibration $f:X
    \to \mathbb{P}^1$ up to isomorphism.
\end{problem}

Oguiso \cite{Oguiso} solved Problem 3 in the case where $X$ is a Kummer
surface of the product of non-isogenous elliptic curves over the
complex number field. Namely he classified the configuration of
singular fibers on such a Kummer surface $X$ into eleven types
$\mathcal{J}_1, \ldots, \mathcal{J}_{11}$, and determined the number
of the isomorphism classes for each type.

Kuwata and Shioda \cite{Kuwa-Shio} solved Problems 1 and 2 for each
member of Oguiso's list. They gave elliptic parameters and Weierstrass
equations for each member of Oguiso's list by using the Legendre
parameters of the two elliptic curves.

Recently, Kumar \cite{Kumar} solved the above problems completely for the
case of a generic Jacobian Kummer surface $X=\Km(J(C))$, for a genus 2
curve $C$ over an algebraically closed field of characteristic $0$. He
showed that there are exactly 25 distinct Jacobian fibrations on such
a generic Kummer surface, and gave an elliptic parameter, the
Weierstrass equation and the Mordell-Weil lattice for each type.

\subsection{Main results}

In this paper, we focus on the case of a $K3$ surface over $k$ that is
the minimal resolution of the double covering of $\mathbb{P}^2$
ramified along \emph{generic} six lines. When we say that six lines
are \emph{generic}, we mean that the rank of the N{\'e}ron-Severi
group $NS(X)$ is $16$.

In this case, Problem 3 has been partially solved by
Kloosterman~\cite{Kloosterman}. He classified all configurations of
singular fibers of Jacobian fibrations on such a $K3$ surface over a
field of characteristic $0$ (although most of his results hold over
fields of any characteristic, not $2$ or $3$) into sixteen
classes. However he did not give the classification of the isomorphism
classes for each type. Thus, this gives a partial solution to Problem 3
on such a $K3$ surface.

Table~\ref{tb:list} shows a summary of Kloosterman's results in the
generic case. The first column shows the class of Jacobian fibration
following Kloosterman's notation. The second column shows the
configuration of singular fibers. Here, for example, by $\I_2^*+8\I_2$
we mean that the surface has two fibers of type $\I_2^*$ (Kodaira's notation
\cite{Kodaira}) and eight fibers of type $\I_2$. The third
column shows the Mordell-Weil group(MWG) of the fibration.

\begin{table}[htbp]
    \centering
    \caption{List of possible configurations}\label{tb:list}
    \begin{tabular}{|c|c|c|}\hline
        Class & Configuration of singular fibers & MWG \\ \hline
        1.1 & $\I_{10}+ \I_{2}+ a \II+ b\I_1$& $\mathbb{Z}^4$ \\ \hline
        1.2 & $\I_{8}+ \I_{4}+ a \II+ b\I_1$&  $\mathbb{Z}^4$ \\  \hline
        1.3 & $2\I_{6} +  a \II+b\I_1$& $\mathbb{Z}^4$ \\ \hline
        1.4 & $\IV^*+ \I_{4}+  a \II+ b\I_1$& $\mathbb{Z}^5$  \\ \hline 
        2.1 & $\II^* + 6\I_2 + 2\I_1$ & $\{0\}$ \\ \hline
        2.2 & $\III^* +7\I_2 + \I_1$ & $\mathbb{Z}/2\mathbb{Z}$ \\ \hline
        2.3 & $\III^* + \I_0^* + 3\I_2 + 3\I_1$ & $\{0\}$ \\ \hline
        2.4 & $\I_6^* + 4\I_2 + 4\I_1$ & $\{0\}$ \\ \hline
        2.5 & $\I_4^* + 6\I_2 + 2\I_1$ & $\mathbb{Z}/2\mathbb{Z}$ \\ \hline 
        2.6 & $\I_4^* + \I_0^* + 2\I_2 + 4\I_1$ & $\{0\}$ \\ \hline
        2.7 & $\I_2^* + 8\I_2$ & $(\mathbb{Z}/2\mathbb{Z})^2$ \\ \hline
        2.8 & $\I_2^* + \I_0^* + 4\I_2$ & $\mathbb{Z}/2\mathbb{Z}$ \\ \hline
        2.9 & $2\I_2^* + 2\I_2 + 4\I_1$ & $\{0\}$ \\ \hline
        2.10 & $\I_2^* + 2\I_0^* + 8\I_1$ & $\{0\}$ \\ \hline
        2.11 & $2\I_0^* + 6\I_2$ & $(\mathbb{Z}/2\mathbb{Z})^2$ \\ \hline
        2.12 & $3\I_0^* + 2\I_2 + 2\I_1$ & $\mathbb{Z}/2\mathbb{Z}$  \\ \hline
    \end{tabular}
\end{table}

Our main results are as follows: we solve (i) and (ii) of Problem 1
for each class of Table~\ref{tb:list}. We give an elliptic parameter
and its Weierstrass equation for one Jacobian fibration in each class
of the table, although there may exist nonisomorphic Jacobian
fibrations belonging to the same class of the list. More details will
be given in \S~\ref{sec:notations} after we fix the notation.

\subsection{Notation}\label{sec:notations}

Fix generic six lines $L_i \subset \mathbb{P}^2$. Denote by $P_{i,j}$
the point of intersection of $L_i$ and $L_j$.

Let $\varphi': Y \to \mathbb{P}^2$ be the double cover ramified along
the six lines $L_i$. Then $Y$ has 15 double points of type $A_1$,
which correspond to $P_{i,j}$. Blowing up these points gives a $K3$
surface $X$, with 15 exceptional divisors $\ell_{i,j}$ and a rational
map $\varphi:X \to \mathbb{P}^2$. For a curve $C$ on $\mathbb{P}^2$,
we call the strict transform of $\varphi'^*(C)$ \textit{the pull-back
  of $C$}, for short. Let $\ell_i$ be the divisor on $X$ such that
$2\ell_i$ is the pull-back of $L_i$. Let $\mu^{i,j}_{k,m}$ be the
pull-back of the line $M^{i,j}_{k,m}$ connecting $P_{i,j}$ and
$P_{k,m}$ with $i,j,k,m$ pairwise distinct. With this notation, which
is the same as Kloosterman~\cite{Kloosterman}, divisors $\ell_i,
\ell_{i,j}, \mu^{i,j}_{k,m}$ are $(-2)$-curves on $X$. We have the
following intersection numbers.

\begin{equation}
    \begin{aligned}
        &\ell_{i} \cdot \ell_{j}=
        \begin{cases}
            -2 & i=j\\
            0 & i \neq j
        \end{cases}, \quad \ell_{i,j} \cdot \ell_{k,m}=
        \begin{cases}
            -2 & \{i,j\} = \{k,m\}\\
            0 & \text{otherwise}
        \end{cases},\\
        & \\
        & \ell_i \cdot \ell_{k,m}=
        \begin{cases}
            1 & i \in \{k,m \}\\
            0 & \text{otherwise}
        \end{cases}, \quad \ell_p \cdot \mu^{i,j}_{k,m}=
        \begin{cases}
            1 & p \notin \{ i,j,k,m \}\\
            0 & \text{otherwise}\\
        \end{cases},\\
        & \\
        & \ell_{i,j} \cdot \mu^{k,m}_{p,q}=
        \begin{cases}
            2 & \{i,j\}=\{k,m\} \, \text{ or } \, \{i,j\}=\{p,q\}\\
            0 & \text{otherwise}\\
        \end{cases},\\
        & \\
        & \mu^{i,j}_{k,m} \cdot \mu^{p,q}_{r,s}=
        \begin{cases}
            -2 &  \{ \{i,j\}, \, \{k,m\} \} = \{ \{p,q\}, \, \{r,s\} \}\\
            2 &\{ \{i,j\}, \, \{k,m\} \} \cap \{ \{p,q\}, \, \{r,s\} \} = \varnothing  \\
            0 &  \text{otherwise} \\
        \end{cases}.
    \end{aligned}
\end{equation}

Moreover, for some rational plane curves $C$, we name the pull-back of
$C$ as in the following table. Note that all of them are $(-2)$-curves on
$X$.

\begin{center}
    \begin{tabular}{|c|c|}\hline
        Divisor & $C$\\ \hline
        $\eta^{(i_1 j_1)(i_2 j_2)(i_3 j_3)}_{(i_4 j_4)(i_5 j_5)}$ & the conic curve through $P_{i_1, j_1}, \ldots, P_{i_5, j_5}$\\ \hline 
        \multirow{2}{*}{$\xi^{(\overline{i_1 j_1})(i_2 j_2)(i_3 j_4)}_{(i_5 j_5)(i_6 j_6)(i_7 j_7)}$} 
        &\multirow{2}{7cm}{the cubic curve through $P_{i_1, j_1}, \ldots, P_{i_7, j_7}$ with a double point at $P_{i_1,j_1}$}\\
        & \\ \hline
        \multirow{2}{*}{$\nu^{(\overline{i_1 j_1})(\overline{i_2 j_2})(\overline{i_3 j_3})(i_4 j_4)}_{(i_5 j_5)(i_6 j_6)(i_7 j_7)(i_8 j_8)}$}
        &\multirow{2}{7cm}{the quartic curve through $P_{i_1, j_1}, \ldots, P_{i_8, j_8}$ with a double point at
          $P_{i_1, j_1}, P_{i_2, j_2}, P_{i_3, j_3}$}\\
        & \\\hline 
        \multirow{2}{*}{$\gamma^{(\overline{\overline{i_1 j_1}})(\overline{i_2 j_2})(\overline{i_3 j_3})(\overline{i_4 j_4})}_
          {(i_5 j_5)(i_6 j_6)(i_7 j_7)(i_8 j_8)(i_9 j_9)}$}
        &\multirow{3}{7cm}{the quintic curve through $P_{i_1, j_1}, \ldots, P_{i_9, j_9}$ with a triple point at $P_{i_1, j_1}$
          and a double point at $P_{i_2, j_2}, P_{i_3, j_3}, P_{i_4, j_4}$}\\ 
        & \\
        & \\ \hline
    \end{tabular}
\end{center}

We may suppose that the lines $L_i$ are defined by the following equations

\begin{equation}
    \label{eq:Li}
    \begin{aligned}
        &L_1 \; : \; u=0, \quad L_2 \; : \; u=1, \quad L_3 \; : \; v=0, \quad L_4 \; : \; v=1\\
        &L_5 \; : \; au+bv-1=0, \quad L_6 \; : \; cu+dv-1=0, 
    \end{aligned}
\end{equation}
where $u, \, v$ are the affine parameters of projective plane.
\begin{figure}[htbp]
    \centering
    \includegraphics[width=9cm]{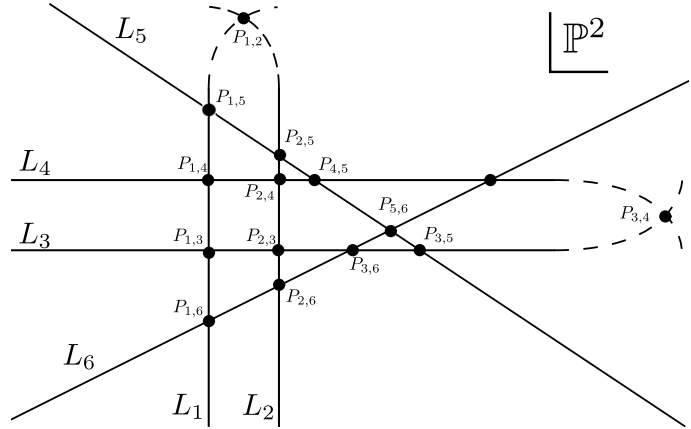}
    \caption{six lines on projective plane}
    \label{fig:six_lines}
\end{figure}
We consider the six lines in the \textit{generic case}, that is, the
rank of the N{\'e}ron-Severi group $NS(X)$ is 16.

Then the singular affine model of $X$ is given by

\begin{equation}
    \label{eq:DoubleCover}
    w^2 = u (u-1) v (v-1) (au+bv-1) (cu+dv-1).
\end{equation}
Under the above notation, we see that the divisors of typical
functions are as follows.

\begin{equation}\label{eq:divisors}
  \begin{aligned}
      &(u) = 2 \ell_1+\ell_{1,3}+\ell_{1,4}+\ell_{1,5}+\ell_{1,6}
      - \left( \mu^{1,2}_{3,4}+\ell_{3,4}  \right)\\
      &(u-1) = 2 \ell_2+\ell_{2,3}+\ell_{2,4}+\ell_{2,5}+\ell_{2,6}
      - \left(  \mu^{1,2}_{3,4}+\ell_{3,4}  \right)\\
      &(v) = 2 \ell_3+\ell_{1,3}+\ell_{2,3}+\ell_{3,5}+\ell_{3,6}
      - \left( \mu^{1,2}_{3,4}+\ell_{1,2}  \right)\\
      &(v-1) = 2 \ell_4+\ell_{1,4}+\ell_{2,4}+\ell_{4,5}+\ell_{4,6}
      - \left( \mu^{1,2}_{3,4}+\ell_{1,2}  \right)\\
      &(au+bv-1) = 2
      \ell_5+\ell_{1,5}+\ell_{2,5}+\ell_{3,5}+\ell_{4,5}+\ell_{5,6}
      -\left( \mu^{1,2}_{3,4}+\ell_{1,2}+\ell_{3,4} \right)\\
      & (cu+dv-1) =
      2\ell_6+\ell_{1,6}+\ell_{2,6}+\ell_{3,6}+\ell_{4,6}+\ell_{5,6}
      -\left( \mu^{1,2}_{3,4}+\ell_{1,2}+\ell_{3,4} \right)\\
      &(w) = \sum_{i=1}^{6}{\ell_i} \, + \, \sum_{1 \leq i < j \leq
        6}\ell_{i,j} - 3 \left( \mu^{1,2}_{3,4}+\ell_{1,2}+\ell_{3,4}
      \right)\\
      &\left( \text{the equation of } M^{i,j}_{k,m} \right) =
      \mu^{i,j}_{k,m}+ \ell_{i,j}+\ell_{k,m} - \left( \mu^{1,2}_{3,4}
        + \ell_{1,2}+\ell_{3,4} \right).
  \end{aligned}
\end{equation}

For each class of Table~\ref{tb:list}, we compute a Weierstrass
equation for one Jacobian fibration belonging to the
class. Theoretically, constructing a Jacobian fibration on a $K3$
surface is to find a divisor that has the same type as a singular
fiber in the Kodaira's list (see \cite{Kodaira}). In practice,
however, we need to find two divisors, one for the fiber at $t=0$, and
the other for the fiber at $t=\infty$, to write down an actual
elliptic parameter. Once an elliptic parameter is found, we would like
to find a change of variables that converts to the defining equation
to a Weierstrass form. In most cases, we encounter an equation of the
form $y^2=$ (quartic polynomial). Then we can transform it to a
Weierstrass form by using a standard algorithm (see for example
\cite{Connell} or \cite{Jacobians}).

In our case, we use two methods to convert defining equation
(\ref{eq:DoubleCover}) to the form $y^2=$ (quartic polynomial). The
first method is an elimination. Since an elliptic parameter $t$ is a
rational map, we can put $t=f/g$ for some $f, g \in k[u,v,w]$. Thus,
we can eliminate one variable from (\ref{eq:DoubleCover}) and the
equation $gt-f=0$. If such an equation can be converted to $y^2=$
(quartic polynomial) by a simple coordinate change, we can get a
Weierstrass equation. We will call this method {\it classical method}
in this paper.

The other method is a {\it 2-neighbor step}, which is designed by Noam
Elkies. This is the technique to transform a Weierstrass equation of a
Jacobian fibration to a Weierstrass equation of a distinct Jacobian
fibration. Using this, we can get an unknown Weierstrass equation of a
class from a known class. We describe a 2-neighbor step in
\S~\ref{sec:2-neighbor}.

\subsection{Results}\label{sec:results}

We state our main theorem. 

\begin{theorem}
    Let $X$ be a $K3$ surface over an algebraically closed field $k$
    with ${\rm char}(k) \neq 2,3$ that is the minimal resolution of
    the double covering of $\mathbb{P}^2$ ramified along generic six
    lines. Suppose that the rank of the N{\'e}ron-Severi group $NS(X)$
    is $16$. Under the singular affine model (\ref{eq:DoubleCover}) of
    $X$, for each class in Table~\ref{tb:list}, an elliptic parameter
    and a Weierstrass equation of a Jacobian fibration belonging to
    the class is given by Table~\ref{tb:11}, \ldots,
    Table~\ref{tb:212}.
\end{theorem}

In the case of using a classical method, we
give an elliptic parameter $t$ in terms of $u, v, w$, a Weierstrass
equation and a picture of the configuration of singular
fibers. Moreover, for the class 2.xx, we also give the correspondence
between the divisors and the torsion sections.

In the case of using a 2-neighbor step from another class (in fact, we
only use a 2-neighbor step from the classes 2.7 or 2.5), we give an
elliptic parameter $s$ in terms of $t, x, y$ used in the equation of
the source class, a picture that shows the way to construct of the
divisor corresponding to the fiber at $s=\infty$ and a picture of the
configuration of singular fibers. In this case, however, we omit
Weierstrass equations, since they are all too long to print in this
section. We give them in \S~\ref{sec:WEQ}.

We explain the detail of some computation. In \S~\ref{sec:class2.7},
we will give an elliptic parameter and a Weierstrass equation of a
Jacobian fibration of the class 2.7 by a classical method. In
\S~\ref{sec:2-neighbor}, we explain a 2-neighbor step. In
\S~\ref{sec:class2.10}, we will give an elliptic parameter and a
Weierstrass equation of the class 2.10 by a 2-neighbor step from the
class 2.7. In \S~\ref{sec:class2.4}, we use a 2-neighbor step from the
class 2.5 for the computation of the class 2.4.

\clearpage

\begin{table}
    \centering
    \begin{tabular}{|c|c|}\hline
        \multicolumn{2}{|c|}{Class 1.1}\\ \hline \hline
        Method & 2-neighbor step from the class 2.5\\ \hline
        Elliptic parameter & $s=\dfrac{y}{t x (act-t-1)}$ \\ \hline
        \multicolumn{2}{|c|}{\includegraphics[width=9cm]{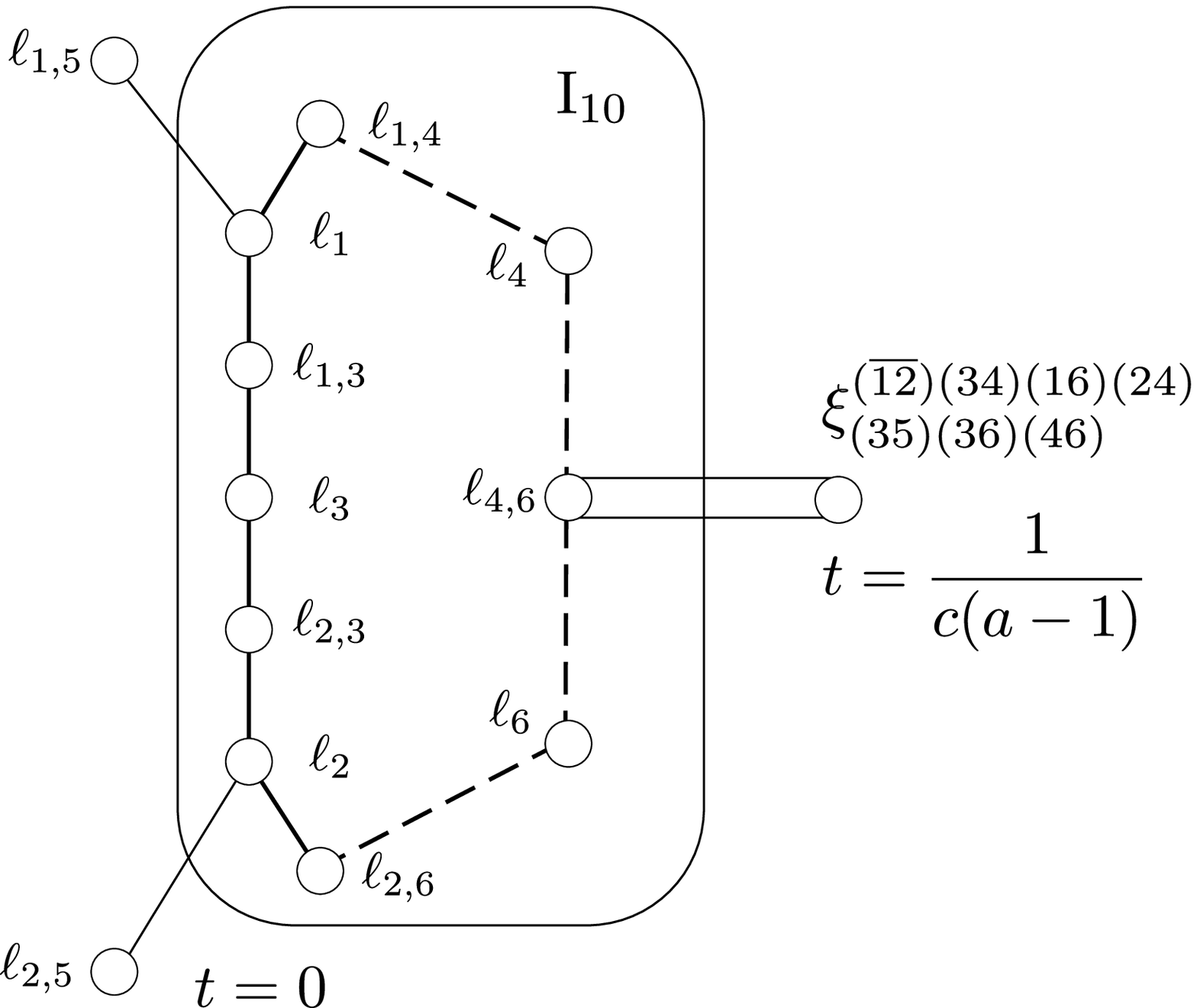} }\\ \hline
        \multicolumn{2}{|c|}{\includegraphics[width=9cm]{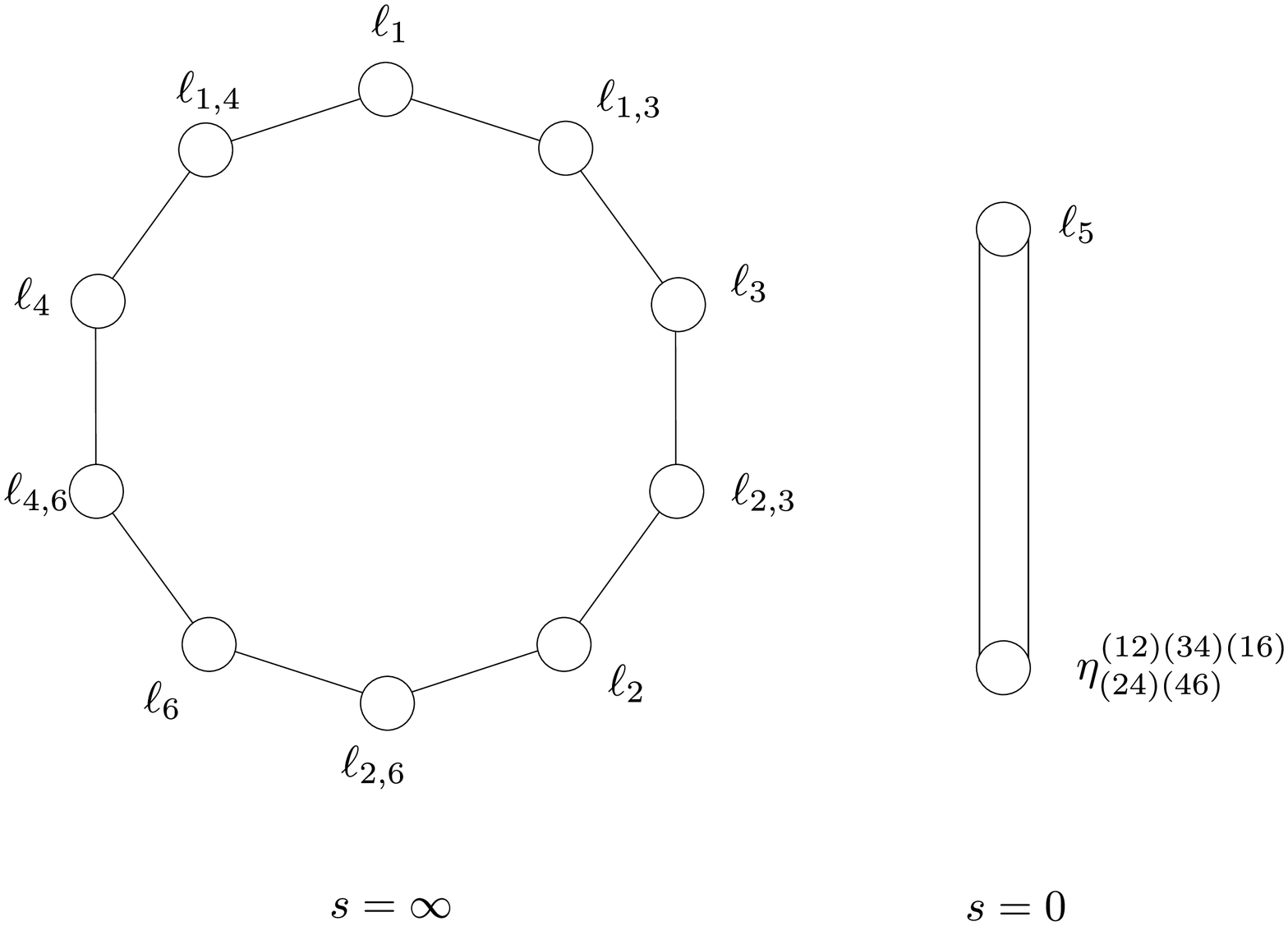} }\\ \hline        
    \end{tabular}
    \vspace{.1in}
    \caption{Class 1.1}\label{tb:11}
\end{table}

\begin{table}[htbp]
    \centering
    \tablinesep=15pt
    \begin{tabular}{|c|c|} \hline \multicolumn{2}{|c|}{Class 1.2}\\ \hline \hline
        Method & Classical\\ \hline
        Elliptic parameter & $t=\dfrac{w}{(au+bv-1)(cu+dv-1)}$ \\ \hline
        \multicolumn{2}{|c|}{\includegraphics[width=9cm]{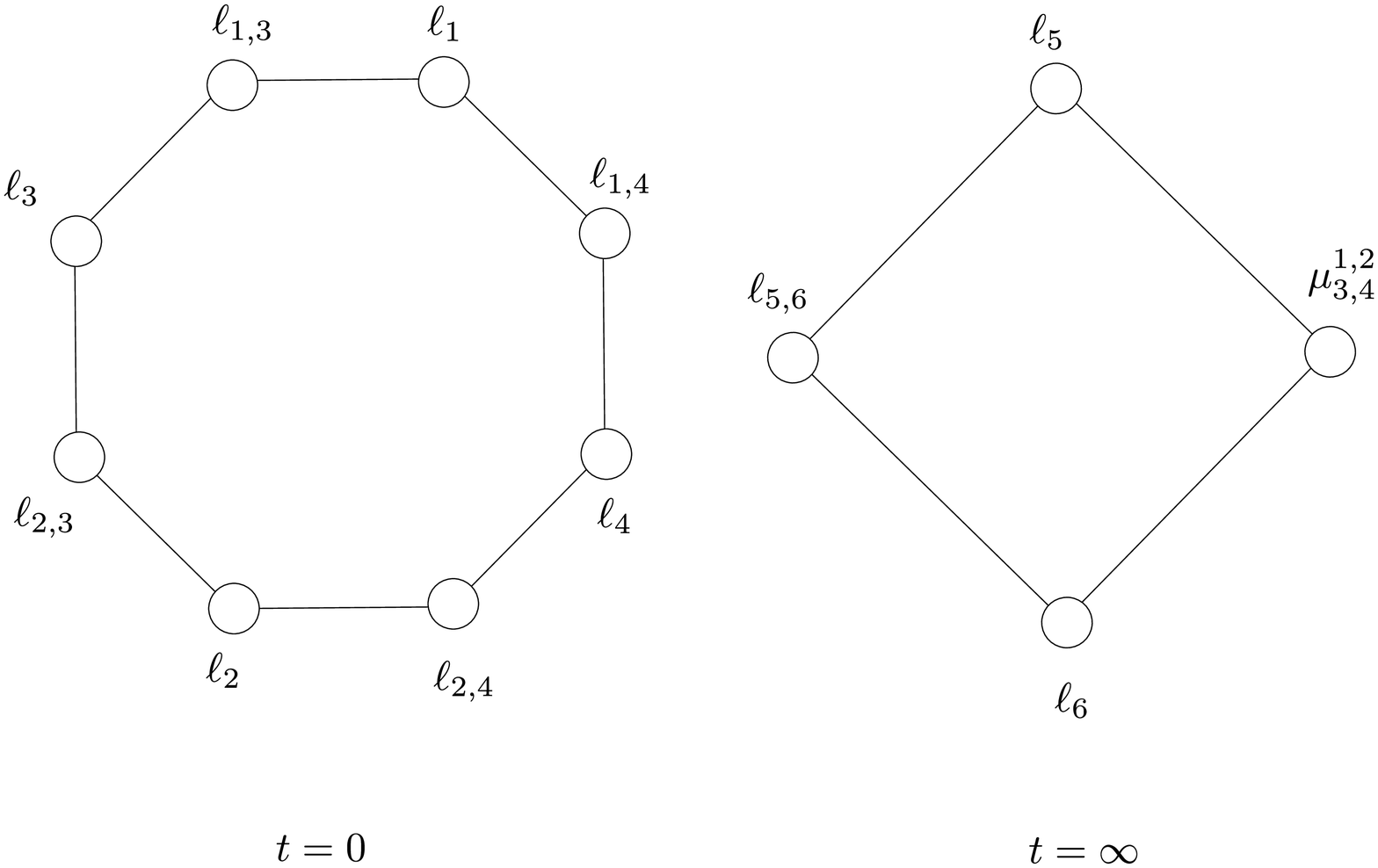}}\\ \hline
        \multicolumn{2}{|c|}{\small
          $
          \begin{aligned}
              &{y}^{2} + \left( 2\, \left( b-d \right) \left( ad-bc
                \right)
                {t}^{2}-4+2\,b+2\,d \right) xy\\
              & \quad + 4(b-d)^4 t^2 \left( \left(2\,a+2\,ac+2\,c
                  -ad-bc\right) {t}^{2}+1 \right) y\\
              & = x^3 -2 \left( (b-d)(
                {b}^{2}+2\,abd-2\,bcd-2\,b-2\,ab+2\,cd -{d}^{2}+2\,d) {t}^{2} \right.\\
              & \qquad \left. +2\, \left( d-1 \right) \left( b-1
                \right)
              \right)x^2 -4  \,  t^4 \, (b-d)^6 \, (4act^2+1) x\\
              & \quad + 8 \, t^4 \, (b-d)^6 \, (4act^2+1) \, \left(
                (b-d)(
                {b}^{2}+2\,abd-2\,bcd-2\,b-2\,ab \right.\\
              & \qquad \left.+2\,cd-{d} ^{2}+2\,d) {t}^{2}+2\, \left(
                  d-1 \right) \left( b-1 \right) \right)
          \end{aligned}
          $}\\ \hline
    \end{tabular}
    \vspace{.1in}
    \caption{Class 1.2}\label{tb:12}
\end{table}

\begin{table}[htbp]
    \centering
    \begin{tabular}{|c|c|} \hline 
        \multicolumn{2}{|c|}{ Class 1.3}\\ \hline \hline
        Method & Classical \\ \hline
        Elliptic parameter & $t=\dfrac{w}{v (u-1) (v-1)}$ \\ \hline
        \multicolumn{2}{|c|}{\includegraphics[width=9cm]{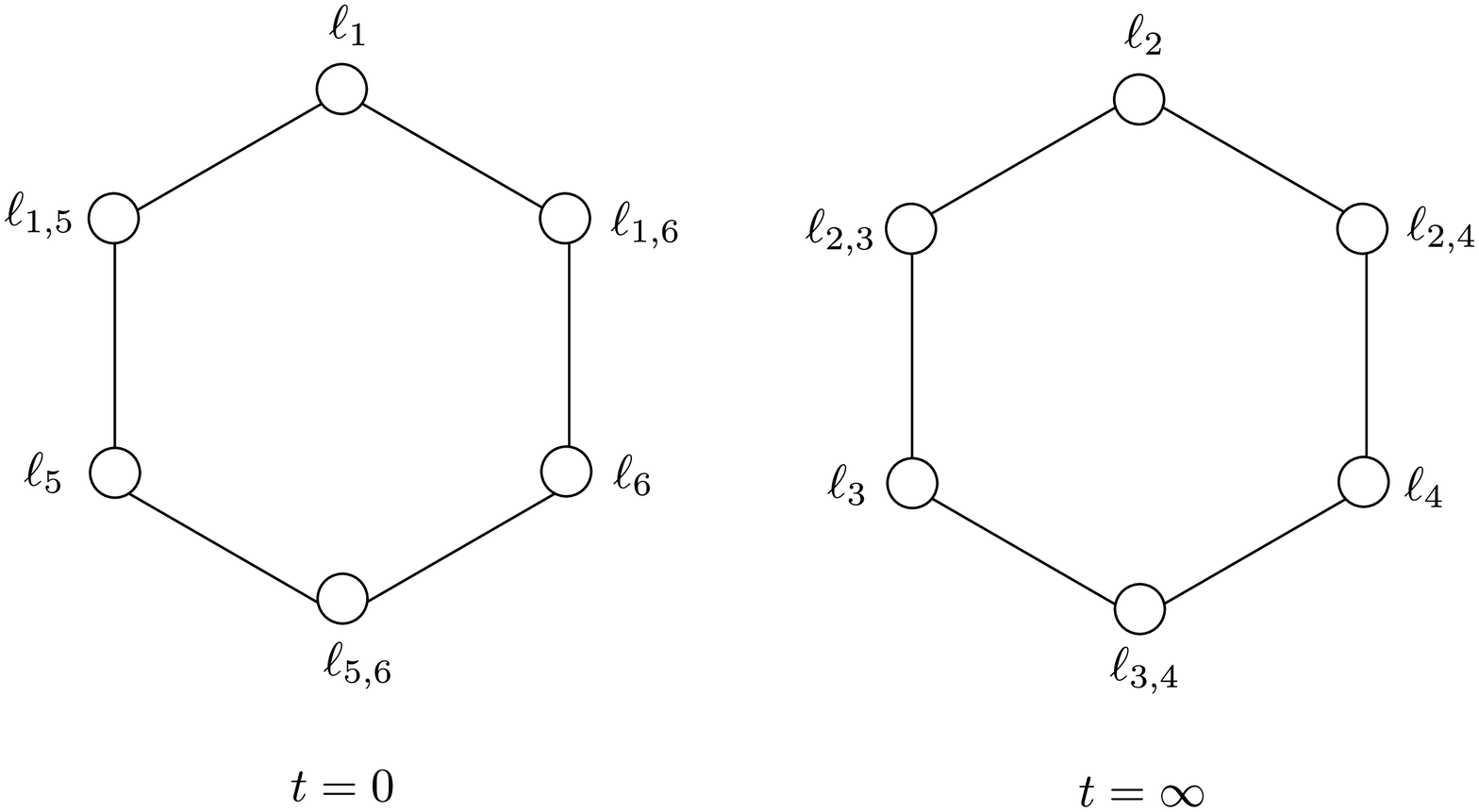}}\\ \hline
        \multicolumn{2}{|c|}{
          $
          \begin{aligned}
              {t}^{2}v-u+ \left( d+b-{t}^{2} \right) uv &+ \left(
                {t}^{2}-bd \right)
              u{v}^{2}- \left( bc+ad \right) {u}^{2}v\\
              &-{t}^{2}{v}^{2}-ac{u}^{3}+\left( a+c \right) {u}^{2}=0
          \end{aligned}
          $}\\
        \multicolumn{2}{|c|}{
          \multirow{2}{*}{This converts to a Weierstrass form (see \cite{Connell} or \cite{Jacobians}).}
        }\\ \hline
    \end{tabular}
    \vspace{.1in}
    \caption{Class 1.3}\label{tb:13} 
\end{table}

\begin{table}[htbp]
    \centering
    \begin{tabular}{|c|c|} \hline
        \multicolumn{2}{|c|}{ Class 1.4}\\ \hline \hline
        Method & 2-neighbor step from the class 2.7\\ \hline
        Elliptic parameter &
          $
          s=\dfrac{y}{t^2 \left(x-a t (bt-dt-1)(bt+ct-t-1) \right)}
          $\\ \hline
        \multicolumn{2}{|c|}{\includegraphics[width=9cm]{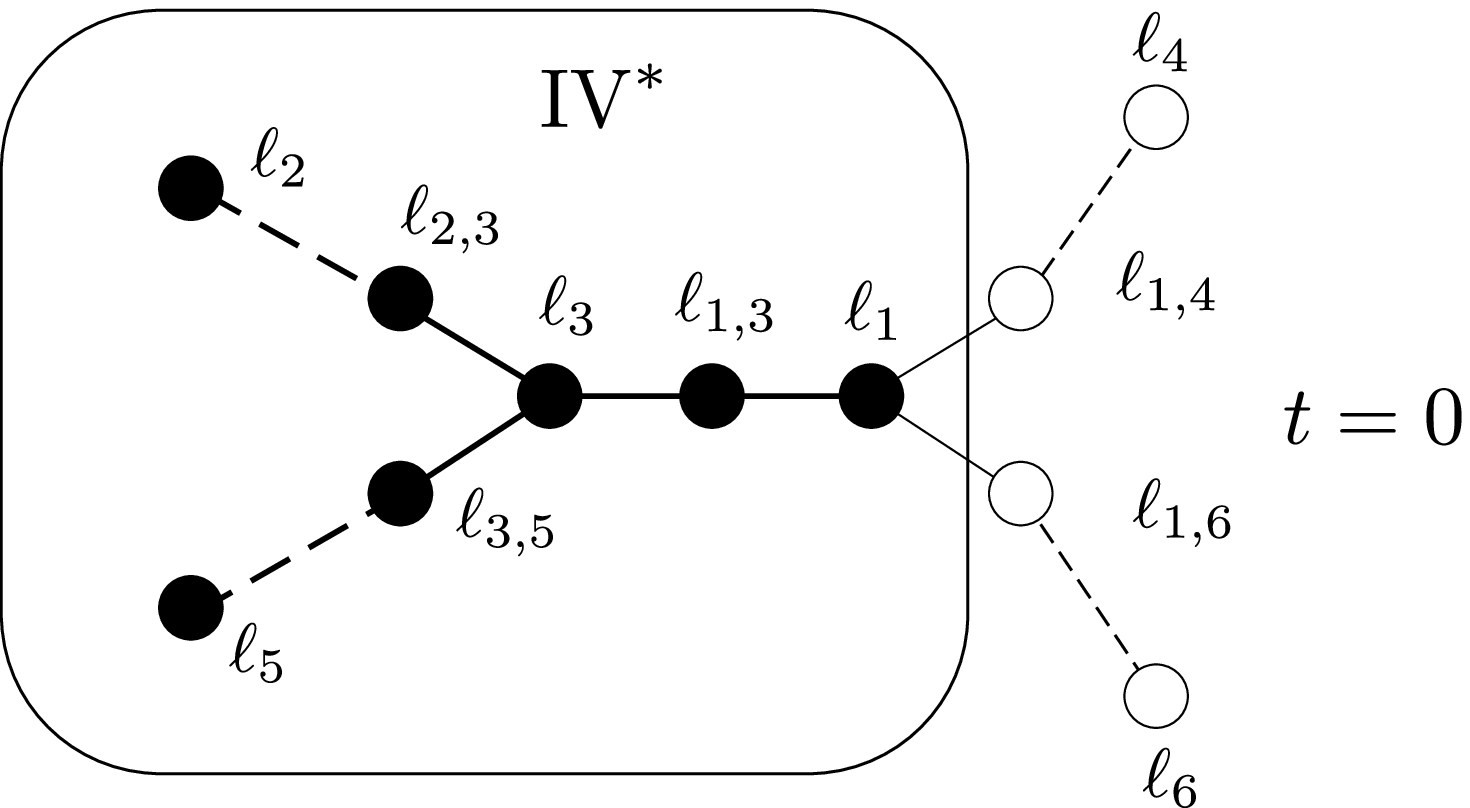}}\\ \hline
        \multicolumn{2}{|c|}{\includegraphics[width=9cm]{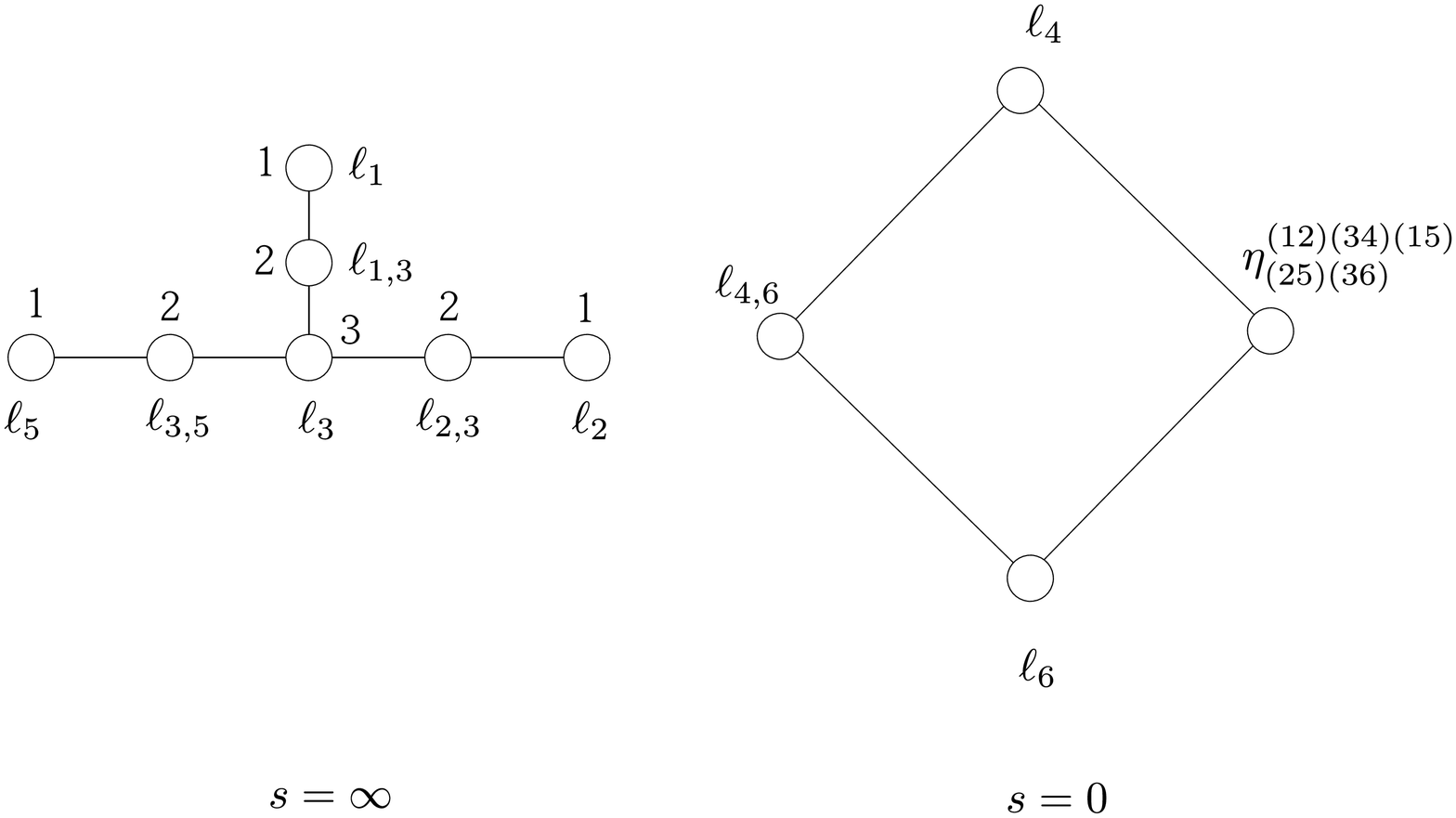}}\\ \hline
    \end{tabular}
    \vspace{.1in}
    \caption{Class 1.4}\label{tb:14}    
\end{table}

\begin{table}[htbp]
    \centering
    \begin{tabular}{|c|c|}\hline
        \multicolumn{2}{|c|}{Class 2.1}\\ \hline \hline
        Method & 2-neighbor step from the class 2.5 \\ \hline
        Elliptic parameter & $s= \dfrac{x-(b-1)(ad-bc+b-d)t^3}{t^4}$\\ \hline
         \multicolumn{2}{|c|}{\includegraphics[width=9cm]{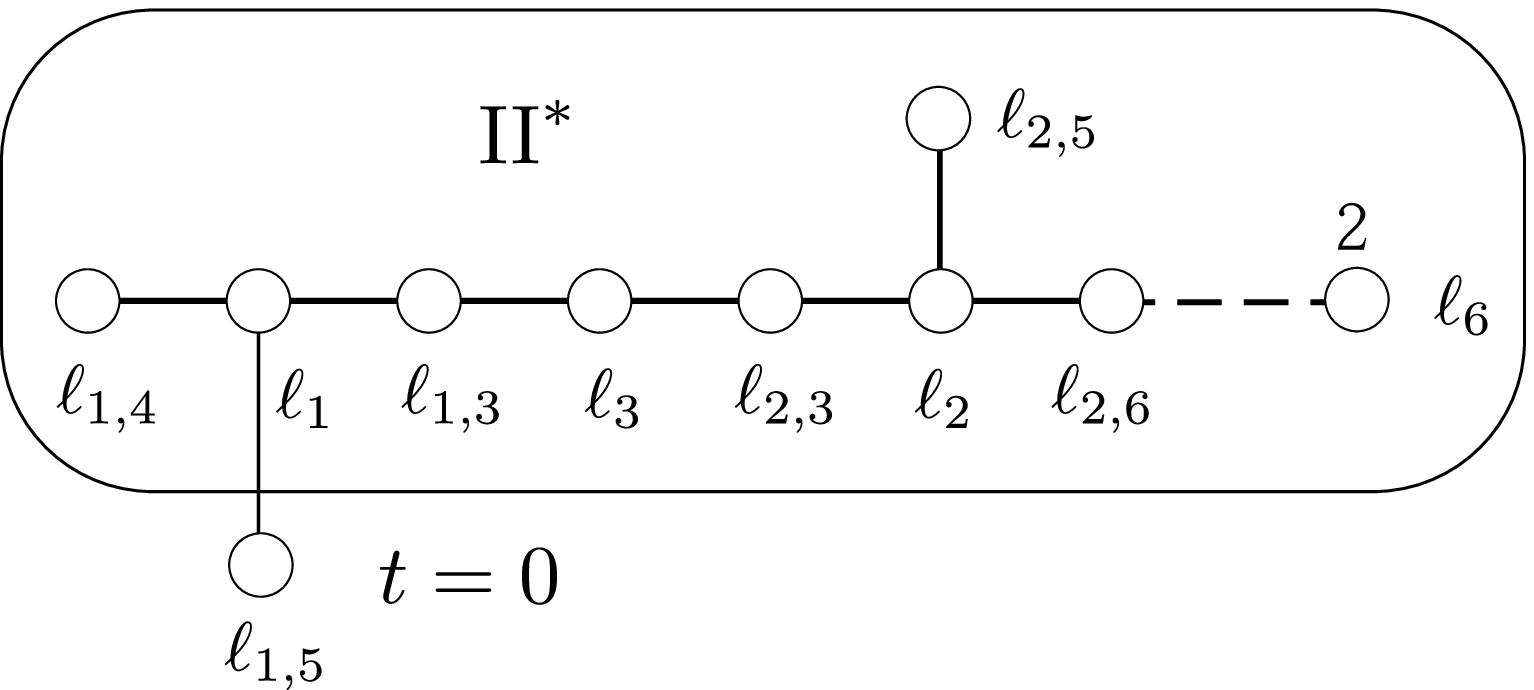}}\\ \hline
         \multicolumn{2}{|c|}{\includegraphics[width=9cm]{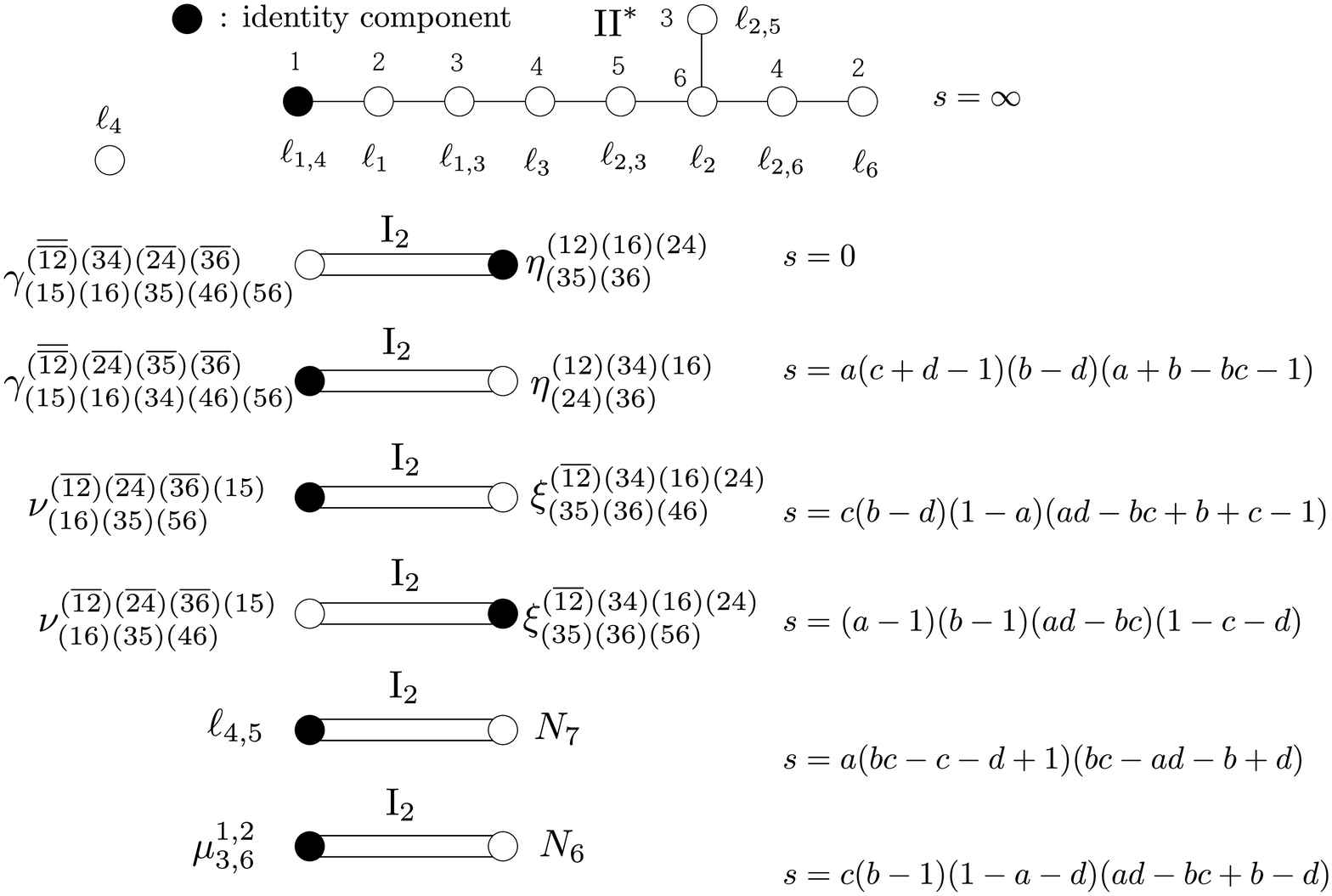}}\\
         \multicolumn{2}{|c|}{$N_i$ : the pull-back of a rational  plane curve of degree $i$}\\ \hline
    \end{tabular}
    \vspace{.1in}
    \caption{Class 2.1}\label{tb:21}    
\end{table}

\begin{table}[htbp]
    \centering
    \begin{tabular}{|c|c|} \hline
        \multicolumn{2}{|c|}{ Class 2.2}\\ \hline \hline
        Method & 2-neighbor step from the class 2.7 \\ \hline
        Elliptic parameter & {
          $
          \begin{aligned}
             & s = \dfrac{x+(1-a) t + A_2 t^2 - A_3 t^3}{t^3 (bt-dt-1)},\\
              &  A_2 = 2ab-bc+ac-a-b-c+1\\
              &  A_3 = (b-d)(ab-bc+ac+ad-a-c-d+1)
          \end{aligned}
          $}\\ \hline
        \multicolumn{2}{|c|}{\includegraphics[width=9cm]{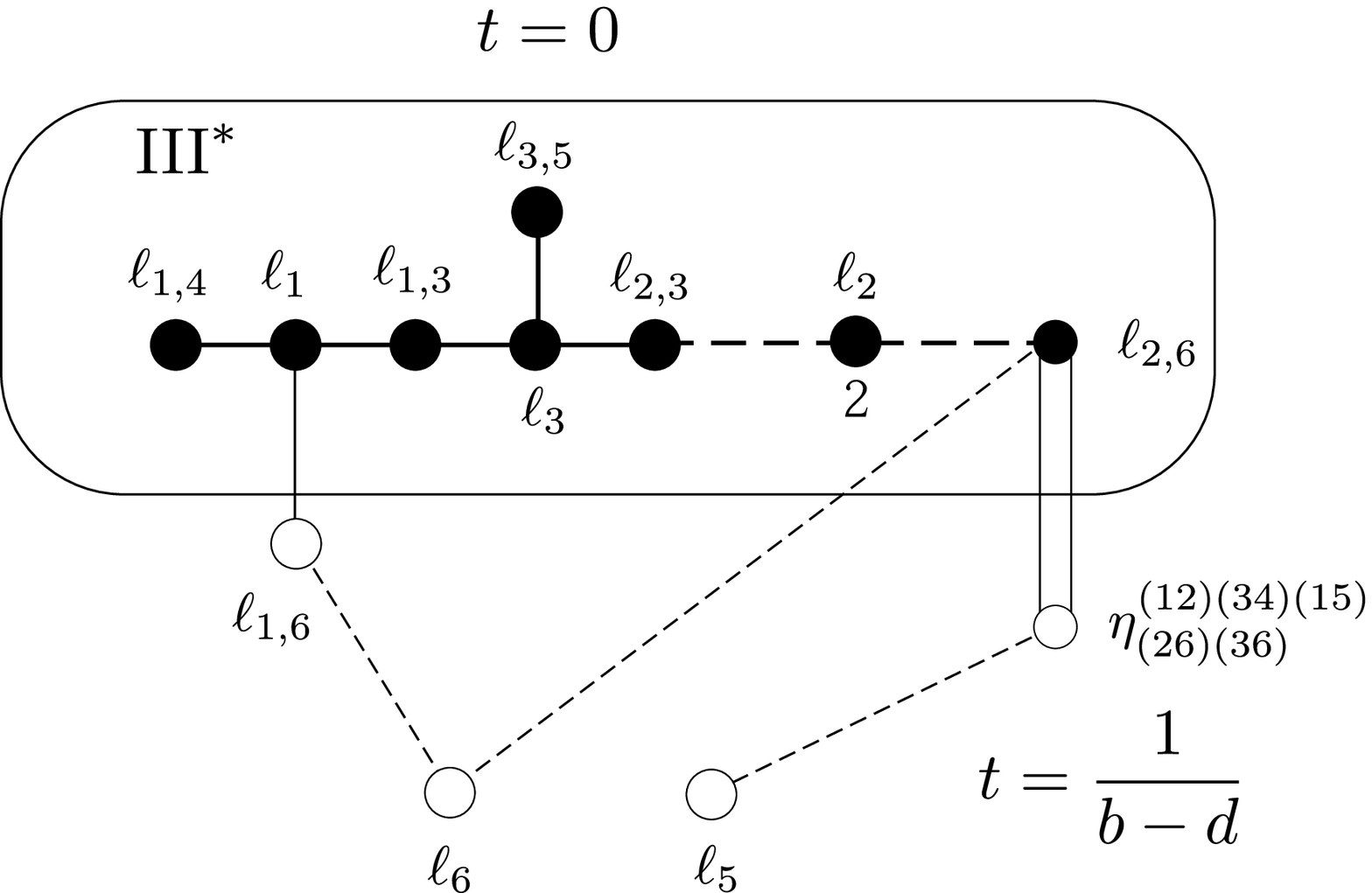}}\\ \hline
        \multicolumn{2}{|c|}{\includegraphics[width=9cm]{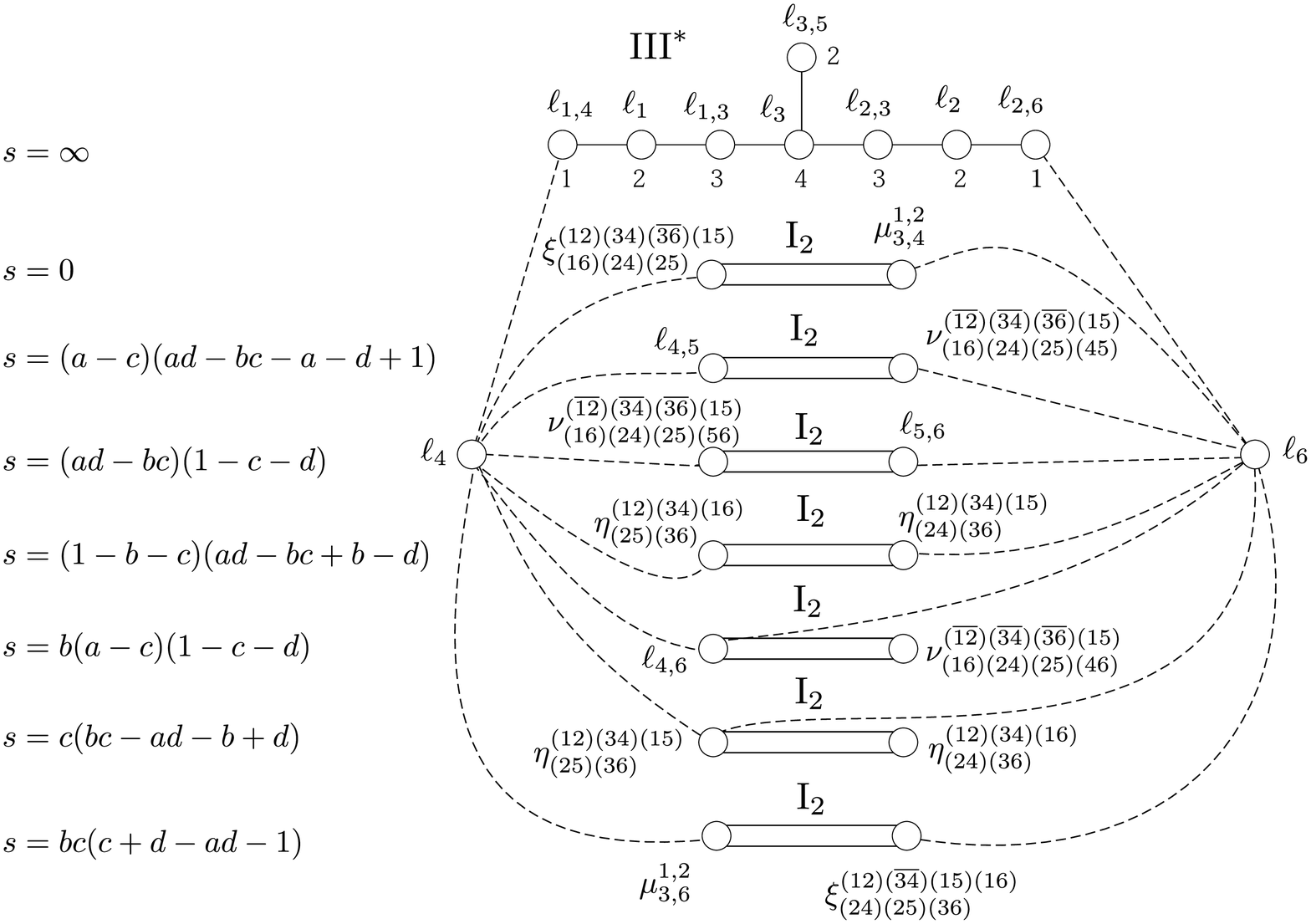}}\\ \hline
    \end{tabular}
    \vspace{.1in}
    \caption{Class 2.2} \label{tb:22}   
\end{table}

\begin{table}[htbp]
    \centering
    \begin{tabular}{|c|c|} \hline
        \multicolumn{2}{|c|}{ Class 2.3}\\ \hline \hline
        Method & 2-neighbor step form the class 2.7\\ \hline
        Elliptic parameter & 
        {\normalsize $
          \begin{aligned}
             & s = \dfrac{x - (a-1) t + A_2 t^2 - b(c+b-1)(a-c) t^3}{t^3 \left( b(a-c)t-a+1 \right)},\\ 
             & A_2 = 2ab-bc+ac-a-b-c+1
          \end{aligned}
          $}\\ \hline
        \multicolumn{2}{|c|}{\includegraphics[width=9cm]{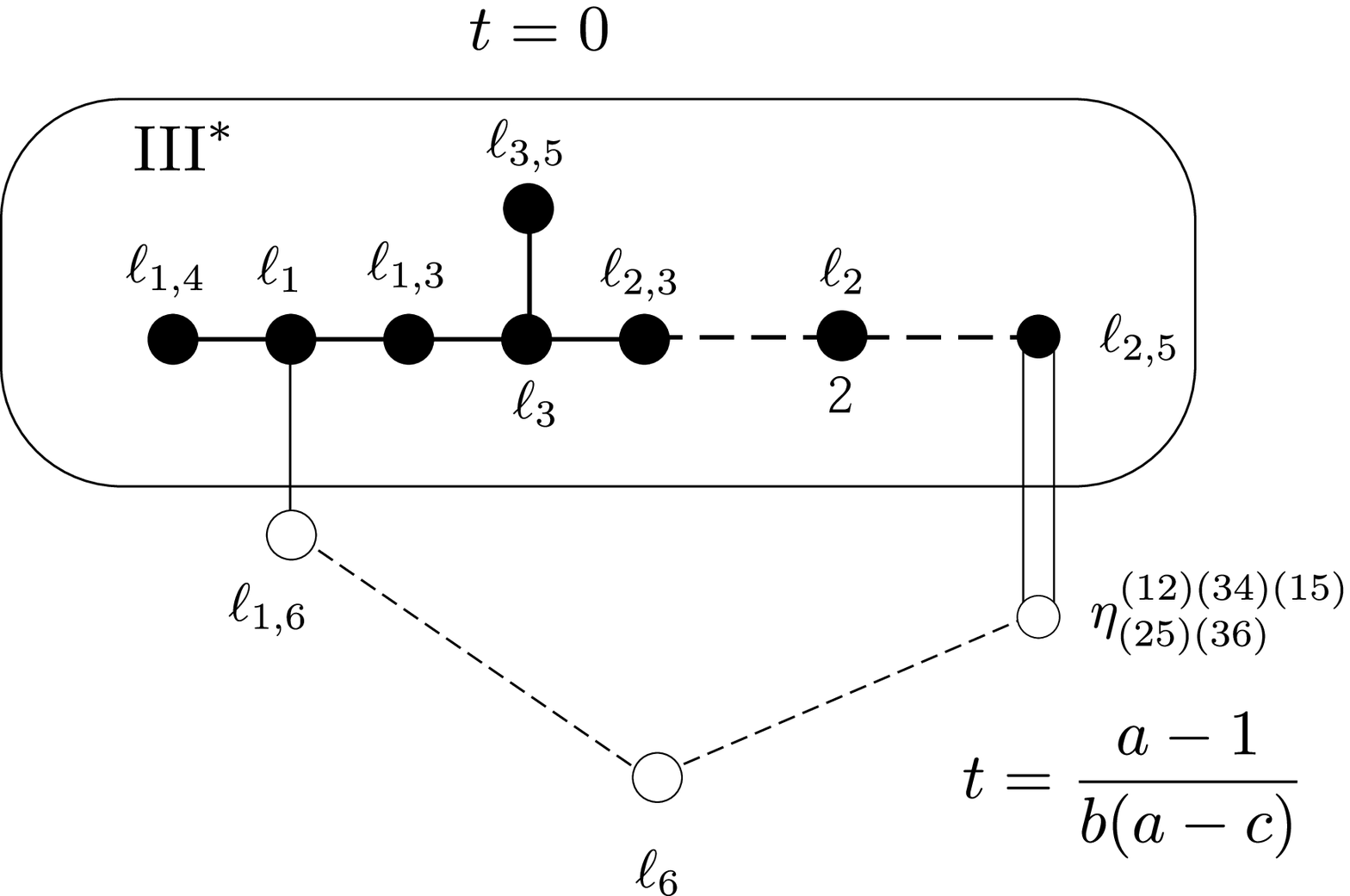}}\\ \hline
        \multicolumn{2}{|c|}{\includegraphics[width=9cm]{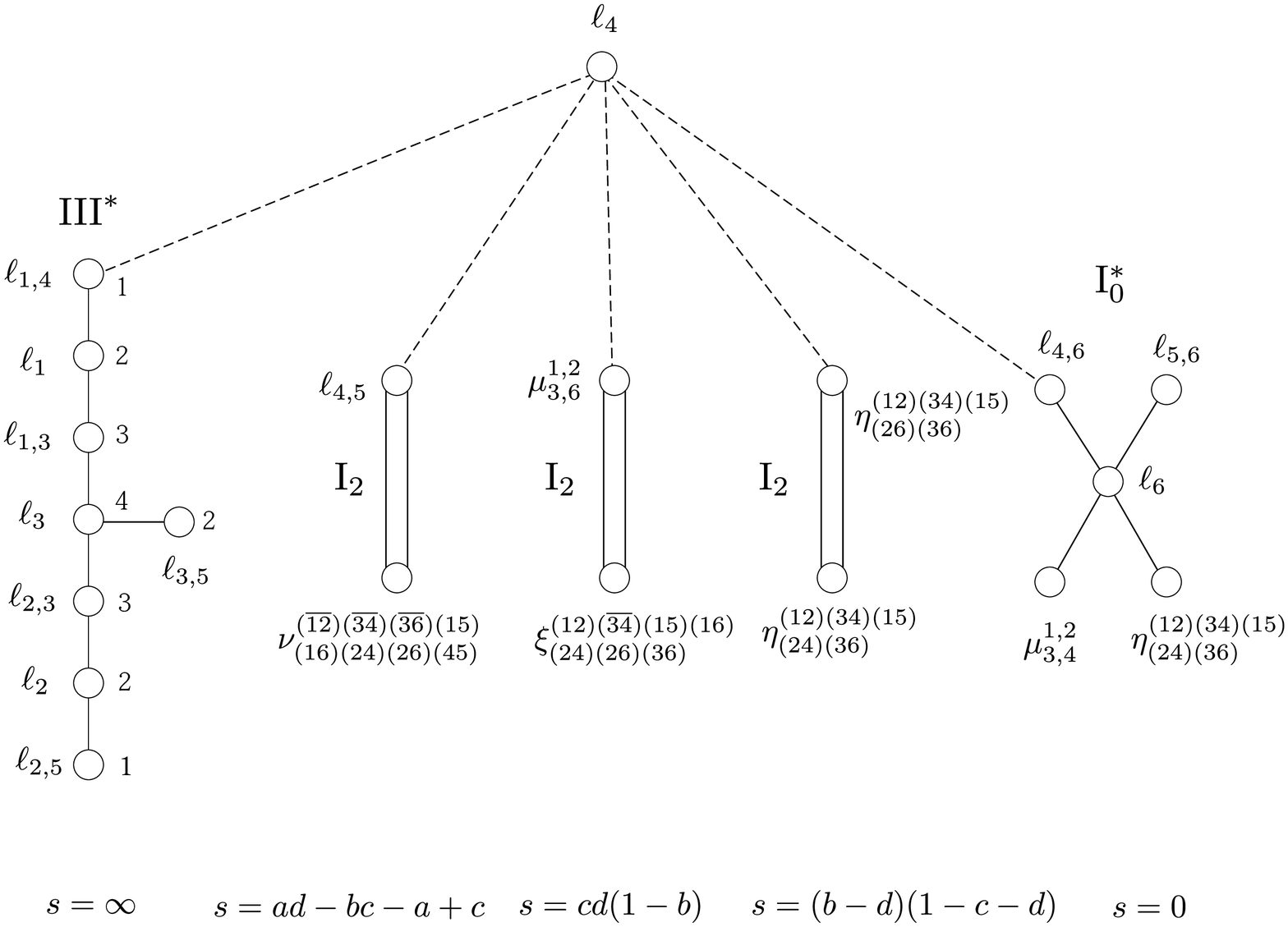}}\\ \hline
    \end{tabular}
    \vspace{.1in}
    \caption{Class 2.3}\label{tb:23}   
\end{table}

\begin{table}[htbp]
    \centering
    \begin{tabular}{|c|c|}\hline
        \multicolumn{2}{|c|}{Class 2.4} \\ \hline \hline
        Method & 2-neighbor step from the class 2.5 \\ \hline
        Elliptic parameter & {\SMALL
        $s = \dfrac{(ad-bc)(x-t)+(a-1)(ad-bc)(c+d-1) t^2}
        {t^2 \left( (a-1)(ad-bc)(c+d-1) t -ad+bc-b+d \right)}$} \\ \hline
        \multicolumn{2}{|c|}{\includegraphics[width=9cm]{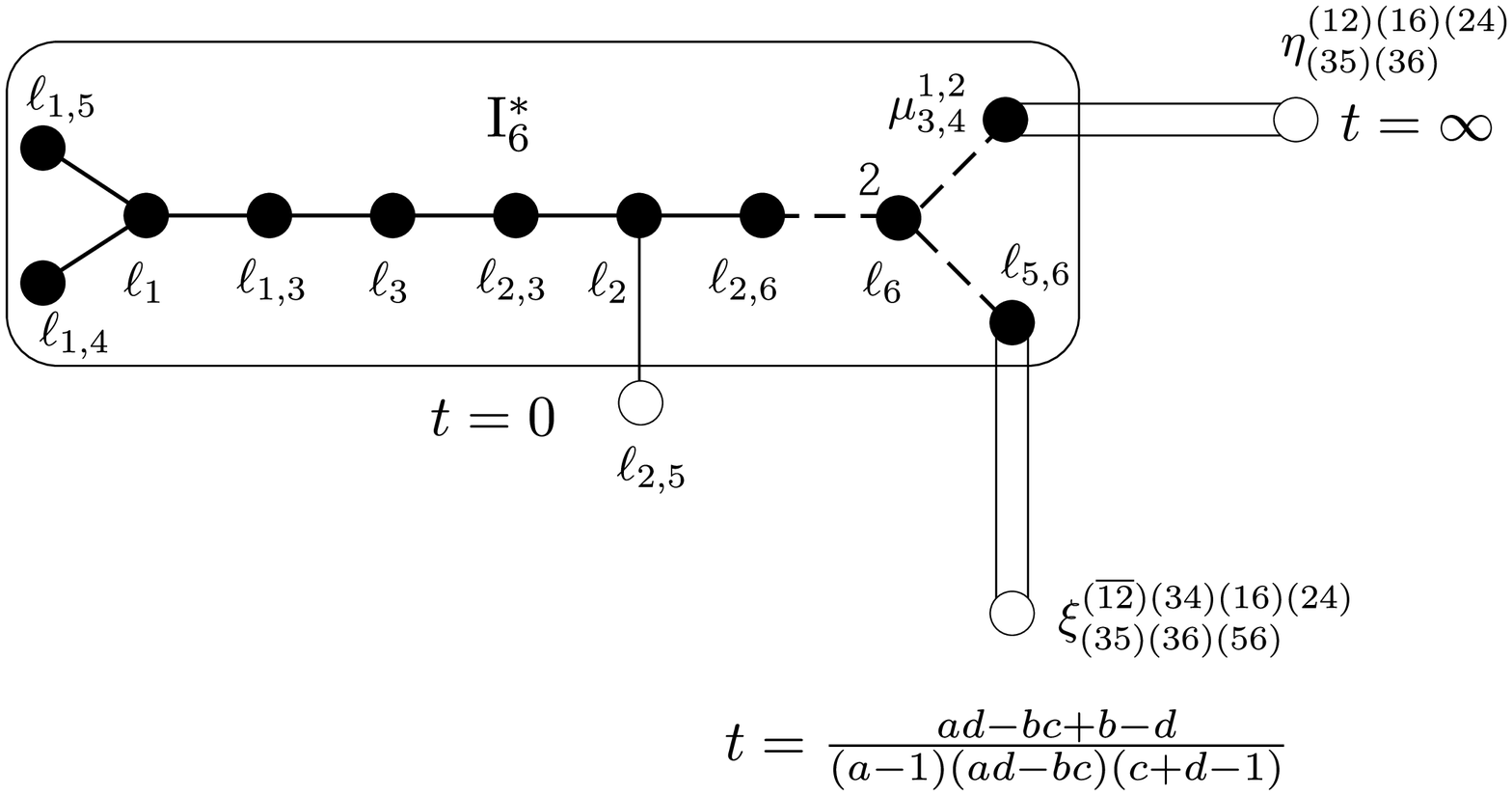}} \\ \hline
        \multicolumn{2}{|c|}{\includegraphics[width=9cm]{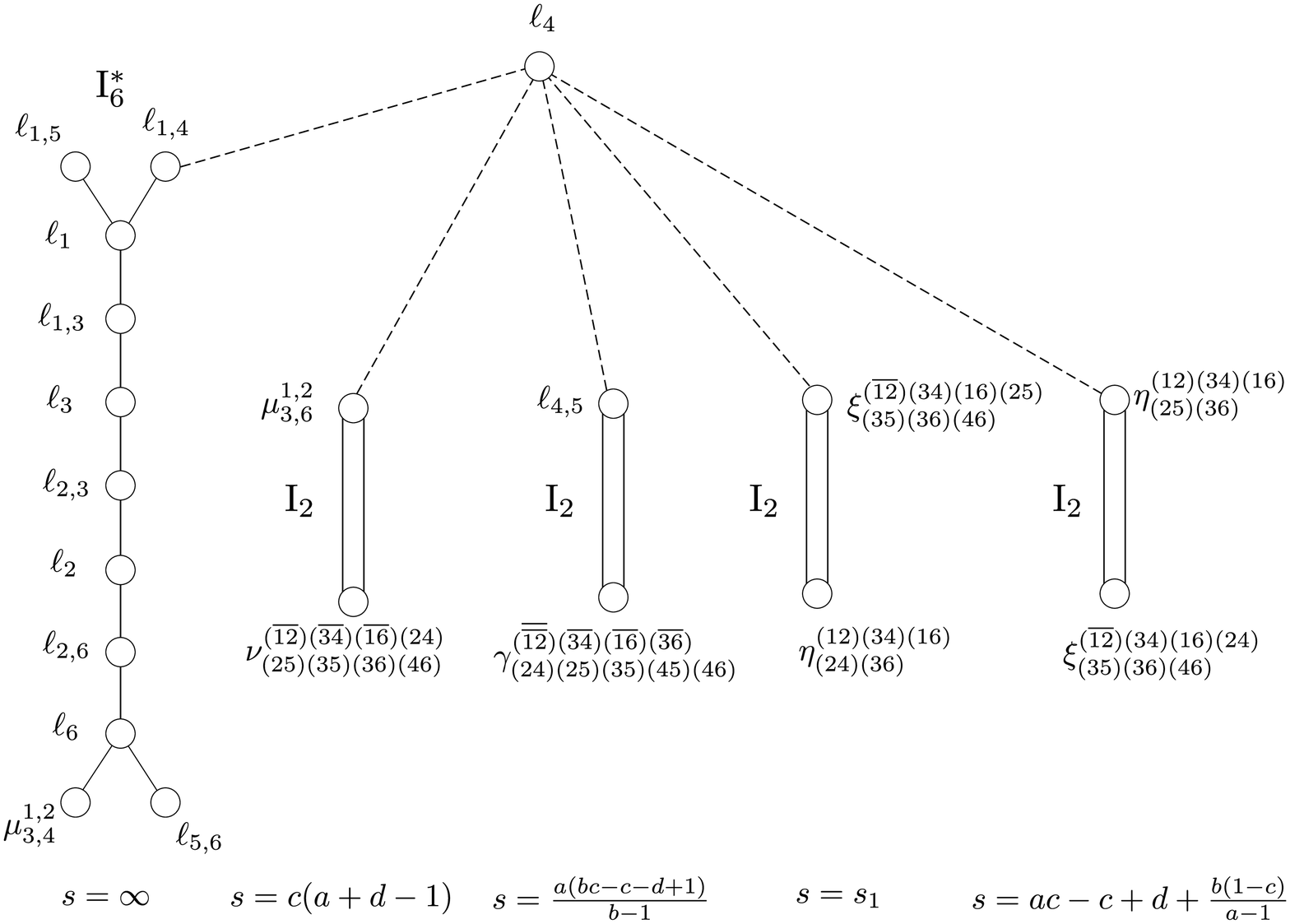}} \\ \hline
        \multicolumn{2}{|c|}{\small $s_1 = \dfrac{(c^2+2cd-2c-d+1)a+(d-1)(bc-b+d)}{c+d-1}$} \\ \hline
    \end{tabular}
    \vspace{.1in}
    \caption{Class 2.4}\label{tb:24}    
\end{table}

\begin{table}[htbp]
    \centering
    \normalsize
    \begin{tabular}{|c|c|} \hline
       \multicolumn{2}{|c|}{ Class 2.5} \\ \hline \hline
       Method & Classical \\ \hline
       Elliptic parameter & {\small $t=\dfrac {-u v (u-1) }
         {(a+c+d-ac-1) \, u \, v +acu^2 -(a+c) \, u -dv+1}$}\\ \hline
         \multicolumn{2}{|c|}{\includegraphics[width=9cm]{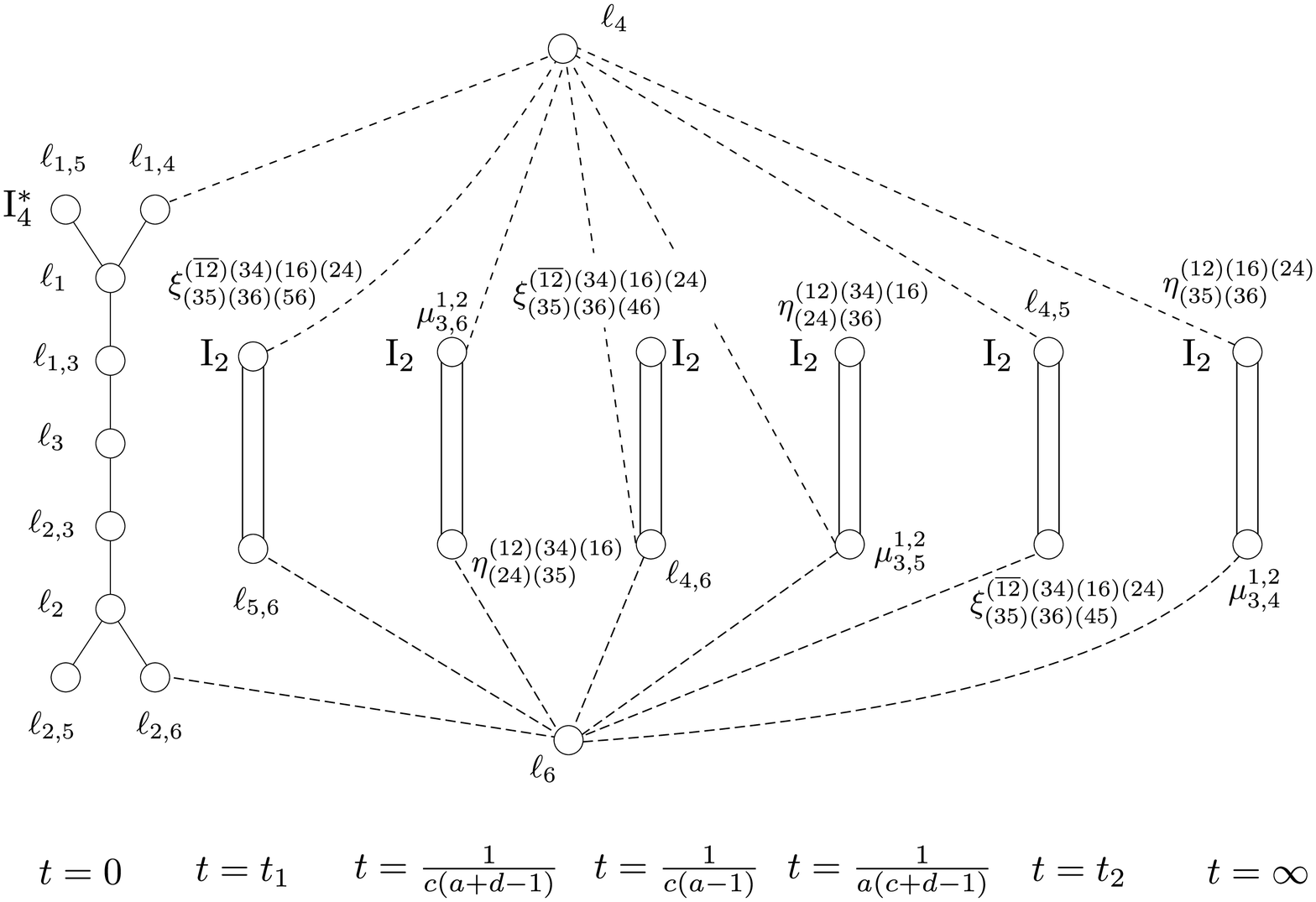} } \\
         \multicolumn{2}{|c|}{\SMALL $t_1 =\dfrac{ad-bc+b-d}{(a-1)(ad-bc)(c+d-1)}, \;
            t_2 = \dfrac{b-1}{a(bc-c-d+1)}$ }\\ \hline
         \multicolumn{2}{|c|}{
           $
           \begin{aligned}
               y^2 = &x^3 + t \, \Big( ac(a-1)(c+d-1)(ac+bc-a-c-d+1) t^3\\
               & + (4ac^2-c^2+a-5ac-cd-a^2+2acd+2bcd-ad^2-3a^2c^2\\
               & +bc^2-2a^2cd+c+3abc-abcd+a^2d-2bc+4a^2c-2bc^2a) t^2\\
               & +(ad+bc+3ac-2a-2b-2c+d+1) t -1 \Big) x^2\\
               & -t^4 \, \Big( c(a+d-1) t-1 \Big) \Big( a(bc-c-d+1) t -b+1 \Big) \Big( (a-1)\\
               & (ad-bc)(c+d-1)t-ad+bc-b+d \Big) x
           \end{aligned}
           $}\\ \hline
         Zero section & $\ell_6$\\ \hline
         2-torsion section & $\ell_4 : (x,y)=(0,0)$\\ \hline
     \end{tabular}
    \vspace{.1in}
    \caption{Class 2.5}\label{tb:25}   
\end{table}

\begin{table}[htbp]
    \centering
    \begin{tabular}{|c|c|} \hline
        \multicolumn{2}{|c|}{ Class 2.6}\\ \hline \hline
        Method & 2-neighbor step from the class 2.7\\ \hline
        Elliptic parameter & {\small
          $s=\dfrac{x-at+a(2b+c-d-1)t^2-a(b-d)(b+c-1)t^3}{t^2(bt-dt-1)(bct-1)}$}\\ \hline
        \multicolumn{2}{|c|}{\includegraphics[width=9cm]{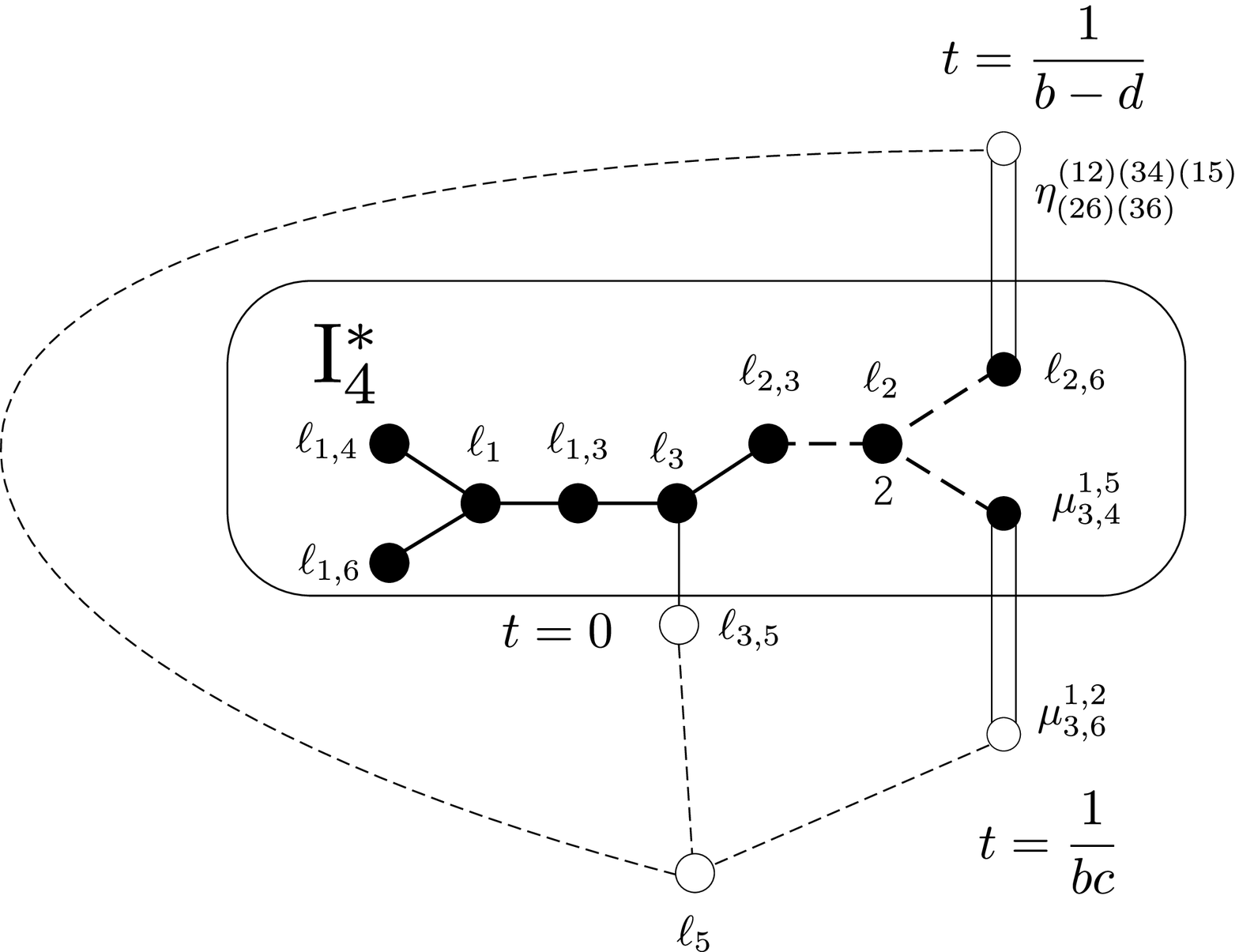}}\\ \hline
        \multicolumn{2}{|c|}{\includegraphics[width=9cm]{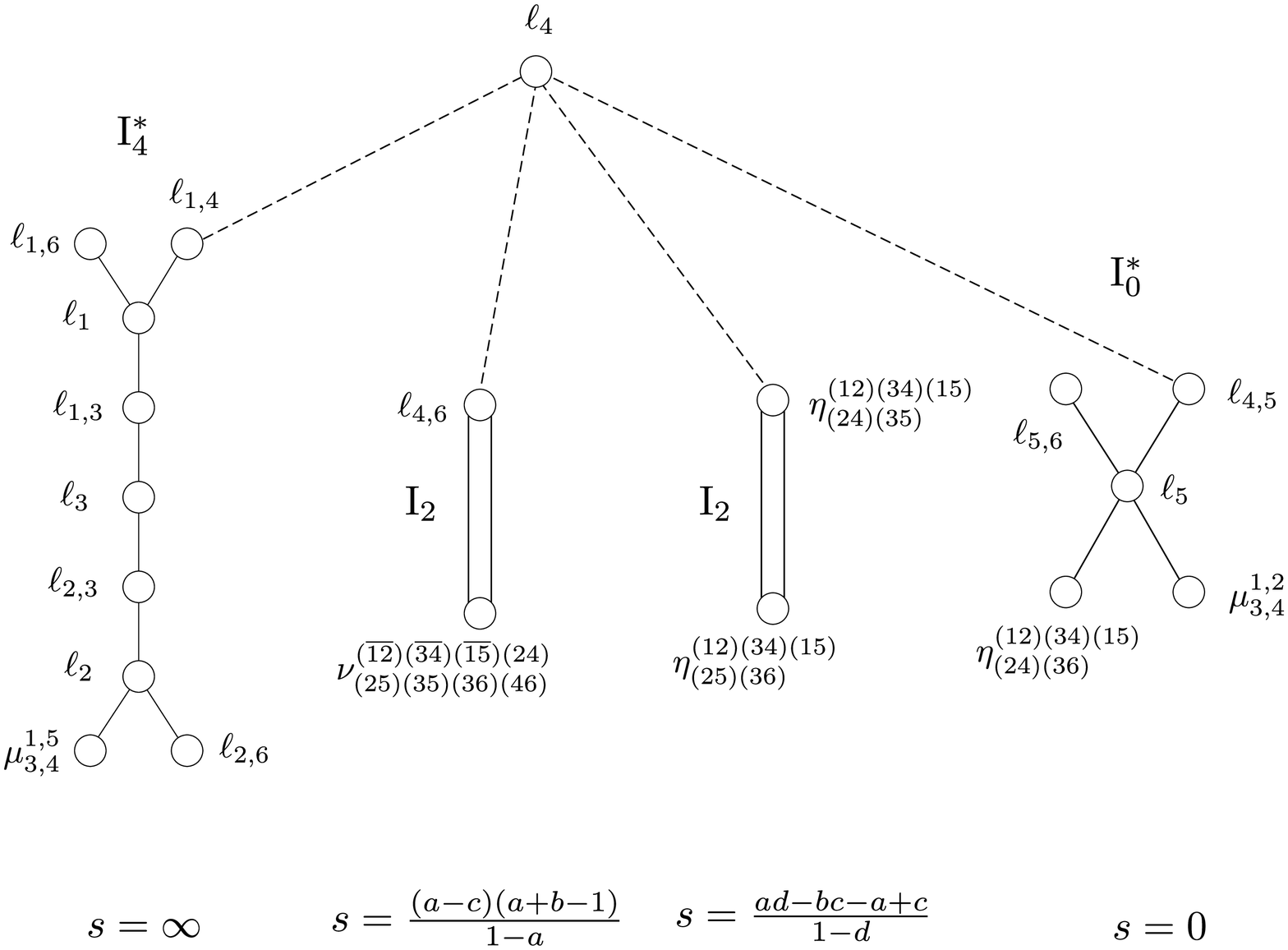}}\\ \hline
    \end{tabular}
    \vspace{.1in}
    \caption{Class 2.6}\label{tb:26}    
\end{table}

\begin{table}[htbp]
    \centering
    \begin{tabular}{|c|c|}\hline
        \multicolumn{2}{|c|}{Class 2.7}\\ \hline \hline
        Method & Classical \\ \hline
        Elliptic parameter & $\dfrac{uv}{cu+bv-1}$ \\ \hline
        \multicolumn{2}{|c|}{
          $
          \begin{aligned}
              y^2 = &\left( x-t(bt+ct-t-1)(abt-bct-a+1) \right)\\
              &  \left( x-a t (bt-dt-1)(bt+ct-t-1) \right)\\
              &  \left( x-t (1-ct)(bt-dt-1)(abt-bct-a+1) \right),
          \end{aligned}
          $
        }\\ \hline
        \multicolumn{2}{|c|}{\includegraphics[width=9cm]{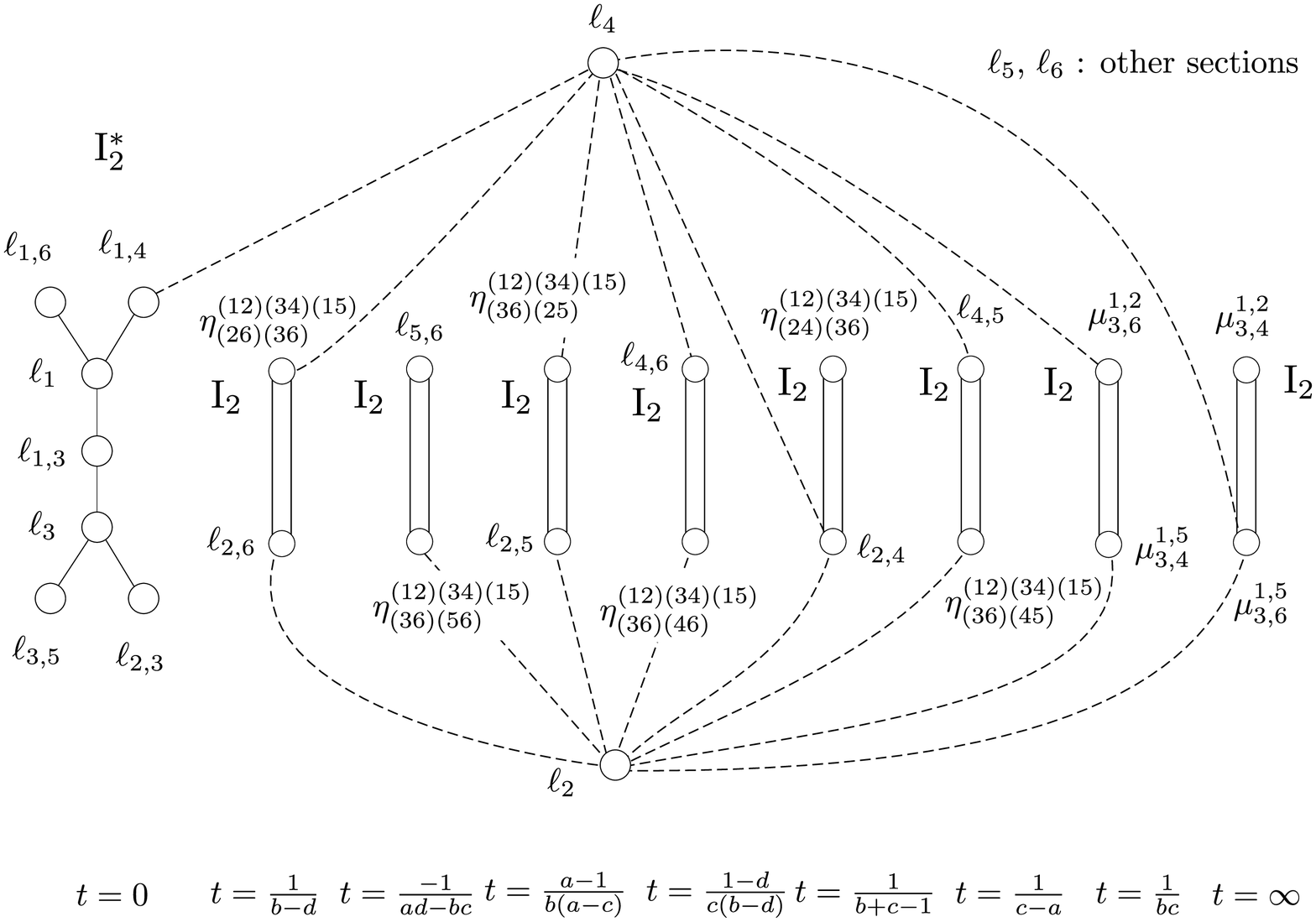}}\\ \hline
        Zero section & $\ell_2$ \\ \hline
        2-torsion section & {\small
        $
        \begin{aligned}
            \ell_4 & : (x,y) = \left( t(1-ct)(abt-bct-a+1)(bt-dt-1), 0 \right)\\
            \ell_5 & : (x,y) = \left( at(bt-dt-1)(bt+ct-t-1),0 \right)\\
            \ell_6 & : (x,y) = \left( t(abt-bct-a+1)(bt+ct-t-1), 0 \right)\\
        \end{aligned}
        $} \\ \hline
  \end{tabular}
  \vspace{.1in}
  \caption{Class 2.7}\label{tb:27}    
\end{table}

\begin{table}[htbp]
    \centering
    \begin{tabular}{|c|c|} \hline
       \multicolumn{2}{|c|}{ Class 2.8} \\ \hline \hline
       Method & Classical\\ \hline
       Elliptic parameter & \multicolumn{1}{|c|}{ $t=\dfrac{u (v-1)}{cu+dv-1}$ }\\ \hline
       \multicolumn{2}{|c|}{\normalsize
          $
          \begin{aligned}
              y^2 = & \, x^3 + t \, \big( (ad+bc+acd-bc^2-2cd)
              \, t^{2} \\
              & \qquad + (ad-2bc-2a+b+c-2d+1)t -b+1 \big) \, x^2\\
              &  + t^3 \, (cdt+d-1) (at-ct-1) (adt-bct+bt-dt-a-b+1)\, x
          \end{aligned}
          $} \\ \hline
        \multicolumn{2}{|c|}{\includegraphics[width=9cm]{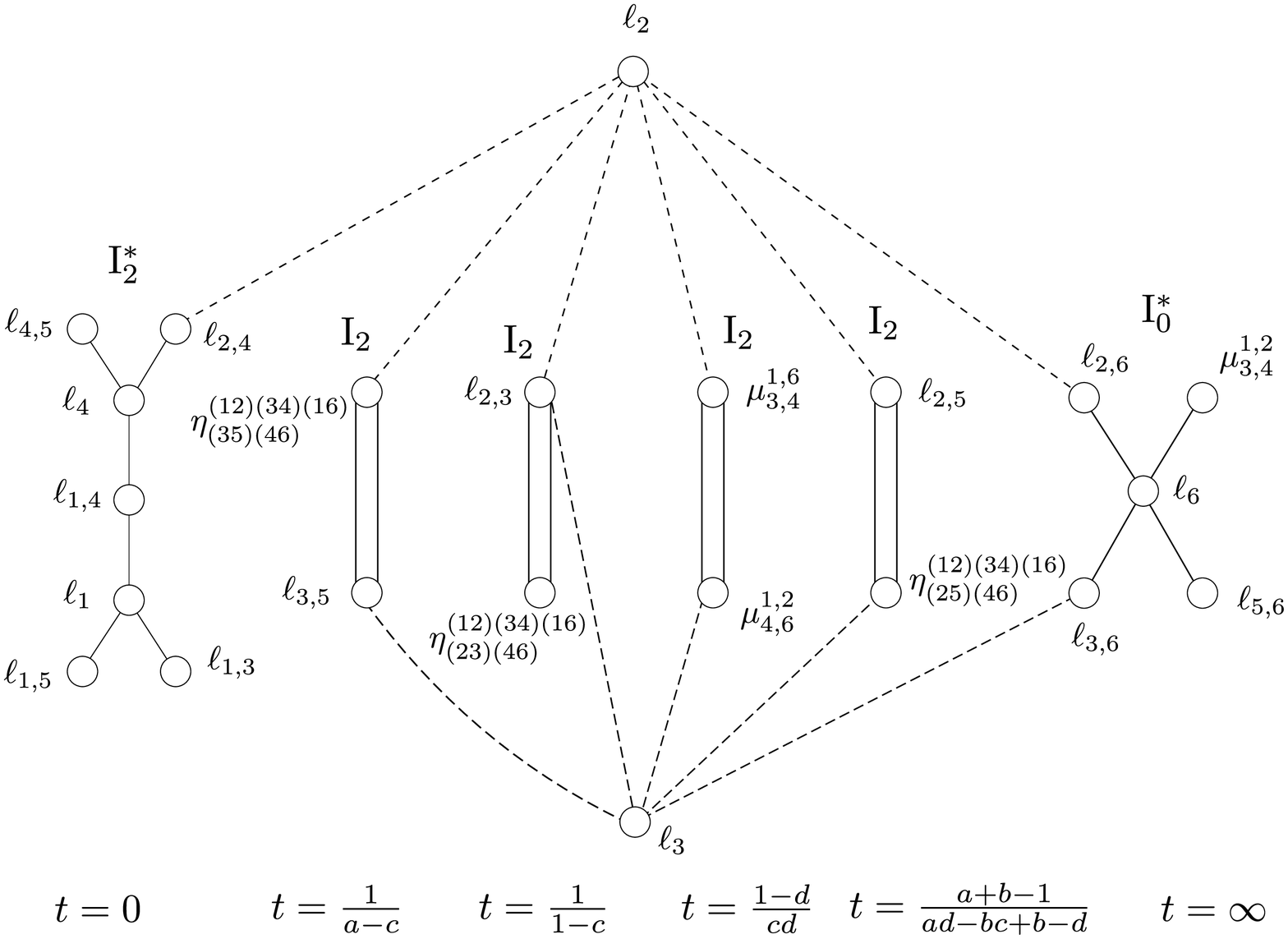}} \\ \hline
        Zero section & $\ell_2$ \\ \hline
        2-torsion section & $\ell_3 : (x,y)=(0,0)$\\ \hline
    \end{tabular}
    \vspace{.1in}
    \caption{Class 2.8}\label{tb:28}    
\end{table}

\begin{table}[htbp]
    \centering
    \begin{tabular}{|c|c|} \hline
       \multicolumn{2}{|c|}{ Class 2.9} \\ \hline \hline
       Method & Classical\\ \hline
       Elliptic parameter & \multicolumn{1}{|c|}{ $t=\dfrac{u (au+bv-1)}{(u-1)(v-1)(au-1)}$ }\\ \hline
        \multicolumn{2}{|c|}{\includegraphics[width=9cm]{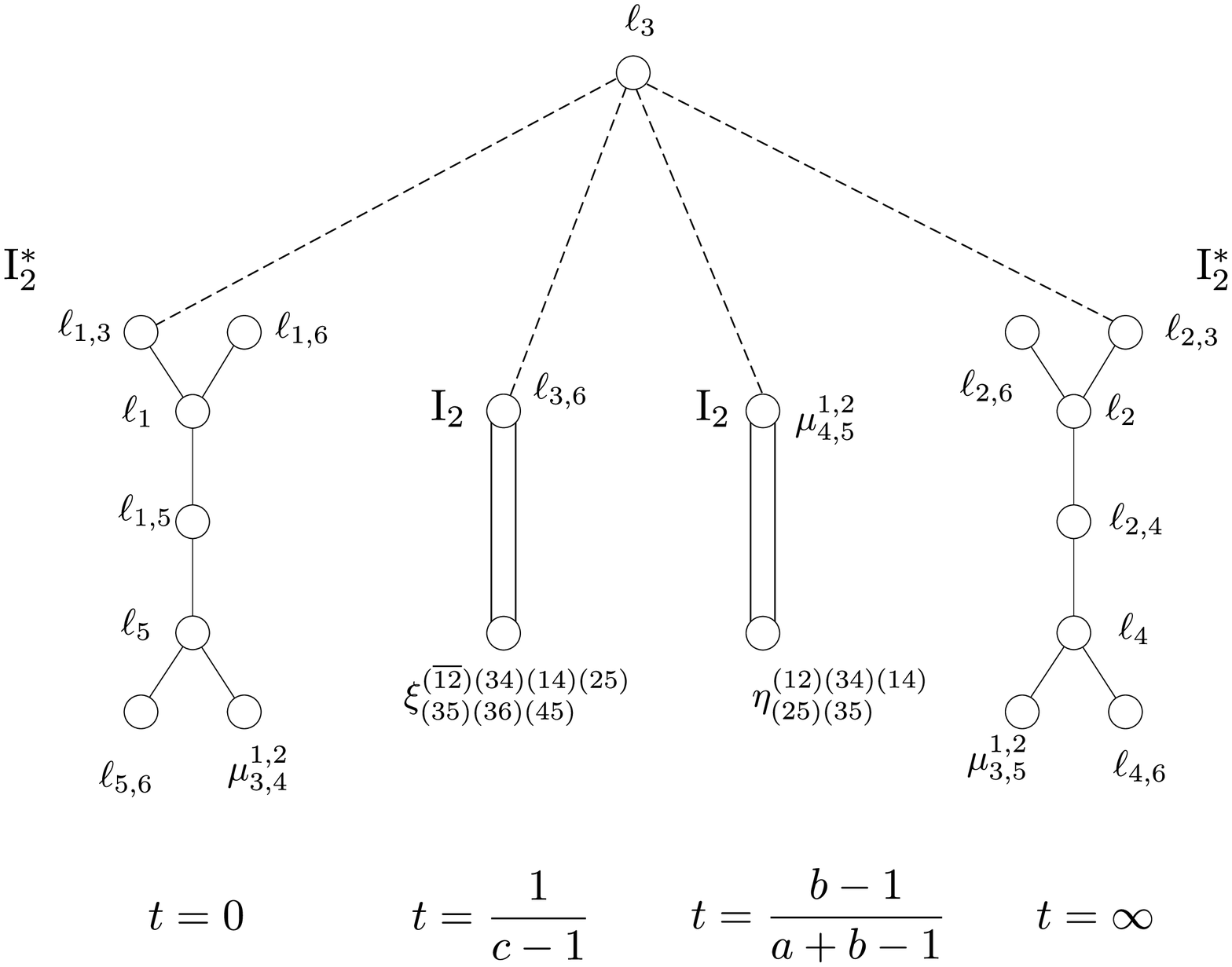}} \\ \hline
       \multicolumn{2}{|c|}{\normalsize
         $
         \begin{aligned}
             {y}^{2} &= {x}^3 +t (t+1) \big( (a-1)(c+d-1) t^3 -(2ac-2ad+2bc-b+3d-2) t^2\\
             & \qquad \quad +(ad-2bc+a+2b+c-3d+1) t^2 +b-d \big) {x}^{2}\\
             & \quad - t^3 (t+1)^2 (ct-t-1) \left( (a+b-1) t +b-1 \right) \big( (2ac+ad-a-c) t^2\\
                 & \qquad \quad -(ac-2ad+bc+a+c) t +ad-bc \big) x\\
             & \quad + ac{t}^{6} \left( t+1 \right) ^{3} \left(
               ct-t-1 \right) ^{2} \left( at+bt-t+b-1 \right) ^2
         \end{aligned}
         $}\\ \hline
         Zero section & $\ell_3$\\ \hline
     \end{tabular}
     \vspace{.1in}
     \caption{Class 2.9}\label{tb:29}    
 \end{table}

\begin{table}[htbp]
    \centering
    \begin{tabular}{|c|c|}\hline
        \multicolumn{2}{|c|}{ Class 2.10} \\ \hline \hline
        Method & 2-neighbor from the class 2.7\\ \hline
        Elliptic parameter & 
        $
        \begin{aligned}
            &s= \dfrac{(ad-bc) \,x -a(ad-bc+b-d) t +A_2 t^2}{t
              \left( (ad-bc) \, t+1\right) \left((b+c-1) \, t-1\right)},\\
            &A_2 = a(c+b-1)(ad-bc+b-d)\\
        \end{aligned}
        $\\ \hline
        \multicolumn{2}{|c|}{\includegraphics[width=9cm]{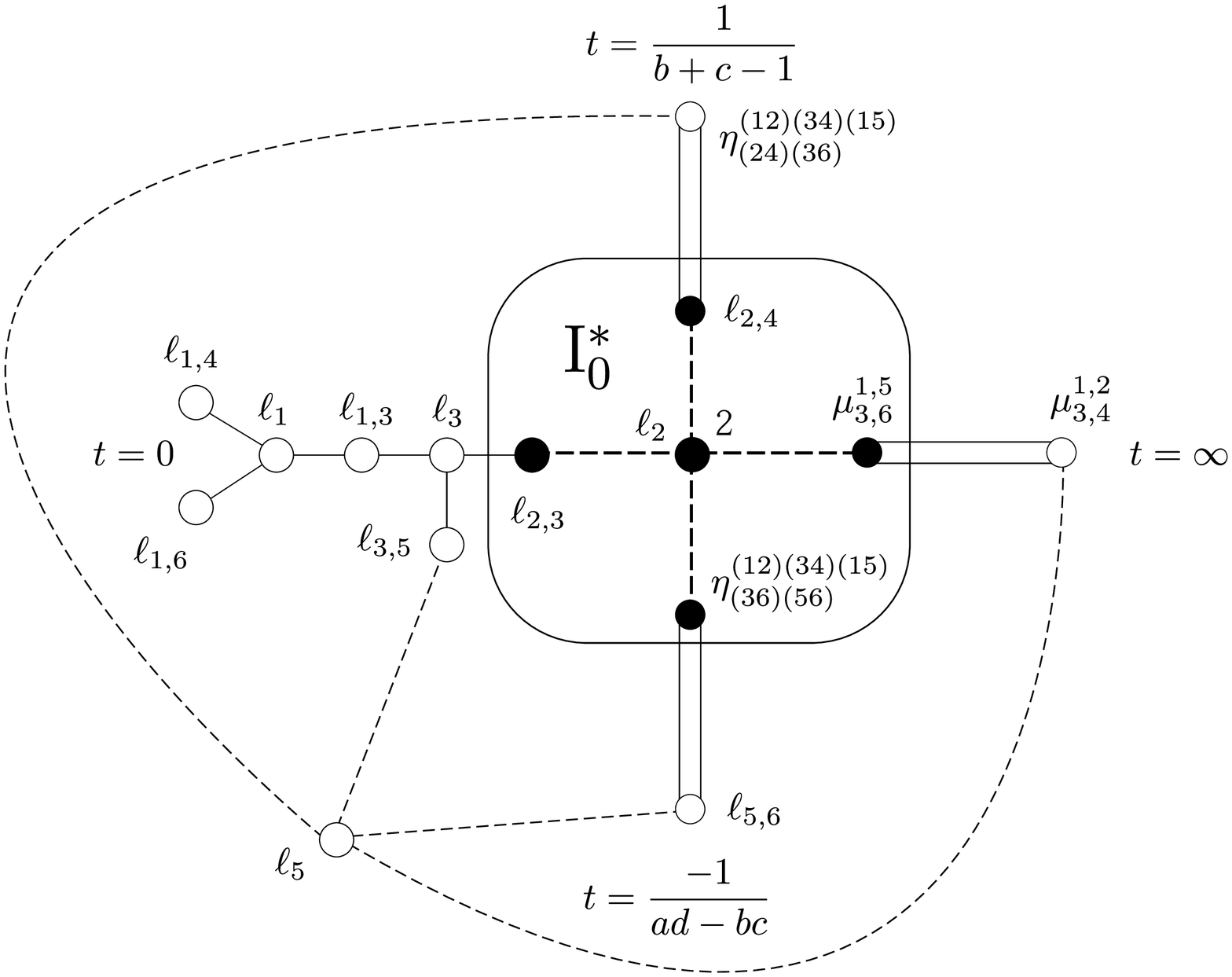}}\\ \hline
        \multicolumn{2}{|c|}{\includegraphics[width=9cm]{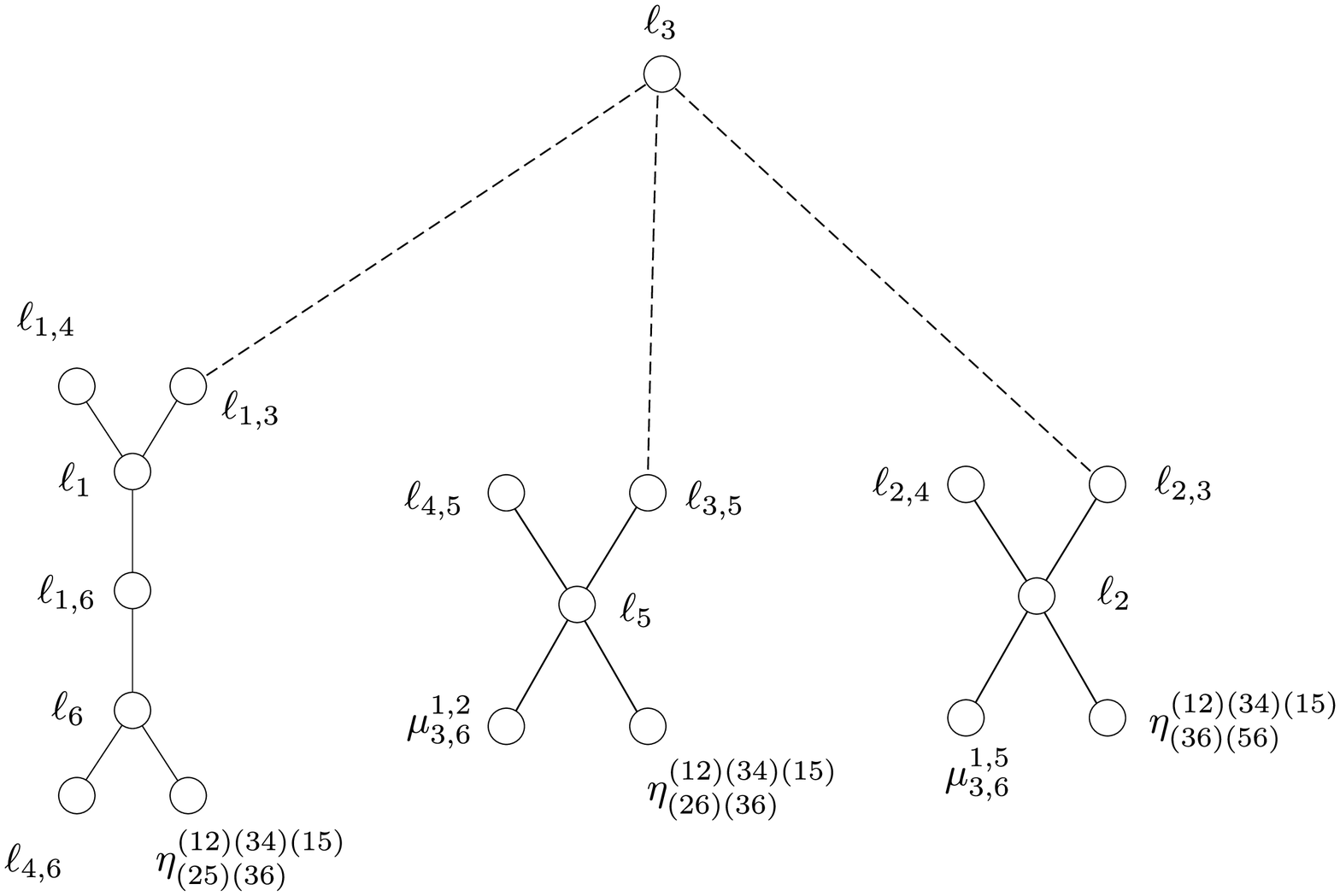}}\\ \hline
    \end{tabular}
    \vspace{.1in}
    \caption{Class 2.10}\label{tb:210}   
\end{table}

\begin{table}[htbp]
    \centering
    \begin{tabular}{|c|c|} \hline
       \multicolumn{2}{|c|}{ Class 2.11} \\ \hline \hline
       Method & Classical \\ \hline
       Elliptic parameter & $t=u$ \\ \hline
        \multicolumn{2}{|c|}{\includegraphics[width=9cm]{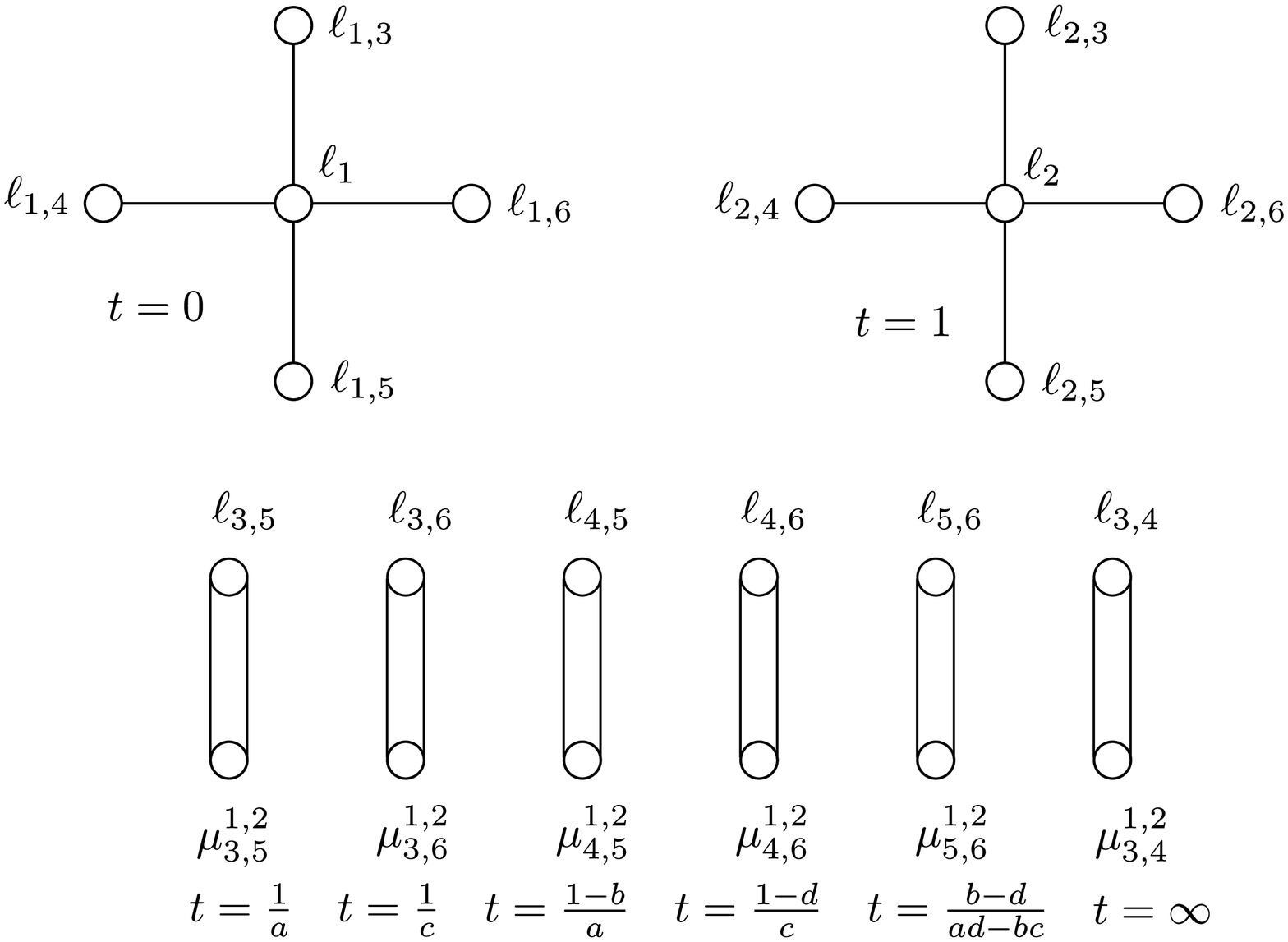}} \\ \hline
        \multicolumn{2}{|c|}{
          $
          \begin{aligned}
          y^2 \; =  & \; ( x- dt \, (t-1) (at-1) )\\
          & \; (x- bt \, (t-1) (ct-1) )\\
          & \; (x-t \, (t-1) (at-1) (ct-1) )
          \end{aligned}
          $} \\ \hline
          Zero section & $\ell_3$ \\ \hline
         2-torsion section &
         $
         \begin{aligned}
             \ell_4 :&  \; (x,y) = \left( t \, (1-t) (ct-1) (at-1) , \, 0 \right)\\
             \ell_5 :& \; (x,y) = \left( bt \, (t-1) (ct-1), \, 0 \right)\\
             \ell_6 :& \; (x,y) = \left( dt \,　(t-1) (at-1), \, 0 \right)\\
         \end{aligned}
         $\\ \hline
    \end{tabular}
    \vspace{.1in}
    \caption{Class 2.11}\label{tb:211}
\end{table}

\begin{table}[htbp]
    \centering
    \begin{tabular}{|c|c|} \hline
       \multicolumn{2}{|c|}{ Class 2.12} \\ \hline \hline
       Method & Classical\\ \hline
       Elliptic parameter & \multicolumn{1}{|c|}{ $t=\dfrac{u (bv+a-1)}{au+bv-1}$ }\\ \hline
         \multicolumn{2}{|c|}{\includegraphics[width=9cm]{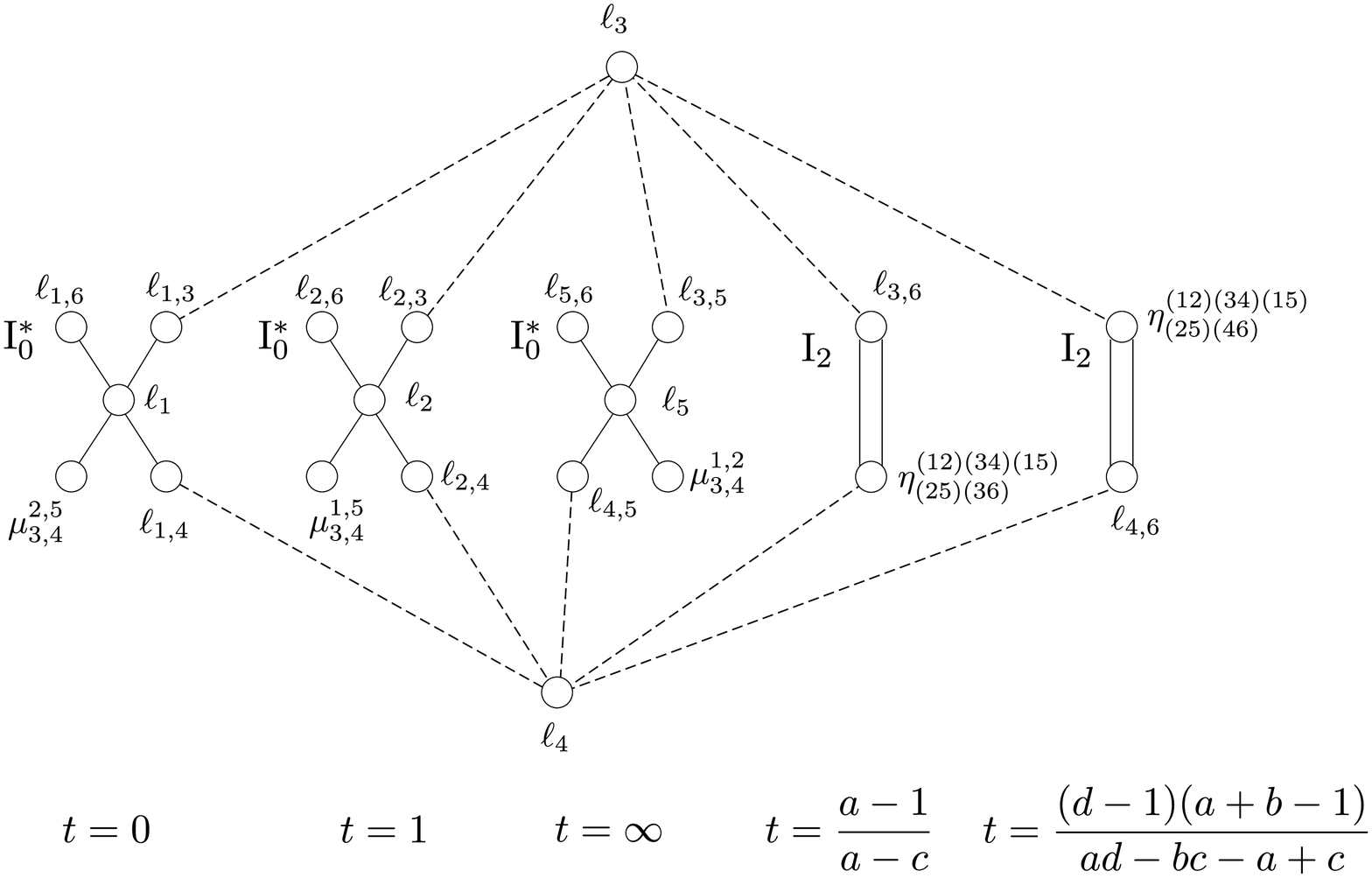}} \\ \hline
         \multicolumn{2}{|c|}{\small
           $
           \begin{aligned}
               & y^2 = \, x^3 +9 \, t \, (t-1) \, \left( (ad-bc-2a+2c) \, t -ad+2a+b+d-2 \right) x^2\\
               & -81 \, t^2 \, (t-1)^2 (at-ct-a+1) \big( (ad-bc-a+c) \, t -(d-1) (a+b-1) \big) x
           \end{aligned}
           $}\\ \hline
         Zero section & $\ell_3$\\ \hline
         2-torsion section & $\ell_4 : (x,y)=(0,0)$ \\ \hline
     \end{tabular}
    \vspace{.1in}
    \caption{Class 2.12}\label{tb:212}   
\end{table}

\clearpage

\section{ Class  2.7}\label{sec:class2.7}

An elliptic parameter for the class 2.7 is given by

\begin{equation}
    \label{eq:par-27}
    t=\frac{uv}{cu+bv-1}.
\end{equation}
It is easy to verify that the divisor of $t$ is given by
\begin{equation}\label{eq:div-27}
        (t) = \ell_{1,6}+\ell_{1,4}+2 \left( \ell_1 +\ell_{1,3}+\ell_3
        \right) +\ell_{2,3}+\ell_{3,5} - \left( \mu^{1,5}_{3,6} +
          \mu^{1,2}_{3,4} \right).
\end{equation}

\noindent
Then the fiber at $t=0$ is of type $\I_2^*$ and $t=\infty$ fiber is of
type $\I_2$. 

Eliminating the variable $v$ from (\ref{eq:DoubleCover}) and
(\ref{eq:par-27}), and making a simple coordinate change, we obtain an
equation of the form $y^2=$ (quartic polynomial). Choosing $\ell_2$
as the zero section of the group structure, we have the Weierstrass
equation for the Jacobian fibration

\begin{equation}
    \label{eq:27}
    \begin{aligned}
        y^2 = &\left( x-t(bt+ct-t-1)(abt-bct-a+1) \right)\\
        &  \left( x-a t (bt-dt-1)(bt+ct-t-1) \right)\\
        &  \left( x-t (1-ct)(bt-dt-1)(abt-bct-a+1) \right),
    \end{aligned}
\end{equation}
where the change of variables is given by
\begin{equation}
    \begin{aligned}
        x& = \frac {t (bt-dt-1) (bt+ct-t-1)(abt-bct-a+1)}{u-1}\\
        y &= \frac {t (bt-dt-1) (bt+ct-t-1) (abt-bct-a+1) (u-bt)^2 w}{u (cu-1) (u-1)^2}.\\
    \end{aligned}
\end{equation}

The remaining of divisors $\ell_4, \ell_5$ and $\ell_6$ are 2-torsion
sections. The correspondence between the divisors and the sections are
as follows.

\begin{equation}
    \begin{aligned}
        \ell_4 & \leftrightarrow (x,y) = \left( t(1-ct)(abt-bct-a+1)(bt-dt-1), 0 \right)\\
        \ell_5 & \leftrightarrow (x,y) = \left( at(bt-dt-1)(bt+ct-t-1),0 \right)\\
        \ell_6 & \leftrightarrow (x,y) = \left( t(abt-bct-a+1)(bt+ct-t-1), 0 \right)\\
    \end{aligned}
\end{equation}
We can check types and positions of other singular fibers by Tate
algorithm~\cite{Tate}. For example, the fiber at $t=\frac{1}{bc}$
is a singular fiber of type $\I_2$. Since we obtain
\begin{equation}
    \left( t- \frac{1}{bc} \right) = \left( \frac{(cu-1)(bv-1)}{(cu+bc-1)cb} \right) =
    \mu^{1,2}_{3,6}+\mu^{1,5}_{3,4} - \left(\mu^{1,2}_{3,4}+\mu^{1,5}_{3,6} \right)
\end{equation}
by (\ref{eq:par-27}) and (\ref{eq:divisors}), we see that the fiber at
$t=\frac{1}{bc}$ consists of $\mu^{1,2}_{3,6}$ and $\mu^{1,5}_{3,4}$.

Similarly, we can determine the other singular fibers. In the
following Figure~\ref{fig:2.7}, we show the complete configuration of
singular fibers for the class 2.7.

\begin{figure}[hbtp]
    \centering
    \includegraphics[width=9cm]{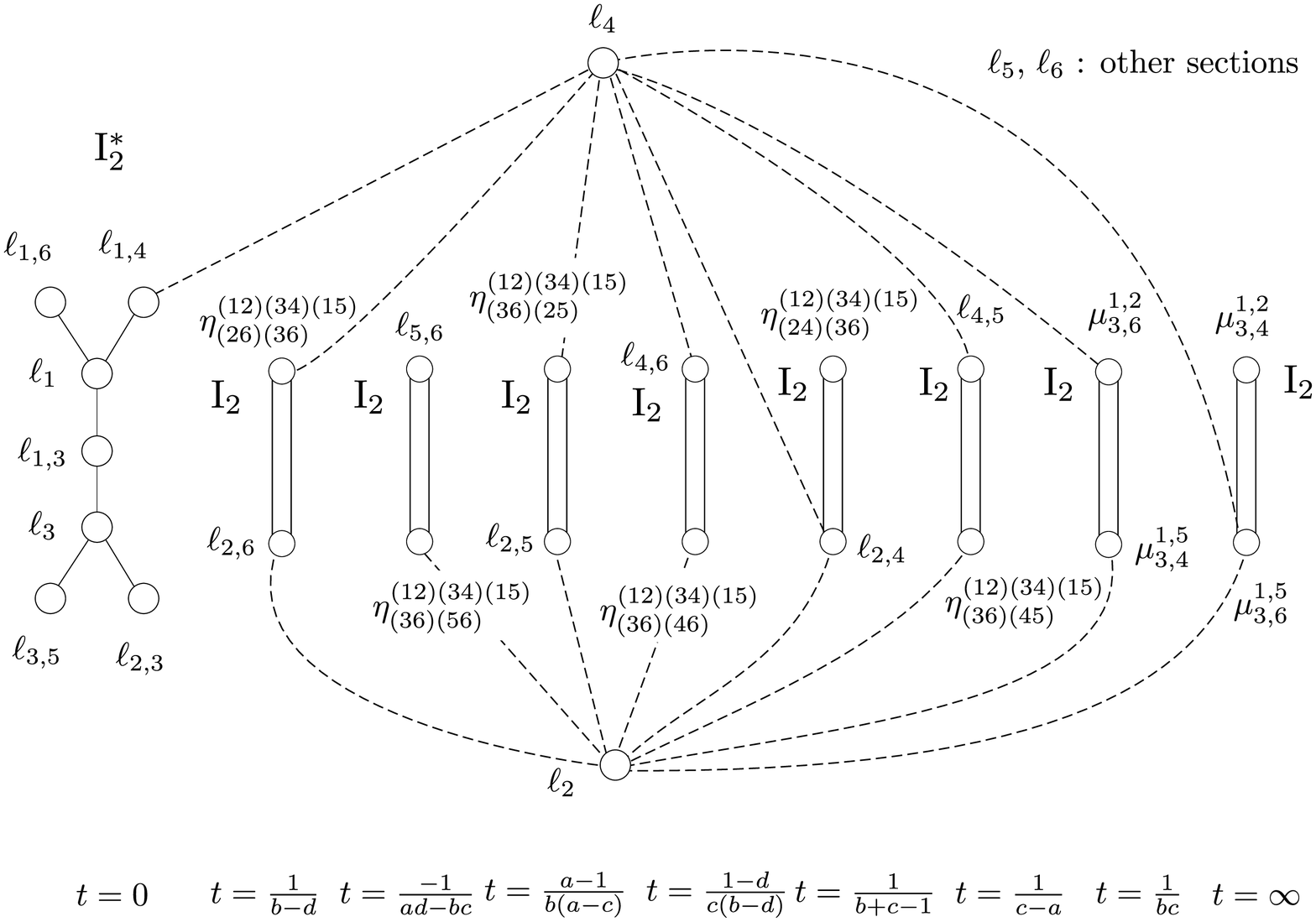}
    \caption{configuration of singular fibers for the class 2.7}
    \label{fig:2.7}
\end{figure}

\section{2-neighbor step }\label{sec:2-neighbor}

In this section, we explain ``2-neighbor step''. The following
description is based on \cite{Kumar}.

Let $X$ be a $K3$ surface over a field $k$ with an Jacobian fibration
over $\mathbb{P}^1$ with a zero section $O$, which defines an elliptic
curve $E$ over $k(t)$. Let $P, Q$ be other sections. Let $F$ be the
class of a fiber. We call $F$ an {\it elliptic divisor} of this
fibration.  Then for an effective divisor $F'=mO+nP+kQ+G$, where $G$
is an effective vertical divisor, we would like to compute the global
sections of $\mathcal{O}_X(F')$. Denote the space of the global
sections of $\mathcal{O}_X(F')$ by $L_X(F')$. Let $\{ s_1, \ldots, s_r
\}$ be a basis of $L_E(mO+nP+kQ)$. Then for any $f \in L_X(F')$, there
exist $b_i(t) \in k(t)$ such that
\begin{equation}
    f = b_1(t)s_1 + \cdots + b_r(t) s_r.
\end{equation}

We suppose that there exists an elliptic divisor $F'$ on $X$ with $F'
\cdot F=2$. Then decomposing $F'$ into horizontal and vertical
components $F'=F'_h+F'_v$, we see that $F'_h \cdot F + F'_v \cdot F =
F' \cdot F =2$. Since $F'_v \cdot F=0$, we have $F'_h \cdot F
=2$. Therefore the possibilities are $F'_h=2P$ or $F'_h=P+Q$ for
sections $P, Q$. By a translation, we can take the first one to be
$2O$, and the second to be $O+T$ with a 2-torsion section $T$ or $O+P$
with a non 2-torsion section $P$ depending on whether the class of the
section $[P-Q]$ is 2-torsion or not.

We suppose that a Weierstrass equation of $E/k(t)$ is given by
\begin{equation}
    \label{eq:Weierstrass}
    y^2=x^3+a_2(t)x^2+a_4(t)x+a_6(t),
\end{equation}
with $a_i(t) \in k[t]$ of degree at most $2i$.

First, we consider the case $D=2O$. Then $1$ and $x$ form a basis for
$L_E(D)$. Therefore, for a new elliptic divisor $F'=2O+G$ with $G$
effective and vertical, we obtain two elements $1$ and
$A(t)+B(t)x$ of $L_X(F')$ for some fixed $A(t), B(t) \in
k(t)$. The ratio of the two global sections  gives the new
elliptic parameter $s=A(t)+B(t)x$. Therefore we set
\begin{equation}
    x=\frac{s-A(t)}{B(t)}
\end{equation}
and substitute this into the Weierstrass equation, to get the form of
\begin{equation}
    y^2=g(t,s).
\end{equation}
Since the generic fiber of the fibration over $\mathbb{P}^1_s$ is a
curve of genus $1$, after absorbing square factors into $y^2$, 
$g$ must be a polynomial of degree $3$ or $4$ in $t$.

Next, we consider the case $D=O+P$ where $P=(x_0,y_0)$ is not a
2-torsion section. Then $1$ and $\frac{y+y_0}{x-x_0}$ form a basis of
$L_E(D)$. Therefore we obtain a new elliptic parameter
\begin{equation}\label{eq:s}
    s=A(t)+B(t)\frac{y+y_0}{x-x_0}.
\end{equation}
Solving (\ref{eq:s}), we get
\begin{equation}
    y=\frac{(s-A(t))(x-x_0)}{B(t)}-y_0.
\end{equation}
Substituting into the Weierstrass equation (\ref{eq:Weierstrass}) we get
\begin{equation}\label{eq:PO}
\left(\frac{(s-A(t)(x-x_0)}{B(t)}-y_0 \right)^2 = x^3+a_2(t)x^2+a_4(t)x+a_6.    
\end{equation}
Since $x_0$ and $y_0$ satisfy the Weierstrass equation
(\ref{eq:Weierstrass}), the difference of the left and right hand
sides of this equation can be divided by $(x-x_0)$. Therefore we get
an equation $g(x,t,s)=0$ that is quadratic in $x$. By completing the
square, we obtain an equation of the form
\begin{equation}
    x^2=h(t,s)
\end{equation}
and, after absorbing square factors into $x$, we have that $h$ is
cubic or quartic in $t$.

Finally, we consider the case $D=O+T$ where $T$ is a 2-torsion
section. In this case, we may assume that $T=(0,0)$ and $a_6=0$ by a
translation. Then $1$ and $\frac{y}{x}$ form a basis of $L_E(D)$. 
Setting a new elliptic parameter $s=A(t)+B(t)\frac{y}{x}$, we obtain
\begin{equation}
    y=\frac{(s-A(t))x}{B(t)}.
\end{equation}
Substituting this into the Weierstrass equation
(\ref{eq:Weierstrass}), we have
\begin{equation}
    \left(\frac{(s-A(t))x}{B(t)}\right)^2=x^3+a_2(t)x^2+a_4(t)x.
\end{equation}
Dividing both sides by $x$, we obtain a quadratic equation, and we can
proceed as in the previous case.

\section{Class 2.10}\label{sec:class2.10}

To obtain the Weierstrass equation for the class 2.10, we use a
2-neighbor step from the class 2.7. We compute explicitly the elements
of $L_X(F')$ where

\begin{equation}
    F'=2\ell_2+\ell_{2,3}+\ell_{2,4}+\mu^{1,5}_{3,6}+\eta^{(12)(34)(15)}_{(36)(56)}
\end{equation}

\begin{figure}[htbp]
    \centering
    \label{fig:27-210}
    \includegraphics[width=9cm]{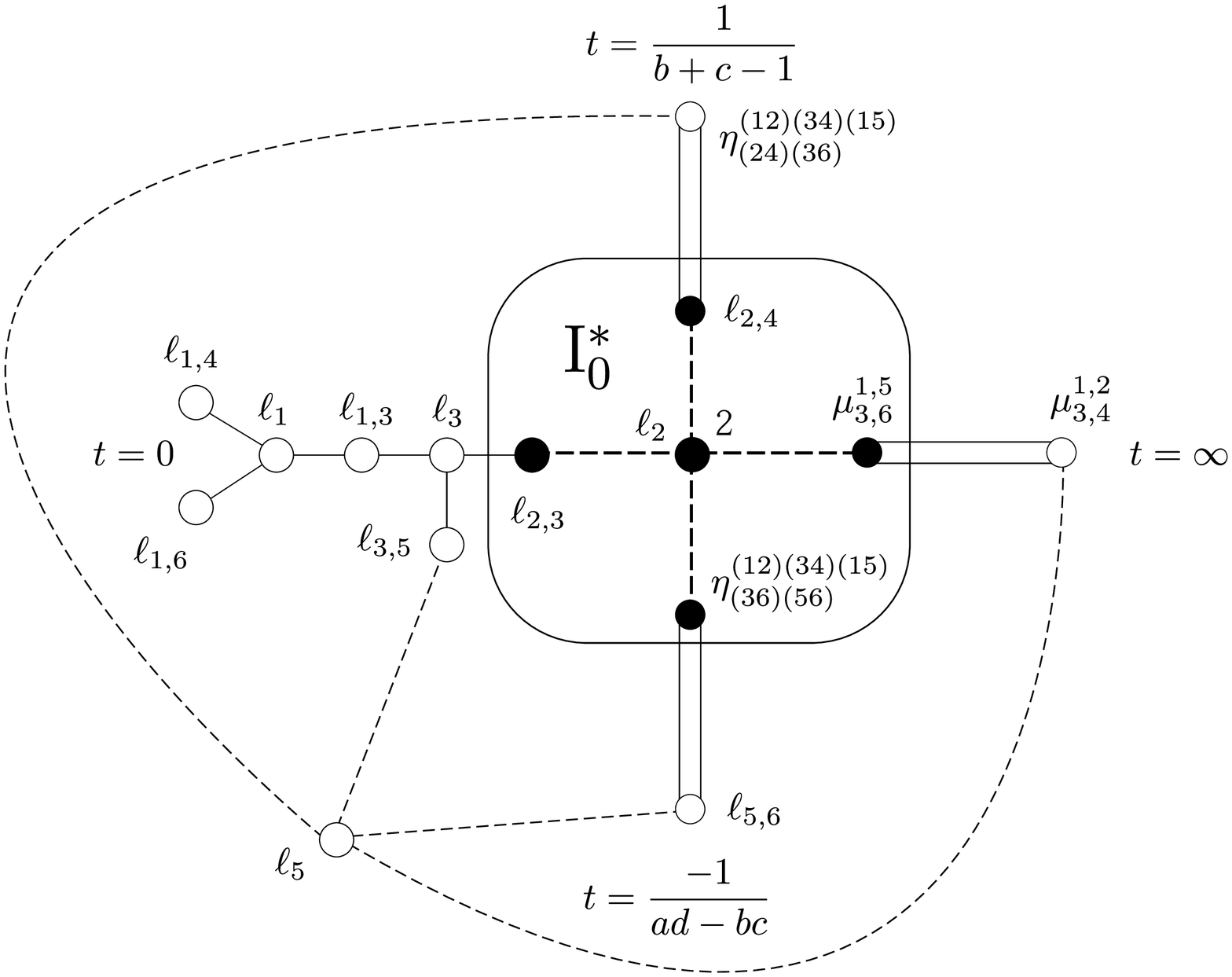}
    \caption{2-neighbor step from class 2.7 to class 2.10}
\end{figure}

\noindent
is the class of the fiber of type $\I_0^*$ we are considering. The
linear space $L(F')$ is 2-dimensional, and the ratio of two linearly
independent elements will be an elliptic parameter for $X$. Thus, we
may find a non-constant rational function on $X$ belonging to $L(F')$,
for which $1$ is an element of $L(F')$. Then it will be an elliptic
parameter of fibration 2.10. Let $s \in L(F')$ be a
non-constant. Notice that $s$ has a pole of order $2$ along $\ell_2$,
which is the zero section of fibration 2.7. Also, it has a simple pole
along $\ell_{2,4}$, the identity component of the fiber at
$t=\frac{1}{b+c-1}$ , a simple pole along $\mu^{1,5}_{3,6}$, the
identity component of the fiber at $t=\infty$, and a simple pole along
$\eta^{(12)(34)(15)}_{(36)(56)}$, the identity component of the fiber
at $t=\frac{-1}{ad-bc}$. Therefore we can put

\begin{equation}
    s = \frac{x+A_0+A_1 t+A_2 t^2+A_3 t^3+A_4 t^4}{t( t+\frac{1}{ad-bc}
      )(t- \frac{1}{b+c-1})}.
\end{equation}

We may subtract a term $A_3 t \left( t+\frac{1}{ad-bc} \right) \left(
  t- \frac{1}{b+c-1} \right)$ from the numerator, since $1$ is an
element of $L(F')$. Thus, assume $A_3=0$. To obtain other coefficients
the $A_i$, we look at the order of vanishing along the non-identity
components of fibers at $t=0, \infty, \frac{-1}{ad-bc},
\frac{1}{b+c-1}$. For example, we look at the fiber at $t=0$. The
rational function $s$ does not have any pole along $\ell_{3,5}$, which
intersects with the section $\ell_5$ of the fibration for the class
2.7 at $t=0$. Hence $s$ has no pole at $t=0$ and $x
=-at(bt-dt+1)(bt+ct-t+1)$, which corresponds to the section $\ell_5$,
and that gives us $A_0$. Similarly, other fibers give remaining
coefficients. After some calculation, we get a new elliptic parameter
\begin{equation}
    \label{eq:s-210}
        s= \frac{(ad-bc) \,x -a(ad-bc+b-d) \, t+a(c+b-1)(ad-bc+b-d) \, t^2}{t
          \left( (ad-bc) \, t+1\right) \left((b+c-1) \, t-1\right)}.
\end{equation}
Solving for $x$ in terms of $s$ and substituting into right hand side
of (\ref{eq:27}), and dividing out some square factors, which can
absorb into $y^2$, we obtained an equation for fibration 2.10 in the
form $y^2=f(t)$, where $f$ is a quartic in $t$ with coefficients in
$k(a,b,c,d,s)$. In fact, since the quartic $f$ factors into
$t$ times cubic factor in $t$, the new equation converts to a
Weierstrass form by a standard algorithm. Finally, we get a
Weierstrass equation of the fibration for the class 2.7.
\begin{equation}
    \label{eq:210}
    \begin{aligned}
        Y^2 = X^3 &+ (ad-bc)(s+ad-ab)(s-ab-bc) \left( (b+c-ad+bc-1) \,
          s \right.\\
        & \quad \left. +b^2c^2+b^2c-2abc-a^2d^2-bcd+a^2d-ad-ab^2
          +ad^2 \right.\\
        & \quad \left. +bc^2+ab-ab^2c+a^2bd \right) \, X^2\\
        &-(s-ab+ad)^2(s-ab+bc)^3(ad-bc)^3
        \left(s^3-3b(a-c)s^2 \right.\\
        & \quad \left. +3b^2(a-c)^2s-b^2(a-c)^3 \right) \, X\\
        &+bc \, (a-c)(b-d)(ad-bc)^4(s-ab+ad)^3(s-ab+bc)^3.
    \end{aligned}
\end{equation}

The new variables $X, Y$ are related to the ones from fibration 2.7 by the following equation.

\begin{equation}
    \begin{aligned}  
        X &= \frac{(ad-bc)(ab-ad-s)(ab-bc-s)^2}{t}\\
        Y &= \frac{(ad-bc)^3(ab-ad-s)(ab-bc-s)^2 y}{t^2}.
      \end{aligned}
\end{equation}

The elliptic parameter $s$ for the class 2.10 is also given in terms
of $u,v$ by

\begin{equation}
    \begin{aligned}
        s &= \frac{g}{\left( u-1 \right) \left( audv-bvcu+cu+bv-1
          \right) \left( cu+bv-1 \right) },\\
    \end{aligned}
\end{equation}
where $g$ is given by
\begin{equation}
    \begin{aligned}
        g= & b(a-c)(b-d)(ad-bc)u^2 v^2 +ac^2(ad-bc+b-d) u^3\\
        &  +c(2ab^2-b^2c^2-abd-b^2c+bcd-ad^2+a^2d^2)u^2 v\\
        &  -b(b^2c^2+ad^2-a^2d^2+b^2c+a^2bd-bcd-ab^2-ab^2c)u v^2\\
        &  +c(bc^2+abc+2ad-2a^2d-2ab)u^2 -b^3 (a-c) v^2\\
        &  +(ad^2-2ab^2c+b^2c-2ab^2+abd+3b^2c^2-a^2d^2-bcd) uv\\
        &  +(ab-ad-2bc^2+abc+a^2d) u + 2b^2 (a-c) v -b(a-c), \\
    \end{aligned}
\end{equation}
which defines the rational plane curve through $P_{1,2}, P_{3,4}, P_{1,5},
P_{2,5}, P_{2,6}, P_{3,6}$ and $P_{5,6}$ with double points at
$P_{1,2},
P_{1,5}$ and $P_{3,6}$.

The reducible fibers are as follows:

\begin{center}
    \begin{tabular}{|l|c|c|}\hline
        Position & Reducible fiber & Type\\ \hline \hline
        $s=\infty$ & $2\ell_2 +\ell_{2,4}+ \mu^{1,5}_{3,6} + \eta^{(12)(34)(15)}_{(36)(56)}$ &$\I_0^*$ \\\hline
        $s=b(c-a)$ & $2\ell_5 +\ell_{4,5}+ \mu^{1,2}_{3,6} + \eta^{(12)(34)(15)}_{(26)(36)}$ &$\I_0^*$ \\ \hline
        $s=a(b-d)$ & $\ell_{1,4}+\ell_{1,3}+2(\ell_1+\ell_{1,3}+\ell_{6})+ \ell_{4,6}+
        \eta^{(12)(34)(15)}_{(25)(36)}$  & $\I_2^*$\\ \hline
    \end{tabular}\\
\end{center}

Then $\ell_3$ is the zero section and $\ell_4$ is a tri-section.

\begin{figure}[hbtp]
    \centering
    \includegraphics[width=9cm]{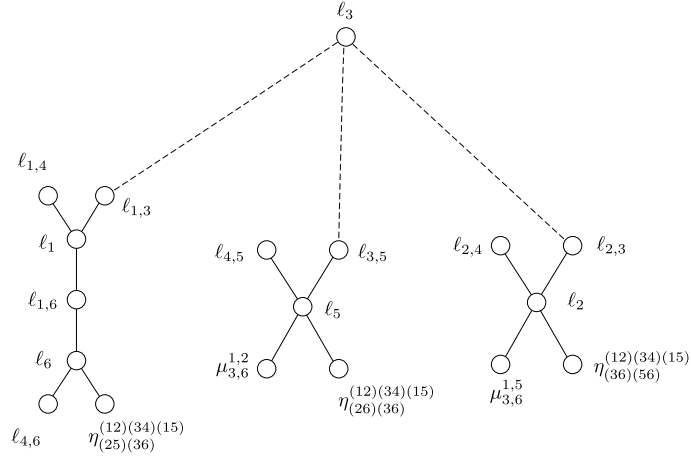}
    \caption{Configuration of singular fibers for the class 2.10}
    \label{fig:2.10}
\end{figure}

\section{Class 2.4}\label{sec:class2.4}

The Weierstrass equation for the class 2.4 is obtained from class 2.5
by using a 2-neighbor step as the following.

\begin{center}
    \includegraphics[width=7cm]{I4s_to_I6s.eps}
\end{center}

\noindent
A new elliptic divisor is 
\begin{equation}
    F' = \ell_{1,5}+\ell_{1,4}+ 2 \left( \ell_1+\ell_{1,3}+\ell_3+\ell_{2,3}+\ell_2+\ell_{2,6}+\ell_6 \right)
    +\mu^{1,2}_{3,4}+\xi^{(\overline{12})(34)(16)(24)}_{(35)(36)(56)}.
\end{equation}
Then, we can put a new elliptic parameter $s$ to
\begin{equation}
    s=\frac{x+A_0+A_1 t+A_2 t^2}{t^2 \left( t- \frac{ad-bc+b-d}{(a-1)(ad-bc)(c+d-1)} \right)}.
\end{equation}
After some scaling for $s$, looking at the order of vanishing along
the non-identity components of fibers at $t=0$ and
$t=\frac{ad-bc+b-c}{(a-1)(ad-bc)(c+d-1)}$, we get
\begin{equation}
    A_0=0, \; A_2=\frac{(1-a)(ad-bc)(c+d-1)}{ad+bc+b-d} A_1.    
\end{equation}
To determine the coefficient $A_1$, we need look at the order of
vanishing along the non-identity component $\ell_{2,5}$ of the fiber
at $t=0$. Now, we denote the equation for the class 2.5, for short, by
\begin{equation}
    \label{eq:25}
    y^2 = x^3 + a_2 x^2 + a_4 x.
\end{equation}
The equation for each component of the fiber at $t=0$ of type $\I_4^*$
is given by the following (see \cite[IV \S9]{Silverman} for detail)

\begin{center}
    \includegraphics[width=7cm]{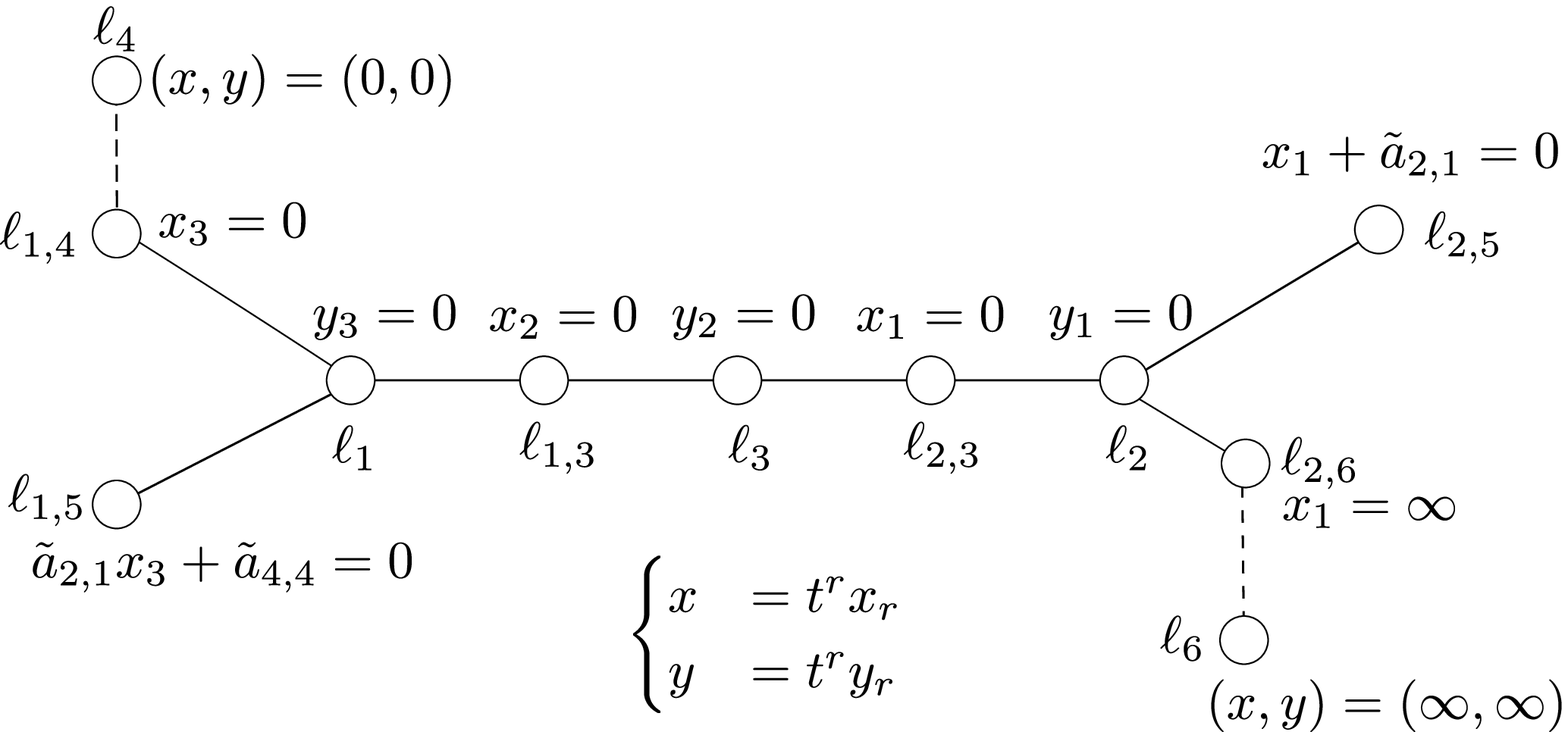}
 \end{center}
\noindent
where $a_{i,r} = t^{-r} a_i$ and $\tilde{a}_{i,r} = a_{i,r}(0)$. In
this case, we see that $\tilde{a}_{2,1}=-1,\,
\tilde{a}_{4,4}=(b-1)(ad-bc+b-d)$. Thus, substituting $x=t x_1$ and $x_1 = -1$, we see that 
\begin{equation}
    s= \frac{1+A_1+A_2t^2}{t \left( (a-1)(ad-bc)(c+d-1)t-(ad-bc+b-d) \right)}
\end{equation}
on $\ell_{2,5}$. Since $s$ has a pole along $\ell_{2,5}$ unless
$A_1=-1$, we get $A_1=-1$.

As a consequence, we have a new elliptic parameter
\begin{equation}
    s = \frac{(ad-bc)(x-t)+(a-1)(ad-bc)(c+d-1) t^2}{t^2 \left( (a-1)(ad-bc)(c+d-1) t -ad+bc-b+d \right)}.
\end{equation}
Therefore, we can compute a Weierstrass equation for the class 2.4 by
using a 2-neighbor step from the class 2.5. However, we omit it, since
it is too long to write down here. The configuration of the class 2.4 is the following.

\begin{figure}[htbp]
    \includegraphics[width=9cm]{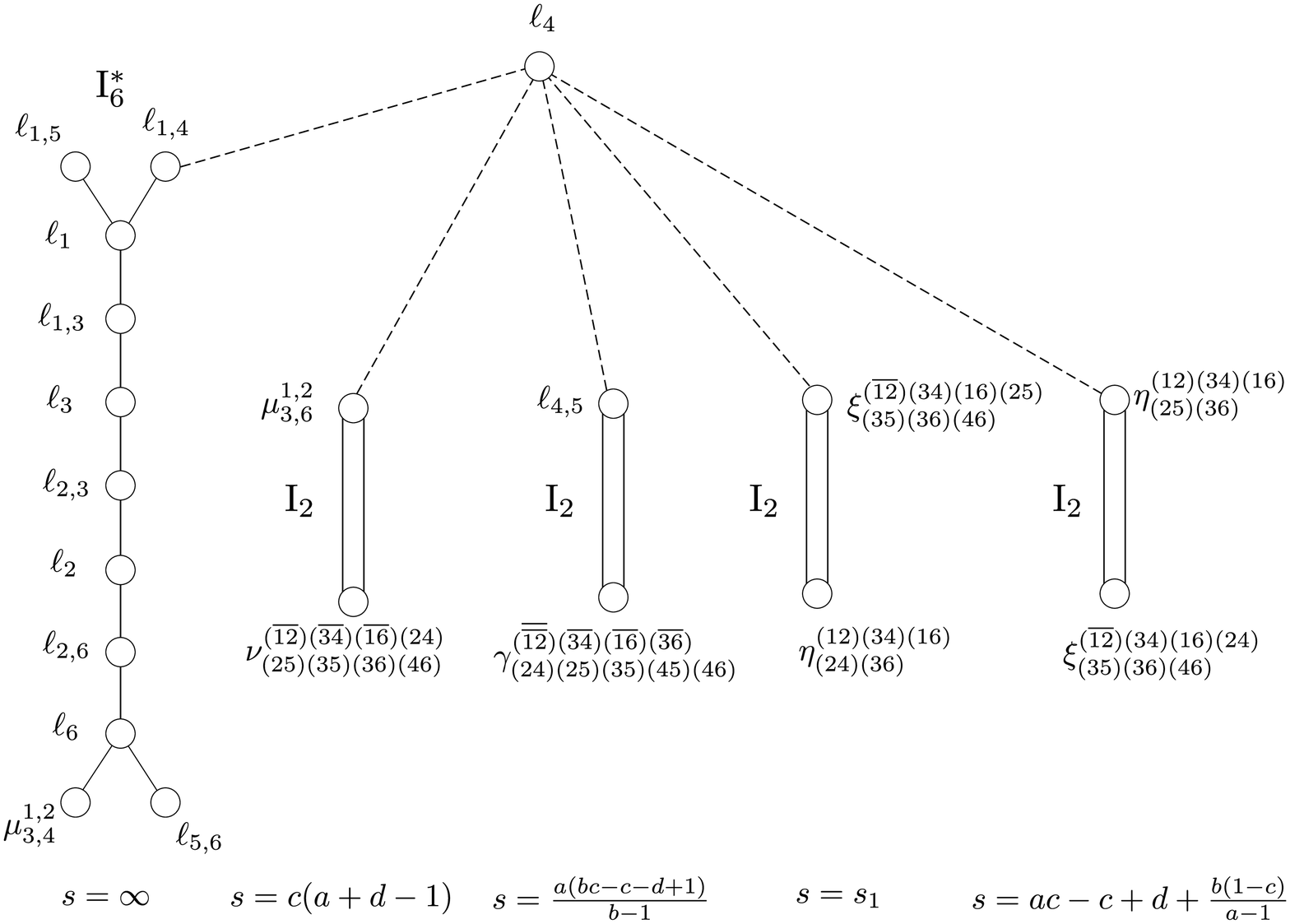}
    \caption{Configuration of singular fibers for the class 2.4}
\end{figure}
\noindent
In this picture, the value of $s_1$ is given by
\begin{equation}
    s_1 = \dfrac{(c^2+2cd-2c-d+1)a+(d-1)(bc-b+d)}{c+d-1}.
\end{equation}

\clearpage

\begin{acknowledgement}
  The author would like to express his gratitude to Professor Masoto
  Kuwata for his helpful comments and beneficial advice. The author
  also would like to thank Professors Hiroyuki Ito and Ichiro Shimada
  for their very valuable suggestions. The computer algebra system
  Maple and Maple Library ``Elliptic Surface Calculator'' written by
  Professor Masato Kuwata~\cite{ESC} were used in the calculation for
  this paper. Furthermore the author would like to thank Professor Matthias
  Sch\"{u}tt for his profitable comments to previous version of this
  paper.
\end{acknowledgement}


\newpage

\section{Weiertstrass equations}\label{sec:WEQ}

\subsection{Class 1.1 : $\I_{10}+\I_2+a\II+b\I_1$ with MWG = $\mathbb{Z}^4$}
\begin{align*}
        y^2 &= x^3 - \Big(27\,{s}^{8}-108\, ( ad+bc-2\,a-2\,b+c+d+1 ) {s}^{6}+54\,( 3\,{a}^{2}{d}^{2}+2\,abcd\\
        & \qquad +3\,{b}^{2}{c}^{2}-8\,{a}^{2}d-16\,abc-16\,abd+14\,acd+2\,a{d}^{2}-8\,{b}^{2}c+2\,b{c}^{2}+14\,bcd\\
        & \qquad +8\,{a}^{2}+24\,ab-8\,ac-2\,ad+8\,{b}^{2}-2\,bc-8\,bd+3\,{c}^{2}-6\,cd+3\,{d}^{2}-8\,a\\
        & \qquad -8\,b+2\,c+2\,d+3 ) {s}^{4}-108\, ( {a}^{3}{d}^{3}-{a}^{2}bc{d}^{2}-a{b}^{2}{c}^{2}d+{b}^{3}{c}^{3}-2\,{a}^{3}{d}^{2}\\
        & \qquad -10\,{a}^{2}bcd-4\,{a}^{2}b{d}^{2}+5\,{a}^{2}c{d}^{2}-{a}^{2}{d}^{3}-4\,a{b}^{2}{c}^{2}-10\,a{b}^{2}cd+12\,ab{c}^{2}d\\
        & \qquad +12\,abc{d}^{2}-2\,{b}^{3}{c}^{2}-{b}^{2}{c}^{3}+5\,{b}^{2}{c}^{2}d+16\,{a}^{2}bc+18\,{a}^{2}bd-8\,{a}^{2}cd-5\,{a}^{2}{d}^{2}\\
        & \qquad +18\,a{b}^{2}c+16\,a{b}^{2}d-10\,ab{c}^{2}-18\,abcd-10\,ab{d}^{2}+5\,a{c}^{2}d-2\,ac{d}^{2}-a{d}^{3}\\
        & \qquad -5\,{b}^{2}{c}^{2}-8\,{b}^{2}cd-b{c}^{3}-2\,b{c}^{2}d+5\,bc{d}^{2}-16\,{a}^{2}b+8\,{a}^{2}d-16\,a{b}^{2}-4\,abc\\
        & \qquad -4\,abd-2\,a{c}^{2}+8\,a{d}^{2}+8\,{b}^{2}c+8\,b{c}^{2}-2\,b{d}^{2}+{c}^{3}-3\,{c}^{2}d-3\,c{d}^{2}+{d}^{3}\\
        & \qquad +14\,ab+4\,ac-5\,ad-5\,bc+4\,bd-{c}^{2}+4\,cd-{d}^{2}-2\,a-2\,b-c-d+1 ) {s}^{2}\\
        & \qquad +27\, ( {a}^{2}{d}^{2}-2\,abcd+{b}^{2}{c}^{2}+4\,abc+4\,abd-2\,acd-2\,a{d}^{2}-2\,b{c}^{2}-2\,bcd\\
        & \qquad -4\,ab+2\,ad+2\,bc+{c}^{2}+2\,cd+{d}^{2}-2\,c-2\,d+1 ) ^{2} \Big) x\\
        & + \Big( 54\,{s}^{12}-324\, ( ad+bc-2\,a-2\,b+c+d+1 ) {s}^{10}+162\,( 5\,{a}^{2}{d}^{2}+6\,abcd+5\,{b}^{2}{c}^{2}\\
        & \qquad -16\,{a}^{2}d-24\,ac-24\,abd+18\,acd+6\,a{d}^{2}-16\,{b}^{2}c+6\,b{c}^{2}+18\,bcd+16\,{a}^{2}\\
        & \qquad +40\,ab-16\,ac-6\,ad+16\,{b}^{2}-6\,bc-16\,bd+5\,{c}^{2}-2\,cd+5\,{d}^{2}-16\,a-16\,b\\
        & \qquad +6\,c+6\,d+5 ) {s}^{8}-216\, ( 5\,{a}^{3}{d}^{3}+3\,{a}^{2}bc{d}^{2}+3\,a{b}^{2}{c}^{2}d+5\,{b}^{3}{c}^{3}-18\,{a}^{3}{d}^{2}\\
        & \qquad-45\,{a}^{2}bcd-33\,{a}^{2}b{d}^{2}+30\,{a}^{2}c{d}^ {2}+3\,{a}^{2}{d}^{3}-33\,a{b}^{2}{c}^{2}-45\,a{b}^{2}cd+39\,ab{c}^{2}d\\
        & \qquad +39\,abc{d}^{2}-18\,{b}^{3}{c}^{2}+3\,{b}^{2}{c}^{3}+30\,{b}^{2}{c}^{ 2}d+24\,{a}^{3}d+72\,{a}^{2}bc+111\,{a}^{2}bd-66\,{a}^{2}cd\\
        & \qquad -15\,{a}^{2 }{d}^{2}+111\,a{b}^{2}c+72\,a{b}^{2}d-45\,ab{c}^{2}-144\,abcd-45\,ab{d}^{2}+30\,a{c}^{2}d+6\,ac{d}^{2}\\
        & \qquad +3\,a{d}^{3}+24\,{b}^{3}c-15\,{b}^{2}{ c}^{2}-66\,{b}^{2}cd+3\,b{c}^{3}+6\,b{c}^{2}d+30\,bc{d}^{2}-16\,{a}^{3}-96\,{a}^{2}b\\
        & \qquad +24\,{a}^{2}c+6\,{a}^{2}d-96\,a{b}^{2}+12\,abc+12\,abd- 18\,a{c}^{2}+33\,acd+6\,a{d}^{2}-16\,{b}^{3}+6\,{b}^{2}c+\\
        & \qquad 24\,{b}^{2}d+ 6\,b{c}^{2}+33\,bcd-18\,b{d}^{2}+5\,{c}^{3}-12\,{c}^{2}d-12\,c{d}^{2}+ 5\,{d}^{3}+24\,{a}^{2}+81\,ab\\
        & \qquad -12\,ac-15\,ad+24\,{b}^{2}-15\,bc-12\,bd+3\,{c}^{2}-3\,cd+3\,{d}^{2}-18\,a-18\,b+3\,c+3\,d\\
        & \qquad +5 ) {s}^{6}+162\, ( 5\,{a}^{4}{d}^{4}-4\,{a}^{3}bc{d}^{3}-2\,{a}^{2}{b}^{2}{c }^{2}{d}^{2}-4\,a{b}^{3}{c}^{3}d+5\,{b}^{4}{c}^{4}-16\,{a}^{4}{d}^{3}\\
        & \qquad - 44\,{a}^{3}bc{d}^{2}-28\,{a}^{3}b{d}^{3}+32\,{a}^{3}c{d}^{3}-4\,{a}^{3}{d}^{4}-40\,{a}^{2}{b}^{2}{c}^{2}d-40\,{a}^{2}{b}^{2}c{d}^{2}+44\,{a} ^{2}b{c}^{2}{d}^{2}\\
        & \qquad +56\,{a}^{2}bc{d}^{3}-28\,a{b}^{3}{c}^{3}-44\,a{b}^{3}{c}^{2}d+56\,a{b}^{2}{c}^{3}d+44\,a{b}^{2}{c}^{2}{d}^{2}-16\,{b}^{4 }{c}^{3}-4\,{b}^{3}{c}^{4}\\
        & \qquad +32\,{b}^{3}{c}^{3}d+16\,{a}^{4}{d}^{2}+144 \,{a}^{3}bcd+124\,{a}^{3}b{d}^{2}-80\,{a}^{3}c{d}^{2}-12\,{a}^{3}{d}^{3}+136\,{a}^{2}{b}^{2}{c}^{2}\\
        & \qquad +376\,{a}^{2}{b}^{2}cd+136\,{a}^{2}{b}^{2 }{d}^{2}-296\,{a}^{2}b{c}^{2}d-372\,{a}^{2}bc{d}^{2}-40\,{a}^{2}b{d}^{3}+94\,{a}^{2}{c}^{2}{d}^{2}\\
        & \qquad -16\,{a}^{2}c{d}^{3}-2\,{a}^{2}{d}^{4}+124 \,a{b}^{3}{c}^{2}+144\,a{b}^{3}cd-40\,a{b}^{2}{c}^{3}-372\,a{b}^{2}{c} ^{2}d-296\,a{b}^{2}c{d}^{2}\\
        & \qquad +44\,ab{c}^{3}d+208\,ab{c}^{2}{d}^{2}+44\,a bc{d}^{3}+16\,{b}^{4}{c}^{2}-12\,{b}^{3}{c}^{3}-80\,{b}^{3}{c}^{2}d-2 \,{b}^{2}{c}^{4}\\
        & \qquad -16\,{b}^{2}{c}^{3}d+94\,{b}^{2}{c}^{2}{d}^{2}-128\,{a }^{3}bc-208\,{a}^{3}bd+64\,{a}^{3}cd+64\,{a}^{3}{d}^{2}-416\,{a}^{2}{b }^{2}c\\
        & \qquad -416\,{a}^{2}{b}^{2}d+144\,{a}^{2}b{c}^{2}+400\,{a}^{2}bcd+144\,{a}^{2}b{d}^{2}-80\,{a}^{2}{c}^{2}d-4\,{a}^{2}c{d}^{2}+4\,{a}^{2}{d}^{ 3}\\
        & \qquad -208\,a{b}^{3}c-128\,a{b}^{3}d+144\,a{b}^{2}{c}^{2}+400\,a{b}^{2}cd+144\,a{b}^{2}{d}^{2}-44\,ab{c}^{3}+44\,ab{c}^{2}d\\
        & \qquad +44\,abc{d}^{2}-44\,ab{d}^{3}+32\,a{c}^{3}d-108\,a{c}^{2}{d}^{2}-4\,a{d}^{4}+64\,{b}^{3}{c}^{2}+64\,{b}^{3}cd+4\,{b}^{2}{c}^{3}\\
        & \qquad -4\,{b}^{2}{c}^{2}d-80\,{b}^{2}c{d }^{2}-4\,b{c}^{4}-108\,b{c}^{2}{d}^{2}+32\,bc{d}^{3}+128\,{a}^{3}b-64 \,{a}^{3}d+296\,{a}^{2}{b}^{2}\\
        & \qquad +32\,{a}^{2}bc+40\,{a}^{2}bd+16\,{a}^{2} {c}^{2}+16\,{a}^{2}cd-66\,{a}^{2}{d}^{2}+128\,a{b}^{3}+40\,a{b}^{2}c+32\,a{b}^{2}d\\
        & \qquad -76\,ab{c}^{2}-212\,abcd-76\,ab{d}^{2}-16\,a{c}^{3}+20\,a{c}^{2}d+40\,ac{d}^{2}+28\,a{d}^{3}-64\,{b}^{3}c\\
        & \qquad -66\,{b}^{2}{c}^{2}+16 \,{b}^{2}cd+16\,{b}^{2}{d}^{2}+28\,b{c}^{3}+40\,b{c}^{2}d+20\,bc{d}^{2 }-16\,b{d}^{3}+5\,{c}^{4}-16\,{c}^{3}d\\
        & \qquad +30\,{c}^{2}{d}^{2}-16\,c{d}^{3} +5\,{d}^{4}-176\,{a}^{2}b-32\,{a}^{2}c+64\,{a}^{2}d-176\,a{b}^{2}+28\,abc+28\,abd\\
        & \qquad +16\,a{c}^{2}-40\,acd+4\,a{d}^{2}+64\,{b}^{2}c-32\,{b}^{2}d +4\,b{c}^{2}-40\,bcd+16\,b{d}^{2}-4\,{c}^{3}+12\,{c}^{2}d\\
        & \qquad +12\,c{d}^{2} -4\,{d}^{3}+16\,{a}^{2}+92\,ab+16\,ac-12\,ad+16\,{b}^{2}-12\,bc+16\,bd -2\,{c}^{2}+8\,cd\\
        & \qquad -2\,{d}^{2}-16\,a-16\,b-4\,c-4\,d+5 ) {s}^{4}- 324\, ( {a}^{5}{d}^{5}-3\,{a}^{4}bc{d}^{4}+2\,{a}^{3}{b}^{2}{c}^{ 2}{d}^{3}+2\,{a}^{2}{b}^{3}{c}^{3}{d}^{2}\\
        & \qquad -3\,a{b}^{4}{c}^{4}d+{b}^{5}{ c}^{5}-2\,{a}^{5}{d}^{4}-2\,{a}^{4}bc{d}^{3}+3\,{a}^{4}c{d}^{4}-3\,{a} ^{4}{d}^{5}+10\,{a}^{3}{b}^{2}{c}^{2}{d}^{2}-6\,{a}^{3}{b}^{2}c{d}^{3}\\
        & \qquad  +2\,{a}^{3}b{c}^{2}{d}^{3}+14\,{a}^{3}bc{d}^{4}-6\,{a}^{2}{b}^{3}{c}^{3}d+10\,{a}^{2}{b}^{3}{c}^{2}{d}^{2}-16\,{a}^{2}{b}^{2}{c}^{3}{d}^{2}- 16\,{a}^{2}{b}^{2}{c}^{2}{d}^{3}\\
        & \qquad -2\,a{b}^{4}{c}^{3}d+14\,a{b}^{3}{c}^{4}d+2\,a{b}^{3}{c}^{3}{d}^{2}-2\,{b}^{5}{c}^{4}-3\,{b}^{4}{c}^{5}+3\,{ b}^{4}{c}^{4}d+8\,{a}^{4}bc{d}^{2}+6\,{a}^{4}b{d}^{3}\\
        & \qquad -4\,{a}^{4}c{d}^{ 3}+{a}^{4}{d}^{4}+72\,{a}^{3}{b}^{2}{c}^{2}d+74\,{a}^{3}{b}^{2}c{d}^{2 }-94\,{a}^{3}b{c}^{2}{d}^{2}-104\,{a}^{3}bc{d}^{3}-6\,{a}^{3}b{d}^{4}\\
        & \qquad -4\,{a}^{3}{c}^{2}{d}^{3}-8\,{a}^{3}c{d}^{4}+2\,{a}^{3}{d}^{5}+74\,{a}^{2}{b}^{3}{c}^{2}d+72\,{a}^{2}{b}^{3}c{d}^{2}-56\,{a}^{2}{b}^{2}{c}^{3 }d-114\,{a}^{2}{b}^{2}{c}^{2}{d}^{2}\\
        & \qquad -56\,{a}^{2}{b}^{2}c{d}^{3}+98\,{a }^{2}b{c}^{3}{d}^{2}+88\,{a}^{2}b{c}^{2}{d}^{3}-16\,{a}^{2}bc{d}^{4}+6\,a{b}^{4}{c}^{3}+8\,a{b}^{4}{c}^{2}d-6\,a{b}^{3}{c}^{4}\\
        & \qquad -104\,a{b}^{3} {c}^{3}d-94\,a{b}^{3}{c}^{2}{d}^{2}-16\,a{b}^{2}{c}^{4}d+88\,a{b}^{2}{ c}^{3}{d}^{2}+98\,a{b}^{2}{c}^{2}{d}^{3}+{b}^{4}{c}^{4}-4\,{b}^{4}{c}^ {3}d\\
        & \qquad +2\,{b}^{3}{c}^{5}-8\,{b}^{3}{c}^{4}d-4\,{b}^{3}{c}^{3}{d}^{2}-8\, {a}^{4}b{d}^{2}+4\,{a}^{4}{d}^{3}-80\,{a}^{3}{b}^{2}{c}^{2}-224\,{a}^{3}{b}^{2}cd-72\,{a}^{3}{b}^{2}{d}^{2}\\
        & \qquad +80\,{a}^{3}b{c}^{2}d+248\,{a}^{3}bc{d}^{2}+80\,{a}^{3}b{d}^{3}+12\,{a}^{3}{c}^{2}{d}^{2}+30\,{a}^{3}c{d}^{3}+12\,{a}^{3}{d}^{4}-72\,{a}^{2}{b}^{3}{c}^{2}\\
        & \qquad -224\,{a}^{2}{b}^{3}cd-80\,{a}^{2}{b}^{3}{d}^{2}+72\,{a}^{2}{b}^{2}{c}^{3}+156\,{a}^{2}{b}^{2}{c}^{2}d+156\,{a}^{2}{b}^{2}c{d}^{2}+72\,{a}^{2}{b}^{2}{d}^{3}-94\,{a}^{2}b{c}^{3}d\\
        & \qquad -174\,{a}^{2}b{c}^{2}{d}^{2}+44\,{a}^{2}bc{d}^{3}+10 \,{a}^{2}b{d}^{4}-4\,{a}^{2}{c}^{3}{d}^{2}+6\,{a}^{2}c{d}^{4}+2\,{a}^{2}{d}^{5}-8\,a{b}^{4}{c}^{2}+80\,a{b}^{3}{c}^{3}\\
        & \qquad +248\,a{b}^{3}{c}^{2}d+80\,a{b}^{3}c{d}^{2}+10\,a{b}^{2}{c}^{4}+44\,a{b}^{2}{c}^{3}d-174\,a{b}^{2}{c}^{2}{d}^{2}-94\,a{b}^{2}c{d}^{3}+2\,ab{c}^{4}d\\
        & \qquad -102\,ab{c}^{3}{d}^{2}-102\,ab{c}^{2}{d}^{3}+2\,abc{d}^{4}+4\,{b}^{4}{c}^{3}+12\,{b}^{3}{c}^{4}+30\,{b}^{3}{c}^{3}d+12\,{b}^{3}{c}^{2}{d}^{2}+2\,{b}^{2}{c}^{5}\\
        & \qquad +6\,{b}^{2}{c}^{4}d-4\,{b}^{2}{c}^{2}{d}^{3}+160\,{a}^{3}{b}^{2}c+152\,{a}^{3}{b}^{2}d-160\,{a}^{3}bcd-154\,{a}^{3}b{d}^{2}-24\,{a}^{3}c{d}^{2}\\
        & \qquad -26\,{a}^{3}{d}^{3}+152\,{a}^{2}{b}^{3}c+160\,{a}^{2}{b}^{3}d-56\,{a}^{2}{b}^{2}{c}^{2}+4\,{a}^{2}{b}^{2}cd-56\,{a}^{2}{b}^{2}{d}^{2}+8\,{a}^{2}b{c}^{3}\\
        & \qquad +134\,{a}^{2}b{c}^{2}d-66\,{a}^{2}bc{d}^{2}-90\,{a}^{2}b{d}^{3}-4\,{a}^{2}{c}^{3}d-18\,{a}^{2}{c}^{2}{d}^{2}-42\,{a}^{2}c{d}^{3}-22\,{a}^{2}{d}^{4}\\
        & \qquad -154\,a{b}^{3}{c}^{2}-160\,a{b}^{3}cd-90\,a{b}^{2}{c}^{3}-66\,a{b}^{2}{c}^{2}d+134\,a{b}^{2}c{d}^{2}+8\,a{b}^{2}{d}^{3}-2\,ab{c}^{4}\\
        & \qquad +72\,ab{c}^{3}d+256\,ab{c}^{2}{d}^{2}+72\,abc{d}^{3}-2\,ab{d}^{4}+3\,a{c}^{4}d+12\,a{c}^{3}{d}^{2}+12\,a{c}^{2}{d}^{3}\\
        & \qquad -3\,a{d}^{5}-26\,{b}^{3}{c}^{3}-24\,{b}^{3}{c}^{2}d-22\,{b}^{2}{c}^{4}-42\,{b}^{2}{c}^{3}d-18\,{b}^{2}{c}^{2}{d}^{2}-4\,{b}^{2}c{d}^{3}-3\,b{c}^{5}\\
        & \qquad +12\,b{c}^{3}{d}^{2}+12\,b{c}^{2}{d}^{3}+3\,bc{d}^{4}-80\,{a}^{3}{b}^{2}+80\,{a}^{3}bd+12\,{a}^{3}{d}^{2}-80\,{a}^{2}{b}^{3}\\
        & \qquad -104\,{a}^{2}{b}^{2}c-104\,{a}^{2}{b}^{2}d-24\,{a}^{2}b{c}^{2}+14\,{a}^{2}bcd+142\,{a}^{2}b{d}^{2}+12\,{a}^{2}{c}^{2}d+48\,{a}^{2}c{d}^{2}\\
        & \qquad +42\,{a}^{2}{d}^{3}+80\,a{b}^{3}c+142\,a{b}^{2}{c}^{2}+14\,a{b}^{2}cd-24\,a{b}^{2}{d}^{2}-4\,ab{c}^{3}-178\,ab{c}^{2}d\\
        & \qquad -178\,abc{d}^{2}-4\,ab{d}^{3}-2 \,a{c}^{4}-10\,a{c}^{3}d-12\,a{c}^{2}{d}^{2}+10\,ac{d}^{3}+14\,a{d}^{4 }+12\,{b}^{3}{c}^{2}\\
        & \qquad +42\,{b}^{2}{c}^{3}+48\,{b}^{2}{c}^{2}d+12\,{b}^{2 }c{d}^{2}+14\,b{c}^{4}+10\,b{c}^{3}d-12\,b{c}^{2}{d}^{2}-10\,bc{d}^{3}-2\,b{d}^{4}\\
        & \qquad +{c}^{5}-{c}^{4}d-8\,{c}^{3}{d}^{2}-8\,{c}^{2}{d}^{3}-c{d}^{4}+{d}^{5}+88\,{a}^{2}{b}^{2}+24\,{a}^{2}bc-54\,{a}^{2}bd-12\,{a}^{2 }cd\\
        & \qquad -26\,{a}^{2}{d}^{2}-54\,a{b}^{2}c+24\,a{b}^{2}d+24\,ab{c}^{2}+132\,abcd+24\,ab{d}^{2}+8\,a{c}^{3}+12\,a{c}^{2}d\\
        & \qquad -12\,ac{d}^{2}-22\,a{d}^{3}-26\,{b}^{2}{c}^{2}-12\,{b}^{2}cd-22\,b{c}^{3}-12\,b{c}^{2}d+12\,bc{d}^{2}+8\,b{d}^{3}-3\,{c}^{4}\\
        & \qquad +6\,{c}^{3}d+18\,{c}^{2}{d}^{2}+6\,c{d}^{3}-3\,{d}^{4}-8\,{a}^{2}b+4\,{a}^{2}d-8\,a{b}^{2}-28\,abc-28\,abd-12\,a{c}^{2}\\
        & \qquad -6\,acd+12\,a{d}^{2}+4\,{b}^{2}c+12\,b{c}^{2}-6\,bcd-12\,b{d}^{2}+2\,{c}^{3}-12\,{c}^{2}d-12\,c{d}^{2}+2\,{d}^{3}\\
        & \qquad +10\,ab+8\,ac+ad+bc+8\,bd+2\,{c}^{2}+10\,cd+2\,{d}^{2}-2\,a-2\,b-3\,c-3\,d+1 ) {s}^{2}\\
        & \qquad +54\, ( {a}^{2}{d}^{2}-2\,abcd+{b}^{2}{c}^{2}+4\,abc+4\,abd-2\,acd-2\,a{d}^{2}-2\,b{c}^{2}-2\,bcd-4\,ab\\
        & \qquad +2\,ad+2\,bc+{c}^{2}+2\,cd+{d}^{2}-2\,c-2\,d+1 ) ^{3} \Big).
\end{align*}

\subsection{Class 1.2 : $\I_{10}+\I_2+a\II+b\I_1$ with MWG = $\mathbb{Z}^4$}
\begin{align*}
    y^2 &= x^3 - \Big( 27\, (ad-bc)^{4}{t}^{8}-108\,({a}^{3}{d}^{3}-{a}^{2}bc{d}^{2}-a{b}^{2}{c}^{2}d+{b}^{3}{c}^{3}-2\,{a}^{3}{d}^{2}-2\,{a}^{2}bcd\\
    & \qquad +4\,{a}^{2}b{d}^{2}+4\,{a}^{2}c{d}^{2}-2\,{a}^{2}{d}^{3}+4\,a{b}^{2}{c}^{2}-2\,a{b}^{2}cd-2\,ab{c}^{2}d-2\,abc{d}^{2}-2\,{b}^{3}{c}^{2}\\
    & \qquad -2\,{b}^{2}{c}^{3}+4\,{b}^{2}{c}^{2}d-6\,{a}^{2}bd+4\,{a}^{2}{d}^{2}-6\,a{b}^{2}c+16\,abcd-6\,ac{d}^{2}+4\,{b}^{2}{c}^{2}-6\,b{c}^{2}d) {t}^{6}\\
    & \qquad +54\, ( 3\,{a}^{2}{d}^{2}+2\,abcd+3\,{b}^{2}{c}^{2}-8\,{a}^{2}d+4\,abc+4\,abd+4\,acd-8\,a{d}^{2}-8\,{b}^{2}c\\
    & \qquad -8\,b{c}^{2}+4\,bcd+8\,{a}^{2}+4\,ab-8\,ac+8\,ad+8\,{b}^{2}+8\,bc-8\,bd+8\,{c}^{2}+4\,cd+8\,{d}^{2}\\
    & \qquad -8\,a-8\,b-8\,c-8\,d+8 ) {t}^{4}-108\, ( ad+bc-2\,a-2\,b-2\,c-2\,d+4 ) {t}^{2}+27 \Big)x \\
    & + \Big( 54\, ( ad-bc ) ^{6}{t}^{12}-324\, ( {a}^{3}{d}^{3}-{a}^{2}bc{d}^{2}-a{b}^{2}{c}^{2}d+{b}^{3}{c}^{3}-2\,{a}^{3}{d}^{2}-2\,{a}^{2}bcd+4\,{a}^{2}b{d}^{2}\\
    & \qquad +4\,{a}^{2}c{d}^{2}-2\,{a}^{2}{d}^{3}+4\,a{b}^{2}{c}^{2}-2\,a{b}^{2}cd-2\,ab{c}^{2}d-2\,abc{d}^{2}-2\,{b}^{3}{c}^{2}-2\,{b}^{2}{c}^{3}+4\,{b}^{2}{c}^{2}d\\
    & \qquad -6\,{a}^{2}bd+4\,{a}^{2}{d}^{2}-6\,a{b}^{2}c+16\,abcd-6\,ac{d}^{2}+4\,{b}^{2}{c}^{2}-6\,b{c}^{2}d)(ad-bc)^{2}{t}^{10}\\
    & \qquad +162\, ( 5\,{a}^{4}{d}^{4}-4\,{a}^{3}bc{d}^{3}-2\,{a}^{2}{b}^{2}{c}^{2}{d}^{2}-4\,a{b}^{3}{c}^{3}d+5\,{b}^{4}{c}^{4}-16\,{a}^{4}{d}^{3}+4\,{a}^{3}bc{d}^{2}\\
    & \qquad +20\,{a}^{3}b{d}^{3}+20\,{a}^{3}c{d}^{3}-16\,{a}^{3}{d}^{4}-8\,{a}^{2}{b}^{2}{c}^{2}d-8\,{a}^{2}{b}^{2}c{d}^{2}-8\,{a}^{2}b{c}^{2}{d}^{2}+4\,{a}^{2}bc{d}^{3}\\
    & \qquad +20\,a{b}^{3}{c}^{3}+4\,a{b}^{3}{c}^{2}d+4\,a{b}^{2}{c}^{3}d-8\,a{b}^{2}{c}^{2}{d}^{2}-16\,{b}^{4}{c}^{3}-16\,{b}^{3}{c}^{4}+20\,{b}^{3}{c}^{3}d+16\,{a}^{4}{d}^{2}\\
    & \qquad +16\,{a}^{3}bcd-52\,{a}^{3}b{d}^{2}-40\,{a}^{3}c{d}^{2}+40\,{a}^{3}{d}^{3}+40\,{a}^{2}{b}^{2}{c}^{2}+152\,{a}^{2}{b}^{2}cd+40\,{a}^{2}{b}^{2}{d}^{2}\\
    & \qquad -64\,{a}^{2}b{c}^{2}d-88\,{a}^{2}bc{d}^{2}-40\,{a}^{2}b{d}^{3}+40\,{a}^{2}{c}^{2}{d}^{2}-52\,{a}^{2}c{d}^{3}+16\,{a}^{2}{d}^{4}-52\,a{b}^{3}{c}^{2}\\
    & \qquad +16\,a{b}^{3}cd-40\,a{b}^{2}{c}^{3}-88\,a{b}^{2}{c}^{2}d-64\,a{b}^{2}c{d}^{2}+16\,ab{c}^{3}d+152\,ab{c}^{2}{d}^{2}+16\,abc{d}^{3}\\
    & \qquad +16\,{b}^{4}{c}^{2}+40\,{b}^{3}{c}^{3}-40\,{b}^{3}{c}^{2}d+16\,{b}^{2}{c}^{4}-52\,{b}^{2}{c}^{3}d+40\,{b}^{2}{c}^{2}{d}^{2}+48\,{a}^{3}bd-40\,{a}^{3}{d}^{2}\\
    & \qquad -96\,{a}^{2}{b}^{2}c-96\,{a}^{2}{b}^{2}d-16\,{a}^{2}bcd+104\,{a}^{2}b{d}^{2}+104\,{a}^{2}c{d}^{2}-40\,{a}^{2}{d}^{3}+48\,a{b}^{3}c\\
    & \qquad +104\,a{b}^{2}{c}^{2}-16\,a{b}^{2}cd-16\,ab{c}^{2}d-16\,abc{d}^{2}-96\,a{c}^{2}{d}^{2}+48\,ac{d}^{3}-40\,{b}^{3}{c}^{2}-40\,{b}^{2}{c}^{3}\\
    & \qquad +104\,{b}^{2}{c}^{2}d+48\,b{c}^{3}d-96\,b{c}^{2}{d}^{2}+72\,{a}^{2}{b}^{2}-96\,{a}^{2}bd+40\,{a}^{2}{d}^{2}-96\,a{b}^{2}c+160\,abcd\\
    & \qquad -96\,ac{d}^{2}+40\,{b}^{2}{c}^{2}-96\,b{c}^{2}d+72\,{c}^{2}{d}^{2} ) {t}^{8}-216\, ( 5\,{a}^{3}{d}^{3}+3\,{a}^{2}bc{d}^{2}+3\,a{b}^{2}{c}^{2}d\\
    & \qquad +5\,{b}^{3}{c}^{3}-18\,{a}^{3}{d}^{2}-3\,{a}^{2}bcd+9\,{a}^{2}b{d}^{2}+9\,{a}^{2}c{d}^{2}-18\,{a}^{2}{d}^{3}+9\,a{b}^{2}{c}^{2}-3\,a{b}^{2}cd\\
    & \qquad -3\,ab{c}^{2}d-3\,abc{d}^{2}-18\,{b}^{3}{c}^{2}-18\,{b}^{2}{c}^{3}+9\,{b}^{2}{c}^{2}d+24\,{a}^{3}d-12\,{a}^{2}bc-15\,{a}^{2}bd\\
    & \qquad -24\,{a}^{2}cd+48\,{a}^{2}{d}^{2}-15\,a{b}^{2}c-12\,a{b}^{2}d-24\,ab{c}^{2}-60\,abcd-24\,ab{d}^{2}-12\,a{c}^{2}d\\
    & \qquad -15\,ac{d}^{2}+24\,a{d}^{3}+24\,{b}^{3}c+48\,{b}^{2}{c}^{2}-24\,{b}^{2}cd+24\,b{c}^{3}-15\,b{c}^{2}d-12\,bc{d}^{2}-16\,{a}^{3}\\
    & \qquad -12\,{a}^{2}b+24\,{a}^{2}c-36\,{a}^{2}d-12\,a{b}^{2}+96\,abc+96\,abd+24\,a{c}^{2}+96\,acd-36\,a{d}^{2}\\
    & \qquad -16\,{b}^{3}-36\,{b}^{2}c+24\,{b}^{2}d-36\,b{c}^{2}+96\,bcd+24\,b{d}^{2}-16\,{c}^{3}-12\,{c}^{2}d-12\,c{d}^{2}\\
    & \qquad -16\,{d}^{3}+24\,{a}^{2}-24\,ab-96\,ac-36\,ad+24\,{b}^{2}-36\,bc-96\,bd+24\,{c}^{2}-24\,cd\\
    & \qquad +24\,{d}^{2}+24\,a+24\,b+24\,c+24\,d-16 ) {t}^{6}+162\, ( 5\,{a}^{2}{d}^{2}+6\,abcd+5\,{b}^{2}{c}^{2}\\
    & \qquad -16\,{a}^{2}d-4\,abc-4\,abd-4\,acd-16\,a{d}^{2}-16\,{b}^{2}c-16\,b{c}^{2}-4\,bcd+16\,{a}^{2}\\
    & \qquad +20\,ab+8\,ac+40\,ad+16\,{b}^{2}+40\,bc+8\,bd+16\,{c}^{2}+20\,cd+16\,{d}^{2}-40\,a\\
    & \qquad -40\,b-40\,c-40\,d+40 ) {t}^{4}-324\, ( ad+bc-2\,a-2\,b-2\,c-2\,d+4 ) {t}^{2}+54 \Big).
\end{align*}

\subsection{Class 1.3 : $2 \I_6+a\II+b\I_1$ with MWG = $\mathbb{Z}^4$}
    \begin{align*}
        y^2 & = x^3 -\Big(27t^8+(216bc+216+216ad+648ac-432c-108b-108d-432a)t^6\\
        & \quad +(432-432a-432d-432b-432c+432bc+432ad+108bd-432ad^2-432b^2c\\
        & \quad +432a^2+432d^2a^2-432da^2+216abc432ac+216cd+216bcd+216abd+216ab\\
        & \quad -432bc^2+432b^2c^2+432c^2+162d^2+162b^2-432abcd+216acd)t^4\\
        & \quad -108(b-d)(b^2-2b^2c-2b+4bc+4abd-2ab-4bcd-4ad+2cd+2d+2ad^2\\
        & \quad -d^2)t^2+27(b-d)^4 \Big)x\\
        & \quad + \Big(54t^{12}+(1944ac+648bc+648ad-324b-324d-1296a+648-1296c)t^{10}\\
        & \quad +(2592-6480a-2592d-2592b-6480c-15552a^2c+7776abc^2+6480bc\\
        & \quad +6480ad+972bd-2592ad^2-2592b^2c+6480a^2+2592a^2d^2-6480a^2d-8424abc\\
        & \quad -15552ac^2+11664a^2c^2+16848ac+3240cd-648bcd-648abd+3240ba-6480bc^2\\
        & \quad +2592b^2c^2+6480c^2+810d^2+810b^2+1296abcd-8424acd+7776a^2cd)t^8\\
        & \quad +(3456-5184a-5184d-5184b-5184c+2592a^2bc+2592ab^2c^2-5184a^2c+5184abc^2\\
        & \quad +7776bc+7776ad+2592bd+2592a^2b-10368ad^2-10368b^2c-1944cd^2+3888ad^3\\
        & \quad +3888b^3c-648bd^2-648b^2d+3240ab^2c-5184a^2+3456a^3+7776a^2d^2-5184a^3d^2\\
        & \quad -5184a^3d-10368bc^2d+7776a^2d+2592b^2c^2d+2592a^2bd^2+5184ab^2cd+5184abcd^2\\
        & \quad -5184ab^2c^2d-5184a^2bcd^2+5184abc^2d+5184a^2bcd-20736cba-5184ac^2-5184b^3c^2\\
        & \quad +3456a^3d^3-5184a^2d^3+3456b^3c^3-5184b^2c^3-5184bc^3+20736ac+2592c^2d+5184cd\\
        & \quad -10368a^2bd+1944bcd+1944abd+5184ab+7776c^2b+7776c^2b^2-5184c^2+3888d^2\\
        & \quad +3888b^2+3456c^3-1080b^3-1080d^3+3888abcd+2592ac^2d-20736acd+2592a^2cd^2\\
        & \quad +3240d^2ac+5184a^2cd-1944ab^2+648b^2cd+648abd^2-1944bcd^2-1944ab^2d)t^6\\
        & \quad +(2592b^4c^2-2592ad^4+2592a^2d^4-2592b^4c-5184bd-324b^2d^2-6480ad^2\\
        & \quad -6480b^2c+1296cd^2+6480ad^3+6480b^3c+2592bd^2+2592b^2d-6480ab^2c\\
        & \quad -648bd^3-648b^3d+6480d^2a^2+2592bc^2d-10368b^2c^2d-10368a^2bd^2+20088ab^2cd\\
        & \quad +20088abcd^2-10368ab^2cd^2-6480b^3c^2+2592a^2b^2-648ab^3-6480a^2d^3+2592c^2d^2\\
        & \quad -648cd^3-6480bc^2d^2+3240ab^3c-6480a^2b^2d+3240ab^3d-1296ab^2d^2+648abd^3\\
        & \quad +648b^3cd-1296b^2cd^2+3240bcd^3+2592b^3c^2d+2592a^2bd^3+3240acd^3+6480b^2c^2d^2\\
        & \quad +6480a^2b^2d^2+2592a^2bd+5184bcd+5184abd+6480c^2b^2-6480abcd^3-6480ab^3cd\\
        & \quad +2592d^2+2592b^2+810b^4-2592b^3-2592d^3+810d^4-10368abcd-6480acd^2\\
        & \quad +1296ab^2+1944b^2cd+1944abd^2-7776bcd^2-7776ab^2d)t^4\\
        & \quad -324(b-d)^3(b^2-2b^2c-2b+4bc+4abd-2ab-4bcd-4ad+2cd+2d+2ad^2\\
        & \quad -d^2)t^2+54(b-d)^6 \Big). 
    \end{align*}

\subsection{Class 1.4 : $\IV^*+\I_4+a\II+b\I_1$ with MWG = $\mathbb{Z}^5$}

\begin{align*}
     y^2 &= x^3 - \Big( 216\, ( 2\,{a}^{2}{d}^{2}-2\,abcd+2\,{b}^{2}{c}^{2}-2\,{a}^{2}d+abc+abd+acd-2\,a{d}^{2}-2\,{b}^{2}c\\
    & \qquad -2\,b{c}^{2}+bcd+2\,{a}^{2}+ab-2\,ac+2\,ad+2\,{b}^{2}+2\,bc-2\,bd+2\,{c}^{2}+cd+2\,{d}^{2}\\
    & \qquad -2\,a-2\,b-2\,c-2\,d+2 ) {s}^{4}-216\, ( 2\,{a}^{3}{b}^{2}{c}^{2}d+{a}^{3}{b}^{2}c{d}^{2}-{a}^{3}{b}^{2}{d}^{3}-{a}^{3}b{c}^{2}{d}^{2}\\
    & \qquad -4\,{a}^{3}
    bc{d}^{3}-{a}^{3}{c}^{2}{d}^{3}-{a}^{2}{b}^{3}{c}^{3}+{a}^{2}{b}^{3}{c
    }^{2}d+2\,{a}^{2}{b}^{3}c{d}^{2}-{a}^{2}{b}^{2}{c}^{3}d+4\,{a}^{2}{b}^{2}{c}^{2}{d}^{2}\\
    & \qquad -{a}^{2}{b}^{2}c{d}^{3}+2\,{a}^{2}b{c}^{3}{d}^{2}+{a}^{2}b{c}^{2}{d}^{3}-4\,a{b}^{3}{c}^{3}d-a{b}^{3}{c}^{2}{d}^{2}+a{b}^{2}{c}^{3}{d}^{2}+2\,a{b}^{2}{c}^{2}{d}^{3}\\
    & \qquad -{b}^{3}{c}^{3}{d}^{2}-{a}^{3}{b}^{2}{c}^{2}-3\,{a}^{3}{b}^{2}cd+{a}^{3}{b}^{2}{d}^{2}-{a}^{3}b{c}^{2}d+6\,{a}^{3}bc{d}^{2}+3\,{a}^{3}b{d}^{3}+2\,{a}^{3}{c}^{2}{d}^{2}\\
    & \qquad +3\,{a}^{3}c{d}^{3}+{a}^{2}{b}^{3}{c}^{2}-3\,{a}^{2}{b}^{3}cd-{a}^{2}{b}^{3}{d}^{2}+2\,{a}^{2}{b}^{2}{c}^{3}-7\,{a}^{2}{b}^{2}{c}^{2}d-7\,{a}^{2}{b}^{2}c{d}^{2}\\
    & \qquad +2\,{a}^{2}{b}^{2}{d}^{3}-{a}^{2}b{c}^{3}d-7\,{a}^{2}b{c}^{2}{d}^{2}+6\,{a}^{2}bc{d}^{3}-{a}^{2}{c}^{3}{d}^{2}+{a}^{2}{c}^{2}{d}^{3}+3\,a{b}^{3}{c}^{3}\\
    & \qquad +6\,a{b}^{3}{c}^{2}d-a{b}^{3}c{d}^{2}+6\,a{b}^{2}{c}^{3}d-7\,a{b}^{2}{c}^{2}{d}^{2}-a{b}^{2}c{d}^{3}-3\,ab{c}^{3}{d}^{2}-3\,ab{c}^{2}{d}^{3}\\
    & \qquad +3\,{b}^{3}{c}^{3}d+2\,{b}^{3}{c}^{2}{d}^{2}+{b}^{2}{c}^{3}{d}^{2}-{b}^{2}{c}^{2}{d}^{3}+2\,{a}^{3}{b}^{2}c+{a}
    ^{3}{b}^{2}d-2\,{a}^{3}bcd\\
    & \qquad -6\,{a}^{3}b{d}^{2}-3\,{a}^{3}c{d}^{2}+{a}^{2}{b}^{3}c+2\,{a}^{2}{b}^{3}d-2\,{a}^{2}{b}^{2}{c}^{2}+10\,{a}^{2}{b}^{2}cd-2\,{a}^{2}{b}^{2}{d}^{2}\\
    & \qquad +10\,{a}^{2}b{c}^{2}d-3\,{a}^{2}b{d}^{3}-2\,{a}^{2}{c}^{2}{d}^{2}-6\,{a}^{2}c{d}^{3}-6\,a{b}^{3}{c}^{2}-2\,a{b}^{3}cd-3\,a{b}^{2}{c}^{3}\\
    & \qquad +10\,a{b}^{2}c{d}^{2}-2\,ab{c}^{3}d+10\,ab{c}^{2}{d}^{2}-2\,abc{d}^{3}+2\,a{c}^{3}{d}^{2}+a{c}^{2}{d}^{3}-3\,{b}^{3}{c}^{2}d\\
    & \qquad -6\,{b}^{2}{c}^{3}d-2\,{b}^{2}{c}^{2}{d}^{2}+b{c}^{3}{d}^{2}+2\,b{c}^{2}{d}^{3}-{a}^{3}{b}^{2}+3\,{a}^{3}bd-{a}^{2}{b}^{3}-2\,{a}^{2}{b}^{2}c\\
    & \qquad -2\,{a}^{2}{b}^{2}d-6\,{a}^{2}bcd+6\,{a}^{2}b{d}^{2}+6\,{a}^{2}c{d}^{2}+3\,a{b}^{3}c+6\,a{b}^{2}{c}^{2}-6\,a{b}^{2}cd\\
    & \qquad -6\,ab{c}^{2}d-6\,abc{d}^{2}-2\,a{c}^{2}{d}^{2}+3\,ac{d}^{3}+6\,{b}^{2}{c}^{2}d+3\,b{c}^{3}d-2\,b{c}^{2}{d}^{2}-{c}^{3}{d}^{2}\\
    & \qquad -{c}^{2}{d}^{3}+2\,{a}^{2}{b}^{2}-3\,{a}^{2}bd-3\,a{b}^{2}c+8\,abcd-3\,ac{d}^{2}-3\,b{c}^{2}d+2\,{c}^{2}{d}^{2} ) {s}^{2}\\
    & \qquad +27\, ( abc+abd-acd-bcd-ab+cd ) ^{4} \Big)x \\
    & + \Big( 11664\,{s}^{8}+864\, ( 4\,{a}^{3}{d}^{3}-6\,{a}^{2}bc{d}^{2}-6\,a{b}^{2}{c}^{2}d+4\,{b}^{3}{c}^{3}-6\,{a}^{3}{d}^{2}+6\,{a}^{2}bcd\\
    & \qquad +3\,{a}^{2}b{d}^{2}+3\,{a}^{2}c{d}^{2}-6\,{a}^{2}{d}^{3}+3\,a{b}^{2}{c}^{2}+6\,a{b}^{2}cd+6\,ab{c}^{2}d+6\,abc{d}^{2}-6\,{b}^{3}{c}^{2}-6\,{b}^{2}{c}^{3}\\
    & \qquad +3\,{b}^{2}{c}^{2}d-6\,{a}^{3}d+3\,{a}^{2}bc-12\,{a}^{2}bd+6\,
    {a}^{2}cd+9\,{a}^{2}{d}^{2}-12\,a{b}^{2}c+3\,a{b}^{2}d+6\,ab{c}^{2}\\
    & \qquad+36\,abcd+6\,ab{d}^{2}+3\,a{c}^{2}d-12\,ac{d}^{2}-6\,a{d}^{3}-6\,{b}^{3}c+9\,{b}^{2}{c}^{2}+6\,{b}^{2}cd-6\,b{c}^{3}\\
    & \qquad -12\,b{c}^{2}d+3\,bc{d}^{2}+4\,{a}^{3}+3\,{a}^{2}b-6\,{a}^{2}c+9\,{a}^{2}d+3\,a{b}^{2}-24\,abc-24\,abd\\
    & \qquad -6\,a{c}^{2}-24\,acd+9\,a{d}^{2}+4\,{b}^{3}+9\,{b}^{2}c-6\,{b}^{2}d+9\,b{c}^{2}-24\,bcd-6\,b{d}^{2}+4\,{c}^{3}\\
    & \qquad +3\,{c}^{2}d+3\,c{d}^{2}+4\,{d}^{3}-6\,{a}^{2}+6\,ab+24\,ac+9\,ad-6\,{b}^{2}+9\,bc+24\,bd-6\,{c}^{2}\\
    & \qquad +6\,cd-6\,{d}^{2}-6\,a-6\,b-6\,c-6\,d+4 ) {s}^{6}+648\, ( 10\,{a}^{4}{b}^{2}{c}^{2}{d}^{2}-4\,{a}^{4}{b}^{2}c{d}^{3}\\
    & \qquad +4\,{a}^{4}{b}^{2}{d}^{4}+4\,{a}^{4}b{c}^{2}{d}^{3}+16\,{a}^{4}bc{d}^{4}+4\,{a}^{4}{c}^{2}{d}^{4}-10\,{a}^{3}{b}^{3}{c}^{3}d+16\,{a}^{3}{b}^{3}{c}^{2}{d}^{2}\\
    & \qquad -10\,{a}^{3}{b}^{3}c{d}^{3}-16\,{a}^{3}{b}^{2}{c}^{3}{d}^{2}-24\,{a}^{3}{b}^{2}{c}^{2}{d}^{3}+4\,{a}^{3}{b}^{2}c{d}^{4}-10\,{a}^{3}b{c}^{3}{d}^{3}-4\,{a}^{3}b{c}^{2}{d}^{4}\\
    & \qquad +4\,{a}^{2}{b}^{4}{c}^{4}-4\,{a}^{2}{b}^{4}{c}^{3}d+10\,{a}^{2}{b}^{4}{c}^{2}{d}^{2}+4\,{a}^{2}{b}^{3}{c}^{4}d-24\,{a}^{2}{b}^{3}{c}^{3}{d}^{2}-16\,{a}^{2}{b}^{3}{c}^{2}{d}^{3}\\
    & \qquad+10\,{a}^{2}{b}^{2}{c}^{4}{d}^{2}+16\,{a}^{2}{b}^{2}{c}^{3}{d}^{3}+10\,{a}^{2}{b}^{2}{c}^{2}{d}^{4}+16\,a{b}^{4}{c}^{4}d+4\,a{b}^{4}{c}^{3}{d}^{2}-4\,a{b}^{3}{c}^{4}{d}^{2}\\
    & \qquad -10\,a{b}^{3}{c}^{3}{d}^{3}+4\,{b}^{4}{c}^{4}{d}^{2}-10\,{a}^{4}{b}^{2}{c}^{2}d-4\,{a}^{4}{b}^{2}c{d}^{2}-6\,{a}^{4}{b}^{2}{d}^{3}-16\,{a}^{4}b{c}^{2}{d}^{2}\\
    & \qquad -32\,{a}^{4}bc{d}^{3}-12\,{a}^{4}b{d}^{4}-10\,{a}^{4}{c}^{2}{d}^{3}-12\,{a}^{4}c{d}^{4}+5\,{a}^{3}{b}^{3}{c}^{3}-{a}^{3}{b}^{3}{c}^{2}d-{a}^{3}{b}^{3}c{d}^{2}\\
    & \qquad +5\,{a}^{3}{b}^{3}{d}^{3}+31\,{a}^{3}{b}^{2}{c}^{3}d+16\,{a}^{3}{b}^{2}{c}^{2}{d}^{2}+39\,{a}^{3}{b}^{2}c{d}^{3}-10\,{a}^{3}{b}^{2}{
      d}^{4}+31\,{a}^{3}b{c}^{3}{d}^{2}\\
    & \qquad +39\,{a}^{3}b{c}^{2}{d}^{3}-32\,{a}^{3}bc{d}^{4}+5\,{a}^{3}{c}^{3}{d}^{3}-6\,{a}^{3}{c}^{2}{d}^{4}-6\,{a}^{2}{b}^{4}{c}^{3}-4\,{a}^{2}{b}^{4}{c}^{2}d-10\,{a}^{2}{b}^{4}c{d}^{2}\\
    & \qquad -10\,{a}^{2}{b}^{3}{c}^{4}+39\,{a}^{2}{b}^{3}{c}^{3}d+16\,{a}^{2}{b}^{3
    }{c}^{2}{d}^{2}+31\,{a}^{2}{b}^{3}c{d}^{3}-16\,{a}^{2}{b}^{2}{c}^{4}d\\
    & \qquad +16\,{a}^{2}{b}^{2}{c}^{3}{d}^{2}+16\,{a}^{2}{b}^{2}{c}^{2}{d}^{3}-16\,{a}^{2}{b}^{2}c{d}^{4}-10\,{a}^{2}b{c}^{4}{d}^{2}-{a}^{2}b{c}^{3}{d}^{3}-4\,{a}^{2}b{c}^{2}{d}^{4}\\
    & \qquad -12\,a{b}^{4}{c}^{4}-32\,a{b}^{4}{c}^{3}d-16\,a{b}^{4}{c}^{2}{d}^{2}-32\,a{b}^{3}{c}^{4}d+39\,a{b}^{3}{c}^{3}{d}^{2}+31\,a{b}^{3}{c}^{2}{d}^{3}\\
    & \qquad -4\,a{b}^{2}{c}^{4}{d}^{2}-a{b}^{2}{c}^{3}{d}^{3}-10\,a{b}^{2}{c}^{2}{d}^{4}-12\,{b}^{4}{c}^{4}d-10\,{b}^{4}{c}^{3}{d}^{2}-6\,{b}^{3}{c}^{4}{d}^{2}+5\,{b}^{3}{c}^{3}{d}^{3}\\
    & \qquad +4\,{a}^{4}{b}^{2}{c}^{2}+16\,{a}^{4}{b}^{2}cd+4\,{a}^{4}{b}^{2}{d}^{2}+4\,{a}^{4}b{c}^{2}d+8\,{a}^{4}bc{d}^{2}+12\,{a}^{4}b{d}^{3}+10\,{a}^{4}{c}^{2}{d}^{2}\\
    & \qquad +36\,{a}^{4}c{d}^{3}+18\,{a}^{4}{d}^{4}-5\,{a}^{3}{b}^{3}{c}^{2}+16\,{a}^{3}{b}^{3}cd-5\,{a}^{3}{b}^{3}{d}^{2}-10\,{a}^{3}{b}^{2}{c}^{3}-28\,{a}^{3}{b}^{2}{c}^{2}d\\
    & \qquad -56\,{a}^{3}{b}^{2}c{d}^{2}-6\,{a}^{3}{b}^{2}{d}^{3}-16\,{a}^{3}b{c}^{3}d-21\,{a}^{3}b{c}^{2}{d}^{2}+40\,{a}^{3}bc{d}^{3}+36\,{a}^{3}b{d}^{4}\\
    & \qquad -10\,{a}^{3}{c}^{3}{d}^{2}-6\,{a}^{3}{c}^{2}{d}^{3}+12\,{a}^{3}c{d}^{4}+4\,{a}^{2}{b}^{4}{c}^{2}+16\,{a}^{2}{b}^{4}cd+4\,{a}^{2}{b}^{4}{d}^{2}-6\,{a}^{2}{b}^{3}{c}^{3}\\
    & \qquad -56\,{a}^{2}{b}^{3}{c}^{2}d-28\,{a}^{2}{b}^{3}c{d}^{2}-10\,{a}^{2}{b}^{3}{d}^{3}+10\,{a}^{2}{b}^{2}{c}^{4}-21\,{a}^{2}{b}^{2}{c}^{3}d+12\,{a}^{2}{b}^{2}{c}^{2}{d}^{2}\\
    & \qquad -21\,{a}^{2}{b}^{2}c{d}^{3}+10\,{a}^{2}{b}^{2}{d}^{4}+4\,{a}^{2}b{c}^{4}d-28\,{a}^{2}b{c}^{3}{d}^{2}-56\,{a}^{2}b{c}^{2}{d}^{3}+8\,{a}^{2}bc{d}^{4}\\
    & \qquad +4\,{a}^{2}{c}^{4}{d}^{2}-5\,{a}^{2}{c}^{3}{d}
    ^{3}+4\,{a}^{2}{c}^{2}{d}^{4}+12\,a{b}^{4}{c}^{3}+8\,a{b}^{4}{c}^{2}d+
    4\,a{b}^{4}c{d}^{2}+36\,a{b}^{3}{c}^{4}\\
    & \qquad +40\,a{b}^{3}{c}^{3}d-21\,a{b}^
    {3}{c}^{2}{d}^{2}-16\,a{b}^{3}c{d}^{3}+8\,a{b}^{2}{c}^{4}d-56\,a{b}^{2
    }{c}^{3}{d}^{2}-28\,a{b}^{2}{c}^{2}{d}^{3}\\
    & \qquad +4\,a{b}^{2}c{d}^{4}+16\,ab{c}^{4}{d}^{2}+16\,ab{c}^{3}{d}^{3}+16\,ab{c}^{2}{d}^{4}+18\,{b}^{4}{c}^{4}+36\,{b}^{4}{c}^{3}d+10\,{b}^{4}{c}^{2}{d}^{2}\\
    & \qquad +12\,{b}^{3}{c}^{4}d-6\,{b}^{3}{c}^{3}{d}^{2}-10\,{b}^{3}{c}^{2}{d}^{3}+4\,{b}^{2}{c}^{4}{
      d}^{2}-5\,{b}^{2}{c}^{3}{d}^{3}+4\,{b}^{2}{c}^{2}{d}^{4}-8\,{a}^{4}{b}
    ^{2}c\\
    & \qquad -6\,{a}^{4}{b}^{2}d+8\,{a}^{4}bcd+12\,{a}^{4}b{d}^{2}-24\,{a}^{4}c{d}^{2}-36\,{a}^{4}{d}^{3}-5\,{a}^{3}{b}^{3}c-5\,{a}^{3}{b}^{3}d\\
    & \qquad +10\,{a}^{3}{b}^{2}{c}^{2}+3\,{a}^{3}{b}^{2}cd+32\,{a}^{3}{b}^{2}{d}^{2}+8\,{a}^{3}b{c}^{2}d+20\,{a}^{3}bc{d}^{2}-36\,{a}^{3}b{d}^{3}+6\,{a}^{3}{c}^{2}{d}^{2}\\
    & \qquad -36\,{a}^{3}c{d}^{3}-36\,{a}^{3}{d}^{4}-6\,{a}^{2}{b}^{4}c-8\,{a}^{2}{b}^{4}d+32\,{a}^{2}{b}^{3}{c}^{2}+3\,{a}^{2}{b}^{3}cd+10\,{a}^{2}{b}^{3}{d}^{2}\\
    & \qquad +6\,{a}^{2}{b}^{2}{c}^{3}+16\,{a}^{2}{b}^{2}{c}^{2}d+16\,{a}^{2}{b}^{2}c{d}^{2}+6\,{a}^{2}{b}^{2}{d}^{3}+8\,{a}^{2}b{c}^{3}d+16\,{a}^{2}b{c}^{2}{d}^{2}\\
    & \qquad +20\,{a}^{2}bc{d}^{3}-24\,{a}^{2}b{d}^{4}+10\,{a}^{2}{c}^{3}{d}^{2}+32\,{a}^{2}{c}^{2}{d}^{3}+12\,{a}^{2}c{d}^{4}+12\,a{b}^{4}{c}^{2}\\
    & \qquad +8\,a{b}^{4}cd-36\,a{b}^{3}{c}^{3}+20\,a{b}^{3}{c}^{2}d+8\,a{b}^{3}c{d}^{2}-24\,a{b}^{2}{c}^{4}+20\,a{b}^{2}{c}^{3}d+16\,a{b}^{2}{c}^{2}{d}^{2}\\
    & \qquad +8\,a{b}^{2}c{d}^{3}+8\,ab{c}^{4}d+3\,ab{c}^{3}{d}^{2}+3\,ab{c}^{2}{d}^{3}+8\,abc{d}^{4}-8\,a{c}^{4}{d}^{2}-5\,a{c}^{3}{d}^{3}-6\,a{c}^{2}{d}^{4}\\
    & \qquad -36\,{b}^{4}{c}^{3}-24\,{b}^{4}{c}^{2}d-36\,{b}^{3}{c}^{4}-36\,{b}^{3}{c}^{3}d+6\,{b}^{3}{c}^{2}{d}^{2}+12\,{b}^{2}{c}^{4}d+32\,{b}^{2}{c}^{3}{d}^{2}\\
    & \qquad +10\,{b}^{2}{c}^{2}{d}^{3}-6\,b{c}^{4}{d}^{2}-5\,b{c}^{3}{d}^{3}-8\,b{c}^{2}{d}^{4}+4\,{a}^{4}{b}^{2}-12\,{a}^{4}bd+18\,{a}^{4}{d}^{2}+5\,{a}^{3}{b}^{3}\\
    & \qquad +10\,{a}^{3}{b}^{2}c-6\,{a}^{3}{b}^{2}d-28\,{a}^{3}bcd-36\,{a}^{3}b{d}^{2}+24\,{a}^{3}c{d}^{2}+72\,{a}^{3}{d}^{3}+4\,{a}^{2}{b}^{4}-6\,{a}^{2}{b}^{3}c\\
    & \qquad +10\,{a}^{2}{b}^{3}d-32\,{a}^{2}{b}^{2}{c}^{2}+15\,{a}^{2}{b}^{2}cd-32\,{a}^{2}{b}^{2}{d}^{2}+32\,{a}^{2}b{c}^{2}d-8\,{a}^{2}bc{d}^{2}+24\,{a}^{2}b{d}^{3}\\
    & \qquad -32\,{a}^{2}{c}^{2}{d}^{2}-36\,{a}^{2}c{d}^{3}+18\,{a}^{2}{d}^{4}-12\,a{b}^{4}c-36\,a{b}^{3}{c}^{2}-28\,a{b}^{3}cd+24\,a{b}^{2}{c}^{3}\\
    & \qquad -8\,a{b}^{2}{c}^{2}d+32\,a{b}^{2}c{d}^{2}-28\,ab{c}^{3}d+15\,ab{c
    }^{2}{d}^{2}-28\,abc{d}^{3}+10\,a{c}^{3}{d}^{2}-6\,a{c}^{2}{d}^{3}\\
    & \qquad -12\,ac{d}^{4}+18\,{b}^{4}{c}^{2}+72\,{b}^{3}{c}^{3}+24\,{b}^{3}{c}^{2}d+18\,{b}^{2}{c}^{4}-36\,{b}^{2}{c}^{3}d-32\,{b}^{2}{c}^{2}{d}^{2}\\
    & \qquad -12\,b{c}^{4}d-6\,b{c}^{3}{d}^{2}+10\,b{c}^{2}{d}^{3}+4\,{c}^{4}{d}^{2}+5\,{c}^{3}{d}^{3}+4\,{c}^{2}{d}^{4}-10\,{a}^{3}{b}^{2}+36\,{a}^{3}bd\\
    & \qquad -36\,{a}^{3}{d}^{2}-10\,{a}^{2}{b}^{3}+6\,{a}^{2}{b}^{2}c+6\,{a}^{2}{b}^{2}d-20\,{a}^{2}bcd+24\,{a}^{2}b{d}^{2}+24\,{a}^{2}c{d}^{2}\\
    & \qquad -36\,{a}^{2}{d}^{3}+36\,a{b}^{3}c+24\,a{b}^{2}{c}^{2}-20\,a{b}^{2}cd-20\,ab{c}^{2}d-
    20\,abc{d}^{2}+6\,a{c}^{2}{d}^{2}\\
    & \qquad +36\,ac{d}^{3}-36\,{b}^{3}{c}^{2}-36
    \,{b}^{2}{c}^{3}+24\,{b}^{2}{c}^{2}d+36\,b{c}^{3}d+6\,b{c}^{2}{d}^{2}-10\,{c}^{3}{d}^{2}-10\,{c}^{2}{d}^{3}\\
    & \qquad +10\,{a}^{2}{b}^{2}-24\,{a}^{2}bd+18\,{a}^{2}{d}^{2}-24\,a{b}^{2}c+40\,abcd-24\,ac{d}^{2}+18\,{b}^{2}{c}^{2}-24\,b{c}^{2}d\\
    & \qquad +10\,{c}^{2}{d}^{2} ) {s}^{4}-648\, ( 2\,{a}^{3}{b}^{2}{c}^{2}d+{a}^{3}{b}^{2}c{d}^{2}-{a}^{3}{b}^{2}{d}^{3}-{a}^{3}b{c}^{2}{d}^{2}-4\,{a}^{3}bc{d}^{3}-{a}^{3}{c}^{2}{d}^{3}\\
    & \qquad -{a}^{2}{b}^{3}{c}^{3}+{a}^{2}{b}^{3}{c}^{2}d+2\,{a}^{2}{b}^{3}c{d}^{2}-{a}^{2}{b}^{2}{c}^{3}d+4\,{a}^{2}{b}^{2}{c}^{2}{d}^{2}-{a}^{2}{b}^{2}c{d}^{3}+2\,{a}^{2}b{c}^{3}{d}^{2}\\
    & \qquad +{a}^{2}b{c}^{2}{d}^{3}-4\,a{b}^{3}{c}^{3}d-a{b}^{3}{c}^{2}{d}^{2}+a{b}^{2}{c}^{3}{d}^{2}+2\,a{b}^{2}{c}^{2}{d}^{3}-{b}^{3}{c}^{3}{d}^{2}-{a}^{3}{b}^{2}{c}^{2}\\
    & \qquad -3\,{a}^{3}{b}^{2}cd+{a}^{3}{b}^{2}{d}^{2}-{a}^{3}b{c}^{2}d+6\,{a}^{3}bc{d}^{2}+3\,{a}^{3}b{d}^{3}+2\,{a}^{3}{c}^{2}{d}^{2}+3\,{a}^{3}c{d}^{3}+{a}^{2}{b}^{3}{c}^{2}\\
    & \qquad -3\,{a}^{2}{b}^{3}cd-{a}^{2}{b}^{3}{d}^{2}+2\,{a}^{2}{b}^{2}{c}^{3}-7\,{a}^{2}{b}^{2}{c}^{2}d-7\,{a}^{2}{b}^{2}c{d}^{2}+2\,{a}^{2}{b}^{2}{d}^{3}-{a}^{2}b{c}^{3}d\\
    & \qquad -7\,{a}^{2}b{c}^{2}{d}^{2}+6\,{a}^{2}bc{d}^{3}-{a}^{2}{c}^{3}{d}^{2}+{a}^{2}{c}^{2}{d}^{3}+3\,a{b}^{3}{c}^{3}+6\,a{b}^{3}{c}^{2}d-a{b}^{3}c{d}^{2}+6\,a{b}^{2}{c}^{3}d\\
    & \qquad -7\,a{b}^{2}{c}^{2}{d}^{2}-a{b}^{2}c{d}^{3}-3\,ab{c}^{3}{d}^{2}-3\,ab{c}^{2}{d}^{3}+3\,{b}^{3}{c}^{3}d+2\,{b}^{3}{c}^{2}{d}^{2}+{b}^{2}{c}^{3}{d}^{2}-{b}^{2}{c}^{2}{d}^{3}\\
    & \qquad +2\,{a}^{3}{b}^{2}c+{a}^{3}{b}^{2}d-2\,{a}^{3}bcd-6\,{a}^{3}b{d}^{2}-3\,{a}^{3}c{d}^{2}+{a}^{2}{b}^{3}c+2\,{a}^{2}{b}^{3}d-2\,{a}^{2}{b}^{2}{c}^{2}\\
    & \qquad +10\,{a}^{2}{b}^{2}cd-2\,{a}^{2}{b}^{2}{d}^{2}+10\,{a}^{2}b{c}^{2}d-3\,{a}^{2}b{d}^{3}-2\,{a}^{2}{c}^{2}{d}^{2}-6\,{a}^{2}c{d}^{3}-6\,a{b}^{3}{c}^{2}\\
    & \qquad -2\,a{b}^{3}cd-3\,a{b}^{2}{c}^{3}+10\,a{b}^{2}c{d}^{2}-2\,ab{c}^{3}d+10\,ab{c}^{2}{d}^{2}-2\,abc{d}^{3}+2\,a{c}^{3}{d}^{2}+a{c}^{2}{d}^{3}\\
    & \qquad -3\,{b}^{3}{c}^{2}d-6\,{b}^{2}{c}^{3}d-2\,{b}^{2}{c}^{2}{d}^{2}+b{c}^{3}{d}^{2}+2\,b{c}^{2}{d}^{3}-{a}^{3}{b}^{2}+3\,{a}^{3}bd-{a}^{2}{b}^{3}-2\,{a}^{2}{b}^{2}c\\
    & \qquad -2\,{a}^{2}{b}^{2}d-6\,{a}^{2}bcd+6\,{a}^{2}b{d}^{2}+6\,{a}^{2}c{d}^{2}+3\,a{b}^{3}c+6\,a{b}^{2}{c}^{2}-6\,a{b}^{2}cd-6\,ab{c}^{2}d\\
    & \qquad -6\,abc{d}^{2}-2\,a{c}^{2}{d}^{2}+3\,ac{d}^{3}+6\,{b}^{2}{c}^{2}d+3\,b{c}^{3}d-2\,b{c}^{2}{d}^{2}-{c}^{3}{d}^{2}-{c}^{2}{d}^{3}+2\,{a}^{2}{b}^{2}\\
    & \qquad -3\,{a}^{2}bd-3\,a{b}^{2}c+8\,abcd-3\,ac{d}^{2}-3\,b{c}^{2}d+2\,{c}^{2}{d}^{2} )  ( abc+abd-acd-bcd\\
    & \qquad -ab+cd ) ^{2}{s}^{2}+54\, ( abc+abd-acd-bcd-ab+cd ) ^{6} \Big).
\end{align*}

\subsection{Class 2.1 : $\II^*+6\I_2+2\I_1$ with MWG = $\{0\}$}
\begin{equation*}
    \begin{aligned}
        y^2 & + \Big(2s^2+4b-4ac-2ad-2-2d+2c+4a-2cb \Big)xy- \Big( 4c(a-1)s^4 -4(3a^2c^2-a^2d\\
        & -2acd-c+a^2+c^2+abcd+ad^2-bc^2+5ac+2a^2cd-4a^2c-3abc-2bcd-a+2abc^2\\
        & +cd-4ac^2+2bc)s^2+4a(c+d-1+)(1-3a+a^2cd+abcd+6ac+ad^2-a^2d-3ac^2\\
        & +2cd+3abc^2-2bcd-4a^2c+2a^2c^2-5abc+b^2c^2-2bad+2ad+2bc-2d-2c-2acd\\
        & +2a^2-2bc^2+c^2+d^2+2ab) \Big)y\\
        & = \ x^3 + \Big(2c(1-a)s^2-2a+2a^2cd+2abcd+4ac+2ad^2-2a^2d-2ac^2+4cd+4b^2c+2abc^2\\
        & -4bcd-4a^2c+2a^2c^2-2abc-4bad+4bd+4ad-4bc-4b^2-4d+4b-4acd+2a^2 \Big)x^2\\
        & - \Big( 4c^2(a-1)^2s^4 -8ac(c+d-1)(a-1)(ac-a-c+cb-d+1)s^2+4a^2(c+d-1)^2(ac\\
        & -a-c+cb-d+1)^2 \Big)x -8 \Big(c(1-a)s^2-a+a^2cd+abcd+2ac+ad^2-a^2d-ac^2+2cd\\
        & +2b^2c+abc^2-2bcd-2a^2c+a^2c^2-abc-2bad+2bd+2ad-2bc-2b^2-2d+2b-2acd\\
        & +a^2 \Big) \Big(c^2(a-1)^2s^4 -2ac(c+d-1)(a-1)(ac-a-c+cb-d+1)s^2 + a^2(c+d-1)^2(ac\\
        & -a-c+cb-d+1)^2 \Big).
    \end{aligned}
\end{equation*}

\subsection{Class 2.2 : $\III^*+7\I_2+\I_1$ with MWG= $\mathbb{Z}/2\mathbb{Z}$}
\begin{equation*}
    \begin{aligned}
        y^2 = & \ x^3 -9 \Big( (1-2a-2b+c+bc+d+ad)s^2 + (4bc^2b+c^2b^2d-ab^2c^2+ab^2c-4b^2c^2\\
        & +2b^2c^3-abcd^2-4abcd-abc^2d+a^2bcd+ab^2cd+cd^2-cd+2bc^3+2ad-2bc^2+2b^2c\\
        & -2b^2cd-bcd+bcd^2-ac^2d-acd+2a^2cd-a^2cd^2+2a^2d^2-2ad^2-2a^2d+c^2d-3abc^2\\
        & +4abc-ab+abd^2+a^2bd-a^2bd^2)s -bc (c+d-1)(a-c) (ad-c-bc-d+1)(ad-bc\\
        & +b-d) \Big)x^2\\
        & +81s(s+acd+ad^2-bc^2-bcd-ad+bc)(s+abd+acd-b^2c-bc^2-ad+b^2+2bc\\
        & -bd-cd-b+d+)(s+abcd-bc^2-bcd+bc)(s-a+a^2+c-ac+abc-bc^2+ad\\
        & -a^2d-cd+acd) x.
    \end{aligned}
\end{equation*}

\subsection{Class 2.3 : $\III^*+\I_0^*+3\I_2+3\I_1$ with MWG = $\{0\}$}
\begin{align*}
    y^2 = & \ x^3 +27s^2 \Big( 3(a-1)s^3 + (a+d+b+c+abcd-a^2d^2-b^2c^2-a^2-b^2-c^2-d^2-ad-bc\\
    & +a^2d+7abc+ad^2+bc^2+b^2c+7cd+bd-8ab+ac+7abd-8acd-8bcd-1)s^2+(7bc^2d^2\\
    & -b^3c^2d-2acd^2+9bc^2d+b^3c^2d-2a^2cd^2+3a^2bd^2+2a^2cd^3-2a^2bd^3+7ac^2d^2-3b^2c^2d^2\\
    & -5b^2c^2d-3bc^3d-acd^3-3abc^2d^2+3ab^2c^2d+3abcd^3-a^2bcd^2-8abcd^2+2ab^2cd+3ab^2cd^2\\
    & -8abc^2d+5abcd-a^2d^3+3a^2d^2+3b^3c^2+3b^2c^3-2b^3c^3-3b^2c^2-4c^2d^2+ad+bc-2a^2d\\
    & +4abc-2ac^2d-2abd^3+4ab^2d^2+4ab^2c^2-b^2cd^2+2a^2cd-a^2bd-7ab^2c-abc^2-7ab^2d\\
    & -3ad^2+2ad^3+cd^2+cd-3c^2d-b^3c+2c^3d-bc^3-a^2bcd+3ab^2-3ab-abd^2+6abd+acd\\
    & -7bcd+7b^2cd-b^3cd)s + (a^2bc^2d^3-a^2b^2c^2d^2-a^2b^2cd^3-a^2b^2d^4-a^2bcd^4-a^2c^2d^4\\
    & +2ab^3c^3d+2ab^3c^2d^2+2ab^3cd^3-2ab^2c^3d^2+2ab^2c^2d^3+2abc^3d^3-b^4c^4-b^4c^3d-b^4c^2d^2\\
    & +b^3c^4d-b^3c^3d^2-b^2c^4d^2+a^2b^2cd^2+2a^2b^2d^3+a^2bc^2d^2+4a^2bcd^3+2a^2bd^4+a^2c^2d^3+a^2cd^4\\
    & -3ab^3c^2d-3ab^3cd^2-3ab^2c^3d-4ab^2c^2d^2-5ab^2cd^3+abc^3d^2-5abc^2d^3-2ac^3d^3+2b^4c^3\\
    & +b^4c^2d+2b^3c^4+3b^3c^2d^2-2b^2c^4d+4b^2c^3d^2+2bc^4d^2-a^2b^2d^2-3a^2bcd^2-4a^2bd^3-a^2c^2d^2\\
    & -3a^2cd^3-a^2d^4+ab^3cd+5ab^2c^2d+7ab^2cd^2+abc^3d+4abc^2d^2+4abcd^3+ac^3d^2+3ac^2d^3\\
    & -b^4c^2-4b^3c^3-b^3c^2d-b^2c^4+3b^2c^3d-4b^2c^2d^2+bc^4d-5bc^3d^2-c^4d^2+2a^2bd^2+2a^2cd^2\\
    & +2a^2d^3-2ab^2cd-2abc^2d-5abcd^2-2ac^2d^2-acd^3+2b^3c^2+2b^2c^3-b^2c^2d-2bc^3d+3bc^2d^2\\
    & +2c^3d^2-a^2d^2+abcd+acd^2-b^2c^2+bc^2d-c^2d^2) \Big)x \\
    &+ 27s^3 \Big( (18a^2-9a^2d-9abc+18ab-9ac+27ad-18bc-27a+9b+9c-18d+9)s^4 +(2a^3d^3\\
    & -3a^2bcd^2-3ab^2c^2d+2b^3c ^3-3a^3d^2+3a^2bcd-30a^2bd^2+33a^2cd^2-3a^2d^3-30ab^2c^2+3ab^2cd\\
    &+3abc^2d+66abcd^2-3b^3c^2-3b^2c^3-30b^2c^2d-3a^3d+33a^2bc+57a^2bd-60a^2cd-27a^2d^2\\
    & +57ab^2c+33ab^2d+3abc^2-45abcd+3abd^2+33ac^2d-69acd^2-3ad^3-3b^3c-27b^2c^2+3b^2cd\\
    & -3bc^3+57bc^2d-30bcd^2+2a^3-30a^2b-3a^2c+36a^2d-30ab^2-75abc-75abd-3ac^2\\
    & +51acd+36ad^2+2b^3+36b^2c-3b^2d+36bc^2-12bcd-3bd^2+2c^3-30c^2d+33cd^2+2d^3\\
    & -3a^2+66ab+12ac-27ad-3b^2-27bc+12bd-3c^2+3cd-3d^2-3a-3b-3c-3d+2)s^3\\
    & +(3a^3bcd^3+6a^3bd^4-6a^3cd^4+3a^2b^2c^2d^2-12a^2b^2cd^3+12a^2bc^2d^3-9a^2bcd^4-12ab^3c^3d\\
    & +3ab^3c^2d^2-3ab^2c^3d^2+27ab^2c^2d^3+6b^4c^4+3b^4c^3d-3b^3c^4d-18b^3c^3d^2-12a^3bcd^2\\
    & -12a^3bd^3+9a^3cd^3+3a^3d^4+27a^2b^2c^2d-3a^2b^2cd^2-33a^2b^2d^3-3a^2bc^2d^2+3a^2bcd^3\\
    & +3a^2bd^4-45a^2c^2d^3+6a^2cd^4-33ab^3c^3-3ab^3c^2d+27ab^3cd^2+39ab^2c^3d+102ab^2c^2d^2\\
    & +45ab^2cd^3+33abc^3d^2-99abc^2d^3-9abcd^4-12b^4c^3-12b^4c^2d-12b^3c^4-42b^3c^3d\\
    & -3b^3c^2d^2-3b^2c^4d+57b^2c^3d^2-18b^2c^2d^3+3a^3bcd+3a^3bd^2+9a^3cd^2+3a^3d^3+15a^2b^2c^2\\
    & +39a^2b^2cd+72a^2b^2d^2-93a^2bc^2d-69a^2bcd^2+9a^2bd^3+72a^2c^2d^2+66a^2cd^3+6a^2d^4\\
    & +72ab^3c^2+39ab^3cd+15ab^3d^2+27ab^2c^3-165ab^2c^2d-252ab^2cd^2-15ab^2d^3-3abc^3d\\
    & -42abc^2d^2+36abcd^3+6abd^4-45ac^3d^2+57ac^2d^3+3acd^4+3b^4c^2+3b^4cd+6b^3c^3\\
    & +54b^3c^2d-12b^3cd^2+3b^2c^4+102b^2c^3d-3b^2c^2d^2+3b^2cd^3+12bc^4d-72bc^3d^2+33bc^2d^3\\
    & +3a^3bd-6a^3cd-12a^3d^2-33a^2b^2c-57a^2b^2d+3a^2bc^2+105a^2bcd+6a^2bd^2+9a^2c^2d\\
    & -78a^2cd^2-30a^2d^3-57ab^3c-33ab^3d-90ab^2c^2+90ab^2cd+12ab^2d^2-12abc^3+132abc^2d\\
    & +165abcd^2+3abd^3+9ac^3d-6ac^2d^2-54acd^3-6ad^4+3b^4c+21b^3c^2+3b^3cd+21b^2c^3\\
    & -102b^2c^2d+18b^2cd^2+3bc^4-75bc^3d+33bc^2d^2-6c^4d+33c^3d^2-15c^2d^3+6a^3d+18a^2b^2\\
    & -12a^2bc-27a^2bd-3a^2cd+30a^2d^2+18ab^3+90ab^2c+30ab^2d+45abc^2-123abcd-24abd^2\\
    & -33ac^2d+15acd^2+12ad^3-15b^3c-39b^2c^2-15bc^3+63bc^2d-18bcd^2+12c^3d-27c^2d^2\\
    & -3cd^3+9a^2b-6a^2d-27ab^2-42abc+6abd+27acd-3ad^2+15b^2c+15bc^2-3bcd-3c^2d\\
    & +12cd^2+9ab-3ad-3bc-3cd)s^2 - (abcd+2abd^2+acd^2-b^2c^2-2b^2cd-bc^2d-2abd\\
    & -2acd-2ad^2+b^2c+bc^2+3bcd+c^2d+2ad-bc-cd) (2abcd+abd^2-acd^2-2b^2c^2-b^2cd\\
    & +bc^2d-abd-acd-ad^2+2b^2c+2bc^2-c^2d+ad-2bc+cd) (abcd-abd^2-2acd^2-b^2c^2\\
    & +b^2cd+2bc^2d+abd+acd+ad^2+b^2c+bc^2-3bcd-2c^2d-ad-bc+2cd) \Big).
\end{align*}

\subsection{Class 2.4 : $\I_6^* + 4\I_2+4\I_1$ with MWG = $\{ 0 \}$}
\begin{align*}
 y^2 = x^3 &+ s\Big( s^4 -(3\,ac+ad+bc-2\,a-2\,b-2\,c+d+1) s^3 + (3\,{a}^{2}{c}^{2}+2\,{a}^{2}cd+2\,ab{c}^{2}\\
 & \quad +abcd-4\,{a}^{2}c-{a}^{2}d-3\,abc+abd-4\,a{c}^{2}-2\,acd+a{d}^{2}-{b}^{2}c-b{c}^{2}-2\,bcd+{a}^{2}\\
 & \quad +5\,ac-ad+{b}^{2}+3\,bc-bd+{c}^{2}+cd-a-b-c+d) s^2 - ({a}^{3}{c}^{3}+{a}^{3}{c}^{2}d+{a}^{2}b{c}^{3}\\
 & \quad +{a}^{2}b{c}^{2}d-2\,{a}^{3}{c}^{2}-{a}^{3}cd-{a}^{2}b{c}^{2}+4\,{a}^{2}bcd-2\,{a}^{2}{c}^{3}-3\,{a}^{2}{c}^{2}d-{a}^{2}c{d}^{2}
 -4\,a{b}^{2}{c}^{2}\\
 & \quad -ab{c}^{3}-ab{c}^{2}d+2\,abc{d}^{2}-2\,{b}^{2}{c}^{2}d+{a}^{3}c+4\,{a}^{2}{c}^{2}-{a}^{2}cd-2\,{a}^{2}{d}^{2}+4\,a{b}^{2}c+5\,ab{c}^{2}\\
 & \quad -4\,abcd+a{c}^{3}+2\,a{c}^{2}d-ac{d}^{2}+2\,{b}^{2}{c}^{2}+2\,{b}^{2}cd+2\,b{c}^{2}d-2\,bc{d}^{2}-2\,{a}^{2}c+2\,{a}^{2}d\\
 & \quad -6\,abc-2\,abd-2\,a{c}^{2}+4\,acd+2\,a{d}^{2}-2\,{b}^{2}c-2\,b{c}^{2}+2\,c{d}^{2}+2\,ab+ac-2\,ad\\
 & \quad +2\,bc-2\,cd) s +3ac(bc-c-d+1)(a+d-1)(ad-bc+b-d) \Big) x^2\\
 & +s^2 (s-ac-cd+c)(s+abc-ac-ad-bs+a) \Big( (a-1)(ad-bc)(c+d-1)s^3\\
 & \quad -(ad-bc+b-d)(3\,{a}^{2}{c}^{2}+2\,{a}^{2}cd+2\,ab{c}^{2}+abcd-4\,{a}^{2}c-{a}^{2}d-3\,abc-4\,a{c}^{2}\\
 & \quad -2\,acd+a{d}^{2}-b{c}^{2}-2\,bcd+{a}^{2}+5\,ac+2\,bc+{c}^{2}+cd-a-c)s^2 + (ad-bc+b\\
 & \quad -d)(2\,{a}^{3}{c}^{3}+2\,{a}^{3}{c}^{2}d+2\,{a}^{2}b{c}^{3}+2\,{a}^{2}b{c}^{2}d-4\,{a}^{3}{c}^{2}-2\,{a}^{3}cd-2\,{a}^{2}b{c}^{2}
 +2\,{a}^{2}bcd\\
 & \quad -4\,{a}^{2}{c}^{3}-6\,{a}^{2}{c}^{2}d-2\,{a}^{2}c{d}^{2}-2\,a{b}^{2}{c}^{2}-2\,ab{c}^{3}-2\,ab{c}^{2}d+abc{d}^{2}-{b}^{2}{c}^{2}d
 +2\,{a}^{3}c\\
 & \quad +8\,{a}^{2}{c}^{2}+4\,{a}^{2}cd-{a}^{2}{d}^{2}+2\,a{b}^{2}c+4\,ab{c}^{2}-2\,abcd+2\,a{c}^{3}+4\,a{c}^{2}d+ac{d}^{2}+{b}^{2}{c}^{2}\\
 & \quad +{b}^{2}cd+b{c}^{2}d-bc{d}^{2}-4\,{a}^{2}c+{a}^{2}d-3\,abc-abd-4\,a{c}^{2}-acd+a{d}^{2}-{b}^{2}c-b{c}^{2}\\
 & \quad +c{d}^{2}+ab+2\,ac-ad+bc-cd)-3\, ac (bc-c-d+1)(a+d-1)(ad-bc+b-d)^{2} \Big) x\\ 
 & +ac(ad-bc+b-d)^2 s^3 (s-ac-cd+c)^2 ((1-b)s+abc-ac-ad+a)^2 \Big( (1-a)(c+d\\
 & \quad -1)(ac+bc-a-c-d+1)s+(bc-c-d+1)(a+d-1)(ad-bc+b-d) \Big). 
\end{align*}

\subsection{Class 2.5 : $\I_4^* + 6\I_2+2\I_1$ with MWG = $\mathbb{Z}/2\mathbb{Z}$}
\begin{equation*}
    \begin{aligned}
        y^2 &= x^3 + t \, \Big( ac(a-1)(c+d-1)(ac+bc-a-c-d+1) t^3 +(4ac^2-c^2+a-5ac\\
        & -cd-a^2+2acd+2bcd-ad^2-3a^2c^2+bc^2-2a^2cd+c+3abc-abcd+a^2d\\
        & -2bc+4a^2c-2abc^2) t^2+(ad+bc+3ac-2a-2b-2c+d+1) t -1 \Big) x^2\\
        & -t^4 \, \Big( c(a+d-1) t-1 \Big) \Big( a(bc-c-d+1) t -b+1 \Big) \Big( (a-1)(ad-bc)(c+d-1)t\\
        & -ad+bc-b+d \Big) x.
    \end{aligned}
\end{equation*}

\subsection{Class 2.6 : $\I_4^* + \I_0^* + 2\I_2+4\I_1$ with MWG = $\{ 0 \}$}
\begin{align*}
  y^2 &= x^3 -27 s^2 \Big(s^4 +( ad+bc+a+b-2\,c-2\,d+1
  ) {s}^{3}+ ( {a}^{2}{d}^{2}-abcd +{b}^{2}{c}^{2}+4\,{a}^{2}d+3\,abc\\
  & \quad   +3\,abd-2\,acd-a{d}^{2}+4\,{b}^{2}c-6\,b{c}^{2}-7\,bcd+{a}^{2}+3\,ab-6\,ac-9\,ad+{b}^{2}+bc-bd+6\,{c}^{2}\\
  & \quad +8\,cd+{d}^{2}+4\,a-b-6\,c-d+1 ) {s}^{2}+ ( 2\,{a}^{3}{d}^{2}-3
  \,{a}^{2}bcd+{a}^{2}b{d}^{2}-{a}^{2}c{d}^{2}+a{b}^{2}{c}^{2}-3\,a{b}^{
    2}cd\\
  & \quad +3\,ab{c}^{2}d+3\,abc{d}^{2}+2\,{b}^{3}{c}^{2}-2\,{b}^{2}{c}^{3}-3
  \,{b}^{2}{c}^{2}d+2\,{a}^{3}d+{a}^{2}bc+4\,{a}^{2}bd-{a}^{2}cd-2\,{a}^
  {2}{d}^{2}\\
  & \quad +4\,a{b}^{2}c+a{b}^{2}d-6\,ab{c}^{2}-12\,abcd-2\,ab{d}^{2}-2\,a{c}^{2}d-ac{d}^{2}+2\,{b}^{3}c-{b}^{2}cd+6\,b{c}^{3}\\
  & \quad +7\,b{c}^{2}d+{a}^{2}b-2\,{a}^{2}c-6\,{a}^{2}d+a{b}^{2}+5\,abc+abd+6\,a{c}^{2}+12\,acd+2\,a{d}^{2}-4\,{b}^{2}c\\
  & \quad -8\,b{c}^{2}-4\,{c}^{3}-4\,{c}^{2}d+2\,{a}^{2}-3\,ab-8\,ac-2\,ad+4\,bc+6\,{c}^{2}+cd+2\,a-2\,c ) s\\
  & \quad +{a}^{4}{d}^{2}-2\,{a}^{3}bcd+{a}^{3}b{d}^{2}-{a}^{3}c{d}^{2}+{a}^{2}{b }^{2}{c}^{2}-2\,{a}^{2}{b}^{2}cd
  +{a}^{2}{b}^{2}{d}^{2}+2\,{a}^{2}b{c}^{2}d\\
  & \quad +{a}^{2}bc{d}^{2}+{a}^{2}{c}^{2}{d}^{2}+a{b}^{3}{c}^{2}-2\,a{b}^{3
  }cd-a{b}^{2}{c}^{3}-2\,a{b}^{2}{c}^{2}d-2\,ab{c}^{3}d+{b}^{4}{c}^{2}\\
  &\quad +{b}^{3}{c}^{3}+{b}^{2}{c}^{4}+{a}^{3}bd+{a}^{3}cd-{a}^{3}{d}^{2}-{a}^{2
  }{b}^{2}c-{a}^{2}{b}^{2}d-{a}^{2}b{c}^{2}-2\,{a}^{2}bcd-2\,{a}^{2}b{d}^{2}\\
  & \quad -3\,{a}^{2}{c}^{2}d-{a}^{2}c{d}^{2}+a{b}^{3}c+3\,a{b}^{2}{c}^{2}+5
  \,a{b}^{2}cd+3\,ab{c}^{3}+5\,ab{c}^{2}d+2\,a{c}^{3}d-3\,{b}^{3}{c}^{2}\\
  & \quad -4\,{b}^{2}{c}^{3}-2\,b{c}^{4}-{a}^{3}d+{a}^{2}{b}^{2}+3\,{a}^{2}bc+2
  \,{a}^{2}bd+{a}^{2}{c}^{2}+4\,{a}^{2}cd+{a}^{2}{d}^{2}-4\,a{b}^{2}c\\
  & \quad -8\,ab{c}^{2}-4\,abcd-2\,a{c}^{3}-3\,a{c}^{2}d+4\,{b}^{2}{c}^{2}+5\,b{c}^{3}+{c}^{4}-2\,{a}^{2}b-2\,{a}^{2}c-{a}^{2}d\\
  & \quad +5\,abc+4\,a{c}^{2}+acd-3\,b{c}^{2}-2\,{c}^{3}+2 +{a}^{2}-2\,ac+{c}^{2} \Big)x\\
  & -27 s^3 \Big( 2\,{s}^{6}+6\, ( ad+bc+a+b-2\,c-2\,d+1 ) {s}^{5}+3\,
 ( 2\,{a}^{2}{d}^{2}+abcd+2\,{b}^{2}{c}^{2}+6\,{a}^{2}d\\
 & \quad +5\,abc+5\,abd-6\,acd-5\,a{d}^{2}+6\,{b}^{2}c-10\,b{c}^{2}-11\,bcd+2\,{a}^{2}+5
 \,ab-10\,ac-11\,ad\\
 & \quad +2\,{b}^{2}-bc-5\,bd+10\,{c}^{2}+16\,cd+5\,{d}^{2}+6
 \,a+b-10\,c-5\,d+2 ) {s}^{4}+ ( 2\,{a}^{3}{d}^{3}\\
 & \quad -3\,{a}^{2}bc{d}^{2}-3\,a{b}^{2}{c}^{2}d+2\,{b}^{3}{c}^{3}+18\,{a}^{3}{d}^{2}+3
 \,{a}^{2}bcd+12\,{a}^{2}b{d}^{2}-9\,{a}^{2}c{d}^{2}-3\,{a}^{2}{d}^{3}\\
 & \quad +12\,a{b}^{2}{c}^{2}+3\,a{b}^{2}cd+3\,ab{c}^{2}d+3\,abc{d}^{2}+18\,{b}^
 {3}{c}^{2}-24\,{b}^{2}{c}^{3}-30\,{b}^{2}{c}^{2}d+18\,{a}^{3}d\\
 & \quad +12\,{a}^{2}bc+36\,{a}^{2}bd-39\,{a}^{2}cd-48\,{a}^{2}{d}^{2}+36\,a{b}^{2}c+12
 \,a{b}^{2}d-60\,ab{c}^{2}-108\,abcd\\
 &\quad -18\,ab{d}^{2}+12\,a{c}^{2}d+15\,ac{d}^{2}-3\,a{d}^{3}+18\,{b}^{3}c-27\,{b}^{2}{c}^{2}-39\,{b}^{2}cd+60\,
 b{c}^{3}+99\,b{c}^{2}d\\
 & \quad +33\,bc{d}^{2}+2\,{a}^{3}+12\,{a}^{2}b-24\,{a}^{
   2}c-27\,{a}^{2}d+12\,a{b}^{2}+9\,abc-33\,abd+60\,a{c}^{2}\\
 & \quad  +114\,acd+57\,a{d}^{2}+2\,{b}^{3}-6\,{b}^{2}c-3\,{b}^{2}d-48\,b{c}^{2}-12\,bcd-3\,
 b{d}^{2}-40\,{c}^{3}-72\,{c}^{2}d\\
 & \quad -30\,c{d}^{2}+2\,{d}^{3}+18\,{a}^{2}+
 3\,ab-72\,ac-48\,ad-3\,{b}^{2}+15\,bc+12\,bd+60\,{c}^{2}+45\,cd\\
 & \quad -3\,{d}^{2}+18\,a-3\,b-24\,c-3\,d+2 ) {s}^{3}+3\, ( 2\,{a}^{4}{d}^
 {3}-4\,{a}^{3}bc{d}^{2}+{a}^{3}b{d}^{3}-{a}^{3}c{d}^{3}\\
 & \quad +{a}^{2}{b}^{2}{c}^{2}d+{a}^{2}{b}^{2}c{d}^{2}-{a}^{2}b{c}^{2}{d}^{2}-6\,{a}^{2}bc{d}
 ^{3}+a{b}^{3}{c}^{3}-4\,a{b}^{3}{c}^{2}d+4\,a{b}^{2}{c}^{3}d+9\,a{b}^{
   2}{c}^{2}{d}^{2}\\
 & \quad +2\,{b}^{4}{c}^{3}-2\,{b}^{3}{c}^{4}-3\,{b}^{3}{c}^{3}
 d+6\,{a}^{4}{d}^{2}-{a}^{3}bcd+9\,{a}^{3}b{d}^{2}-5\,{a}^{3}c{d}^{2}-3
 \,{a}^{3}{d}^{3}+{a}^{2}{b}^{2}{c}^{2}\\
 & \quad -8\,{a}^{2}{b}^{2}cd+{a}^{2}{b}^{2}{d}^{2}+6\,{a}^{2}b{c}^{2}d+12\,{a}^{2}bc{d}^{2}+2\,{a}^{2}b{d}^{3}
 +2\,{a}^{2}{c}^{2}{d}^{2}+4\,{a}^{2}c{d}^{3}+9\,a{b}^{3}{c}^{2}\\
 & \quad -a{b}^{
   3}cd-12\,a{b}^{2}{c}^{3}-24\,a{b}^{2}{c}^{2}d-6\,a{b}^{2}c{d}^{2}-9\,a
 b{c}^{3}d-12\,ab{c}^{2}{d}^{2}+3\,abc{d}^{3}+6\,{b}^{4}{c}^{2}\\
 & \quad -7\,{b}^{3}{c}^{3}-11\,{b}^{3}{c}^{2}d+12\,{b}^{2}{c}^{4}+20\,{b}^{2}{c}^{3}d+
 6\,{b}^{2}{c}^{2}{d}^{2}+2\,{a}^{4}d+{a}^{3}bc+9\,{a}^{3}bd-8\,{a}^{3}cd\\
 & \quad -18\,{a}^{3}{d}^{2}+10\,{a}^{2}{b}^{2}c+10\,{a}^{2}{b}^{2}d-12\,{a}^{2}b{c}^{2}-37\,{a}^{2}bcd-20\,{a}^{2}b{d}^{2}
 +3\,{a}^{2}{c}^{2}d+5\,{a}^{2}c{d}^{2}\\
 & \quad -3\,{a}^{2}{d}^{3}+9\,a{b}^{3}c+a{b}^{3}d-2\,a{b}^{2}{c}
 ^{2}-7\,a{b}^{2}cd+2\,a{b}^{2}{d}^{2}+30\,ab{c}^{3}+70\,ab{c}^{2}d+26\,abc{d}^{2}\\
 & \quad -2\,ab{d}^{3}+4\,a{c}^{3}d+5\,a{c}^{2}{d}^{2}-ac{d}^{3}+2
 \,{b}^{4}c-6\,{b}^{3}{c}^{2}-2\,{b}^{3}cd-9\,{b}^{2}{c}^{3}+3\,{b}^{2}
 {c}^{2}d-{b}^{2}c{d}^{2}\\
 & \quad -20\,b{c}^{4}-33\,b{c}^{3}d-11\,b{c}^{2}{d}^{2}+{a}^{3}b-2\,{a}^{3}c-{a}^{3}d+{a}^{2}{b}^{2}+7\,{a}^{2}bc-5\,{a}^{2}
 bd+12\,{a}^{2}{c}^{2}\\
 & \quad +30\,{a}^{2}cd+27\,{a}^{2}{d}^{2}+a{b}^{3}-11\,a{b}^{2}c-9\,a{b}^{2}d-39\,ab{c}^{2}-24\,abcd+ab{d}^{2}-20\,a{c}^{3}\\
 & \quad -48\,a{c}^{2}d-25\,ac{d}^{2}+2\,a{d}^{3}-5\,{b}^{3}c+9\,{b}^{2}{c}^{2}+10
 \,{b}^{2}cd+34\,b{c}^{3}+23\,b{c}^{2}d+10\,{c}^{4}\\
 & \quad +16\,{c}^{3}d+5\,{c}^{2}{d}^{2}+2\,{a}^{3}-{a}^{2}b-18\,{a}^{2}c-18\,{a}^{2}d+a{b}^{2}+19
 \,abc+11\,abd+36\,a{c}^{2}\\
 & \quad +36\,acd-3\,a{d}^{2}-21\,b{c}^{2}-10\,bcd-20\,{c}^{3}-15\,{c}^{2}d+c{d}^{2}+6\,{a}^{2}-4\,ab-18\,ac\\
 & \quad -3\,ad+5\,bc+12\,{c}^{2}+2\,cd+2\,a-2\,c ) {s}^{2}+3\, ( 2\,{a}^{5}{d}^{3}
 -5\,{a}^{4}bc{d}^{2}+2\,{a}^{4}b{d}^{3}-2\,{a}^{4}c{d}^{3}\\
 & \quad +4\,{a}^{3}{b}^{2}{c}^{2}d-2\,{a}^{3}{b}^{2}c{d}^{2}-{a}^{3}{b}^{2}{d}^{3}+2\,{a}^{3}b{c}^{2}{d}^{2}-10\,{a}^{3}bc{d}^{3}
 -{a}^{3}{c}^{2}{d}^{3}-{a}^{2}{b}^{3}{c}^{3}\\
 & \quad -2\,{a}^{2}{b}^{3}{c}^{2}d+4\,{a}^{2}{b}^{3}c{d}^{2}+2\,{
   a}^{2}{b}^{2}{c}^{3}d+19\,{a}^{2}{b}^{2}{c}^{2}{d}^{2}-3\,{a}^{2}{b}^{
   2}c{d}^{3}+4\,{a}^{2}b{c}^{3}{d}^{2}+3\,{a}^{2}b{c}^{2}{d}^{3}\\
 & \quad +2\,a{b}^{4}{c}^{3}-5\,a{b}^{4}{c}^{2}d-2\,a{b}^{3}{c}^{4}-8\,a{b}^{3}{c}^{3}d
 +6\,a{b}^{3}{c}^{2}{d}^{2}-5\,a{b}^{2}{c}^{4}d-6\,a{b}^{2}{c}^{3}{d}^{2}+2\,{b}^{5}{c}^{3}\\
 & \quad -{b}^{4}{c}^{4}-3\,{b}^{4}{c}^{3}d+2\,{b}^{3}{c}^{
   5}+3\,{b}^{3}{c}^{4}d+2\,{a}^{5}{d}^{2}-{a}^{4}bcd+6\,{a}^{4}b{d}^{2}-{a}^{4}c{d}^{2}-3\,{a}^{4}{d}^{3}-{a}^{3}{b}^{2}{c}^{2}\\
 & \quad-11\,{a}^{3}{b}^{2}cd+5\,{a}^{3}{b}^{2}{d}^{2}+3\,{a}^{3}b{c}^{2}d+19\,{a}^{3}bc{d}^{2}+4\,{a}^{3}b{d}^{3}+3\,{a}^{3}{c}^{2}{d}^{2}
 +8\,{a}{3}c{d}^{3}+5\,{a}^{2}{b}^{3}{c}^{2}\\
 & \quad -11\,{a}^{2}{b}^{3}cd-{a}^{2}{b}^{3}{d}^{2}-2\,{a}
 ^{2}{b}^{2}{c}^{3}-21\,{a}^{2}{b}^{2}{c}^{2}d-12\,{a}^{2}{b}^{2}c{d}^{
   2}+2\,{a}^{2}{b}^{2}{d}^{3}-15\,{a}^{2}b{c}^{3}d\\
 & \quad -38\,{a}^{2}b{c}^{2}{d}^{2}+5\,{a}^{2}bc{d}^{3}-3\,{a}^{2}{c}^{3}{d}^{2}-{a}^{2}{c}^{2}{d}^{3}
 +6\,a{b}^{4}{c}^{2}-a{b}^{4}cd+5\,a{b}^{3}{c}^{3}+10\,a{b}^{3}{c}^{2}d\\
 & \quad -a{b}^{3}c{d}^{2}+12\,a{b}^{2}{c}^{4}+34\,a{b}^{2}{c}^{3}d-7\,a{b}^{
   2}{c}^{2}{d}^{2}+11\,ab{c}^{4}d+11\,ab{c}^{3}{d}^{2}+2\,{b}^{5}{c}^{2}
 -2\,{b}^{4}{c}^{3}\\
 & \quad -{b}^{4}{c}^{2}d-4\,{b}^{3}{c}^{4}+2\,{b}^{3}{c}^{3}d-8\,{b}^{2}{c}^{5}-10\,{b}^{2}{c}^{4}d+2\,{a}^{4}bd-{a}^{4}cd-8\,{a}^
 {4}{d}^{2}+4\,{a}^{3}{b}^{2}c\\
 & \quad +5\,{a}^{3}{b}^{2}d-2\,{a}^{3}b{c}^{2}-12\,{a}^{3}bcd-22\,{a}^{3}b{d}^{2}-2\,{a}^{3}{c}^{2}d-5\,{a}^{3}c{d}^{2}
 -3\,{a}^{3}{d}^{3}+5\,{a}^{2}{b}^{3}c\\
 & \quad +4\,{a}^{2}{b}^{3}d+4\,{a}^{2}{b}^{2}{c}^{2}+24\,{a}^{2}{b}^{2}cd-{a}^{2}{b}^{2}{d}^{2}+12\,{a}^{2}b{c}
 ^{3}+57\,{a}^{2}b{c}^{2}d+28\,{a}^{2}bc{d}^{2}-4\,{a}^{2}b{d}^{3}\\
 & \quad +9\,{a}^{2}{c}^{3}d+15\,{a}^{2}{c}^{2}{d}^{2}-2\,{a}^{2}c{d}^{3}+2\,a{b}^{4}c-19\,a{b}^{3}{c}^{2}
 -4\,a{b}^{3}cd-32\,a{b}^{2}{c}^{3}-26\,a{b}^{2}{c}^{2}d\\
 & \quad +4\,a{b}^{2}c{d}^{2}-20\,ab{c}^{4}-52\,ab{c}^{3}d-5\,ab{c}^{2}{d}^{2}-6\,a{c}^{4}d-5\,a{c}^{3}{d}^{2}
 -7\,{b}^{4}{c}^{2}+5\,{b}^{3}{c}^{3}\\
 & \quad +6\,{b}^{3}{c}^{2}d+21\,{b}^{2}{c}^{4}+10\,{b}^{2}{c}^{3}d+10\,b{c}^{5}
 +11\,b{c}^{4}d+{a}^{4}d-{a}^{3}{b}^{2}+3\,{a}^{3}bc+3\,{a}^{3}bd\\
 & \quad +2\,{a}^{3}{c}^{2}+10\,{a}^{3}cd+17\,{a}^{3}{d}^{2}-{a}^{2}{b}^{3}
 -12\,{a}^{2}{b}^{2}c-11\,{a}^{2}{b}^{2}d-26\,{a}^{2}b{c}^{2}-40\,{a}^{2}bcd\\
 & \quad +5\,{a}^{2}b{d}^{2}-8\,{a}^{2}{c}^{3}-33\,{a}^{2}{c}^{2}d-20\,{a}^{2}c{d}^{2}+2\,{a}^{2}{d}^{3}+36\,a{b}^{2}{c}^{2}
 +15\,a{b}^{2}cd+49\,ab{c}^{3}\\
 & \quad +45\,ab{c}^{2}d-5\,abc{d}^{2}+10\,a{c}^{4}+26\,a{c}^{3}d+6\,a{c}^{2}{d}^{2}+5\,{b}^{3}{c}^{2}-20\,{b}^{2}{c}^{3}-10\,{b}^{2}{c}^{2}d\\
 & \quad -26\,b{c}^{4}-14\,b{c}^{3}d-4\,{c}^{5}-4\,{c}^{4}d-{a}^{3}b-4\,{a}^{3}c-8\,{a}^{3}d+4\,{a}^{2}{b}^{2}+19\,{a}^{2}bc\\
 & \quad +10\,{a}^{2}bd+18\,{a}^{2}{c}^{2}+27\,{a}^{2}cd-3\,{a}^{2}{d}^{2}-10\,a{b}^{2}c-41\,ab{c}^{2}-14\,abcd-24\,a{c}^{3}\\
 & \quad -24\,a{c}^{2}d+2\,ac{d}^{2}+5\,{b}^{2}{c}^{2}+23\,b{c}^{3}+6\,b{c}^{2}d+10\,{c}^{4}+5\,{c}^{3}d+2\,{a}^{3}-5\,{a}^{2}b\\
 & \quad -12\,{a}^{2}c-3\,{a}^{2}d+12\,abc+18\,a{c}^{2}+4\,acd-7\,b{c}^{2}-8\,{c}^{3}-{c}^{2}d+2\,{a}^{2}-4\,ac+2\,{c}^{2} ) s\\
 & \quad + ( {a}^{2}d-abc+2\,abd+acd-2\,{b}^{2}c-b{c}^{2}-ab-ac-2\,ad+3\,bc+{c}^{2}+a-c ) ( 2\,{a}^{2}d\\
 & \quad -2\,abc+abd-acd-{b}^{2}c+b{c}^{2}+ab+ac-ad-{c}^{2}-a+c )  ( {a}^{2}d-abc-abd\\
 & \quad -2\,acd+{b}^{2}c+2\,b{c}^{2}+2\,ab +2\,ac+ad-3\,bc-2\,{c}^{2}-2\,a+2\,c )  \Big).
\end{align*}

\subsection{Class 2.7 : $\I_2^* + 8\I_2$ with MWG = $\left( \mathbb{Z}/2\mathbb{Z}\right)^2$}
\begin{equation*}
    \begin{aligned}
        y^2 = & \Big(x-t(bt+ct-t-1)(abt-bct-a+1) \Big)\\
        &  \times \Big( x-a t (bt-dt-1)(bt+ct-t-1) \Big)\\
        &  \times \Big( x-t (1-ct)(bt-dt-1)(abt-bct-a+1) \Big).
    \end{aligned}
\end{equation*}

\subsection{Class 2.8 : $\I_2^* + \I_0^* + 4\I_2$ with MWG = $\mathbb{Z}/2\mathbb{Z}$}
\begin{equation*}
    \begin{aligned}
        y^2 &= x^3 + t  \Big( (ad+bc+acd-bc^2-2cd)t^{2} + (ad-2bc-2a+b+c-2d+1)t -b+1 \Big)  x^2\\
        &  \quad + t^3 (cdt+d-1) \big((a-c)t-1 \big) \big((ad-bc+b-d)t-a-b+1 \big) x.
    \end{aligned}
\end{equation*}

\subsection{Class 2.9 : $2\I_2^* + 2\I_2+4\I_1$ with MWG = $\{ 0 \}$}
\begin{align*}
    y^2 &= x^3 -27\,t^2\, \Big( (c+d-1)^{2}(a-1)^{2}{t}^{4}-( {a}^
    {2}cd-2\,{a}^{2}{d}^{2}-2\,ab{c}^{2}+abcd+{a}^{2}c+3\,{a}^{2}d\\
    & \quad +3\,abc+abd+a{c}^{2}+acd+6\,a{d}^{2}-b{c}^{2}-4\,bcd-{a}^{2}
    -ab-4\,ac-9\,ad+3\,bc\\
    & \quad +2\,bd-{c}^{2}-2\,cd-4\,{d}^{2}+3\,a-2\,b+3\,c
    +6\,d-2 ) {t}^{3}+ ( {a}^{2}{d}^{2}-abcd+{b}^{2}{c}^{2}-{a}^{2}d\\
    & \quad -2\,abc-2\,abd+3\,acd-6\,a{d}^{2}-{b}^{2}c-b{c}^{2}
    +8\,bcd+{a}^{2}+3\,ab-ac+6\,ad+{b}^{2}\\
    & \quad -4\,bc-6\,bd+{c}^{2}-2\,cd+6\,{d}^{2}-a+4\,b
    -c-6\,d+1 ) {t}^{2}- ( abd+2\,a{d}^{2}+{b}^{2}c\\
    & \quad -4\,bcd-2\,ab-ad-2\,{b}^{2}+bc+6\,bd
    +2\,cd-4\,{d}^{2}-2\,b+2\,d ) t+ ( b-d ) ^{2} \Big) \, x\\
    & + 27\,t^3\,\Big(2\, ( c+d-1 ) ^{3} ( a-1 ) ^{3}{t}^{6}
    -3\, ( a-1 )  ( c+d-1 )  ( {a}^{2}cd-2\,{a}^{2}{d}^{2}-2\,ab{c}^{2}\\
    & \quad +abcd+{a}^{2}c+3\,{a}^{2}d+3\,abc+abd+a{c}^{2}+acd
    +6\,a{d}^{2}-b{c}^{2}-4\,bcd-{a}^{2}\\
    & \quad -ab-4\,ac-9\,ad+3\,bc+2\,bd-{c}^{2}-2\,cd-4\,{d}^{2}
    +3\,a-2\,b+3\,c+6\,d-2 ) {t}^{5}\\
    & \quad -3\, ( {a}^{3}c{d}^{2}-2\,{a}^{3}{d}^{3}+2\,{a}^{2}b{c}^{2}d
    +2\,{a}^{2}bc{d}^{2}-2\,a{b}^{2}{c}^{3}+a{b}^{2}{c}^{2}d
    -4\,{a}^{3}cd+4\,{a}^{3}{d}^{2}\\
    & \quad +3\,{a}^{2}b{c}^{2}-5\,{a}^{2}bcd+3\,{a}^{2}b{d}^{2}
    +{a}^{2}{c}^{2}d-2\,{a}^{2}c{d}^{2}+12\,{a}^{2}{d}^{3}
    +4\,a{b}^{2}{c}^{2}-4\,a{b}^{2}cd\\
    & \quad +2\,ab{c}^{3}+2\,ab{c}^{2}d-15\,abc{d}^{2}-{b}^{2}{c}^{3}
    +5\,{b}^{2}{c}^{2}d+{a}^{3}c-{a}^{3}d-2\,{a}^{2}bc-2\,{a}^{2}bd\\
    & \quad -4\,{a}^{2}{c}^{2}+{a}^{2}cd-24\,{a}^{2}{d}^{2}-a{b}^{2}c
    +a{b}^{2}d-11\,ab{c}^{2}+13\,abcd+4\,ab{d}^{2}+a{c}^{3}\\
    & \quad -3\,a{c}^{2}d-3\,ac{d}^{2}-20\,a{d}^{3}-{b}^{2}{c}^{2}
    -5\,{b}^{2}cd+4\,b{c}^{3}+5\,b{c}^{2}d+16\,bc{d}^{2}-{a}^{3}\\
    & \quad -{a}^{2}b+5\,{a}^{2}c+13\,{a}^{2}d-a{b}^{2}+8\,abc-5\,abd
    +5\,a{c}^{2}+12\,acd+40\,a{d}^{2}+4\,{b}^{2}c\\
    & \quad +2\,{b}^{2}d-10\,b{c}^{2}-26\,bcd-10\,b{d}^{2}-{c}^{3}
    +2\,{c}^{2}d+4\,c{d}^{2}+10\,{d}^{3}-{a}^{2}+ab-10\,ac\\
    & \quad -24\,ad-2\,{b}^{2}+12\,bc+16\,bd-{c}^{2}-9\,cd-20\,{d}^{2}
    +4\,a-6\,b+4\,c+12\,d-2 ) {t}^{4}\\
    & \quad + ( 2\,{a}^{3}{d}^{3}-3\,{a}^{2}bc{d}^{2}-3\,a{b}^{2}{c}^{2}d
    +2\,{b}^{3}{c}^{3}-3\,{a}^{3}{d}^{2}+24\,{a}^{2}bcd-9\,{a}^{2}b{d}^{2}
    +12\,{a}^{2}c{d}^{2}\\
    & \quad -24\,{a}^{2}{d}^{3}-9\,a{b}^{2}{c}^{2}+24\,a{b}^{2}cd
    -18\,ab{c}^{2}d+45\,abc{d}^{2}-3\,{b}^{3}{c}^{2}-3\,{b}^{2}{c}^{3}-
    30\,{b}^{2}{c}^{2}d\\
    & \quad -3\,{a}^{3}d-9\,{a}^{2}bc-6\,{a}^{2}bd-18\,{a}^{2}cd
    +36\,{a}^{2}{d}^{2}-6\,a{b}^{2}c-9\,a{b}^{2}d+24\,ab{c}^{2}
    -24\,abcd\\
    & \quad -18\,ab{d}^{2}+12\,a{c}^{2}d-27\,ac{d}^{2}
    +60\,a{d}^{3}-3\,{b}^{3}c+15\,{b}^{2}{c}^{2}+45\,{b}^{2}cd
    -3\,b{c}^{3}+15\,b{c}^{2}d\\
    & \quad -72\,bc{d}^{2}+2\,{a}^{3}+12\,{a}^{2}b-3\,{a}^{2}c-6\,{a}^{2}d
    +12\,a{b}^{2}-33\,abc+9\,abd-3\,a{c}^{2}+9\,acd\\
    & \quad -90\,a{d}^{2}+2\,{b}^{3}-27\,{b}^{2}c-24\,{b}^{2}d
    +15\,b{c}^{2}+72\,bcd+60\,b{d}^{2}+2\,{c}^{3}-9\,{c}^{2}d+12\,c{d}^{2}\\
    & \quad -40\,{d}^{3}-3\,{a}^{2}+3\,ab+12\,ac+36\,ad+18\,{b}^{2}
    -27\,bc-72\,bd-3\,{c}^{2}+3\,cd+60\,{d}^{2}\\
    & \quad -3\,a+18\,b-3\,c-24\,d+2 ) {t}^{3}-3\, ( {a}^{2}b{d}^{2}
    +2\,{a}^{2}{d}^{3}-4\,a{b}^{2}cd-5\,abc{d}^{2}+{b}^{3}{c}^{2}\\
    & \quad +5\,{b}^{2}{c}^{2}d+2\,{a}^{2}bd-2\,{a}^{2}{d}^{2}
    +3\,a{b}^{2}c+3\,a{b}^{2}d+abcd+4\,ab{d}^{2}+5\,ac{d}^{2}-10\,a{d}^{3}\\
    & \quad +2\,{b}^{3}c-4\,{b}^{2}{c}^{2}-15\,{b}^{2}cd-5\,b{c}^{2}d
    +16\,bc{d}^{2}-2\,{a}^{2}b-{a}^{2}d-5\,a{b}^{2}+2\,abc+abd\\
    & \quad -2\,acd+10\,a{d}^{2}-2\,{b}^{3}+6\,{b}^{2}c+12\,{b}^{2}d
    +b{c}^{2}-6\,bcd-20\,b{d}^{2}+2\,{c}^{2}d-6\,c{d}^{2}+10\,{d}^{3}\\
    & \quad -ab-2\,ad-6\,{b}^{2}+2\,bc+16\,bd+cd-10\,{d}^{2}
    -2\,b+2\,d ) {t}^{2}-3\, ( b-d )  ( abd+2\,a{d}^{2}\\
    & \quad +{b}^{2}c-4\,bcd-2\,ab-ad-2\,{b}^{2}+bc+6\,bd
    +2\,cd-4\,{d}^{2}-2\,b+2\,d ) t+2\, ( b-d ) ^{3} \Big).
\end{align*}

\subsection{Class 2.10 : $\I_2^* + 2\I_0^*+8\I_1$ with MWG = $\{ 0 \}$}
\begin{align*}
    y^2 &= x^3 -27 (ad-bc)^2(s-ab+ad)(s-ab+bc)
    \Big(({a}^{2}{d}^{2}-2\,abcd+{b}^{2}{c}^{2}+abd+acd-{b}^{2}c\\
    & \qquad -b{c}^{2}-ad+{b}^{2}+3\,bc+{c}^{2}-2\,b-2\,c+1 ) {s}^{2}- ( 2\,{a}^{3}b{d}^{2}-2\,{a}^{3}{d}^{3}-4\,{a}^{2}{b}^{2}cd\\
    & \qquad +2\,{a}^{2}bc{d}^{2}+2\,a{b}^{3}{c}^{2}+2\,a{b}^{2}{c}^{2}d-2\,{b}^{3}{c}^{3}+2\,{a}^{3}{d}^{2}
    +2\,{a}^{2}{b}^{2}d+{a}^{2}bcd-{a}^{2}b{d}^{2}-{a}^{2}c{d}^{2}\\
    & \qquad +2\,{a}^{2}{d}^{3}-2\,a{b}^{3}c-3\,a{b}^{2}{c}^{2}-a{b}^{2}cd-ab{c}^{2}d-abc{d}^{2}+2\,{b}^{3}{c}^{2}
    +2\,{b}^{2}{c}^{3}-{b}^{2}{c}^{2}d-{a}^{2}bd\\
    & \qquad -2\,{a}^{2}cd-4\,{a}^{2}{d}^{2}+2\,a{b}^{3}+5\,a{b}^{2}c+4\,ab{c}^{2}+2\,abcd-2\,ab{d}^{2}+ac{d}^{2}-2\,{b}^{3}c-2\,{b}^{2}{c}^{2}\\
    & \qquad +2\,{b}^{2}cd-2\,b{c}^{3}-b{c}^{2}d+2\,{a}^{2}d-4\,a{b}^{2}-6\,abc+2\,abd+2\,acd+2\,a{d}^{2}+2\,{b}^{2}c+2\,b{c}^{2}\\
    & \qquad -2\,bcd+2\,ab-2\,ad ) s + {a}^{4}{b}^{2}{d}^{2}-2\,{a}^{4}b{d}^{3}+{a}^{4}{d}^{4}-2\,{a}^{3}{b}^{3}cd
    +2\,{a}^{3}{b}^{2}c{d}^{2}+{a}^{2}{b}^{4}{c}^{2}\\
    & \qquad +2\,{a}^{2}{b}^{3}{c}^{2}d-2\,{a}^{2}{b}^{2}{c}^{2}{d}^{2}-2\,a{b}^{4}{c}^{3}+{b}^{4}{c}^{4}+2\,{a}^{4}b{d}^{2}
    -2\,{a}^{4}{d}^{3}+{a}^{3}{b}^{3}d-{a}^{3}{b}^{2}{d}^{2}\\
    & \qquad +{a}^{3}bc{d}^{2}+2\,{a}^{3}b{d}^{3}-2\,{a}^{3}{d}^{4}-{a}^{2}{b}^{4}c-2\,{a}^{2}{b}^{3}{c}^{2}-{a}^{2}{b}^{3}cd
    -2\,{a}^{2}{b}^{2}{c}^{2}d+{a}^{2}b{c}^{2}{d}^{2}\\
    & \qquad +2\,{a}^{2}bc{d}^{3}+2\,a{b}^{4}{c}^{2}+3\,a{b}^{3}{c}^{3}-a{b}^{3}{c}^{2}d-a{b}^{2}{c}^{2}{d}^{2}-{b}^{4}{c}^{3}
    -{b}^{3}{c}^{4}+{b}^{3}{c}^{3}d+{a}^{4}{d}^{2}\\
    & \qquad -4\,{a}^{3}bcd-4\,{a}^{3}b{d}^{2}+4\,{a}^{3}{d}^{3}+{a}^{2}{b}^{4}+2\,{a}^{2}{b}^{3}c+4\,{a}^{2}{b}^{2}{c}^{2}
    +{a}^{2}{b}^{2}cd-2\,{a}^{2}{b}^{2}{d}^{2}\\
    & \qquad +2\,{a}^{2}b{c}^{2}d-3\,{a}^{2}bc{d}^{2}+{a}^{2}{d}^{4}-2\,a{b}^{4}c-a{b}^{3}{c}^{2}+2\,a{b}^{3}cd-4\,a{b}^{2}{c}^{3}
    +2\,a{b}^{2}{c}^{2}d\\
    & \qquad +2\,a{b}^{2}c{d}^{2}-ab{c}^{2}{d}^{2}-2\,abc{d}^{3}+{b}^{4}{c}^{2}-{b}^{3}{c}^{3}-2\,{b}^{3}{c}^{2}d+{b}^{2}{c}^{4}
    +{b}^{2}{c}^{3}d+{b}^{2}{c}^{2}{d}^{2}\\
    & \qquad +2\,{a}^{3}bd-2\,{a}^{3}{d}^{2}-2\,{a}^{2}{b}^{3}-4\,{a}^{2}{b}^{2}c+2\,{a}^{2}{b}^{2}d+4\,{a}^{2}bcd
    +2\,{a}^{2}b{d}^{2}-2\,{a}^{2}{d}^{3}\\
    & \qquad +2\,a{b}^{3}c+2\,a{b}^{2}{c}^{2}-4\,a{b}^{2}cd-2\,ab{c}^{2}d+2\,abc{d}^{2}+{a}^{2}{b}^{2}-2\,{a}^{2}bd+{a}^{2}{d}^{2}\Big) \, x\\
    & -27(ad-bc)^3(s-ab+ad)^3(s-ab+bc)^3  \Big((2\,ad-2\,bc+b+c-1)  ( ad-bc+\\
    & \qquad 2\,b+2\,c-2 )  ( ad-bc-b-c+1) {s}^{3}-3\, ( 2\,{a}^{4}b{d}^{3}-2\,{a}^{4}{d}^{4}-6\,{a }^{3}{b}^{2}c{d}^{2}+4\,{a}^{3}bc{d}^{3}\\
    & \qquad +6\,{a}^{2}{b}^{3}{c}^{2}d-2\,a{b}^{4}{c}^{3}-4\,a{b}^{3}{c}^{3}d+2\,{b}^{4}{c}^{4}+2\,{a}^{4}{d}^{3}+3\,{a}^{3}{b}^{2}{d}^{2}
    -2\,{a}^{3}b{d}^{3}-2\,{a}^{3}c{d}^{3}\\
    & \qquad +2\,{a}^{3}{d}^{4}-6\,{a}^{2}{b}^{3}cd-6\,{a}^{2}{b}^{2}{c}^{2}d+{a}^{2}{b}^{2}c{d}^{2}+{a}^{2}b{c}^{2}{d}^{2}
    -3\,{a}^{2}bc{d}^{3}+3\,a{b}^{4}{c}^{2}+4\,a{b}^{3}{c}^{3}\\
    & \qquad +4\,a{b}^{3}{c}^{2}d+4\,a{b}^{2}{c}^{3}d-3\,{b}^{4}{c}^{3}-3\,{b}^{3}{c}^{4}+{b}^{3}{c}^{3}d-{a}^{3}b{d}^{2}
    -{a}^{3}c{d}^{2}-3\,{a}^{3}{d}^{3}-3\,{a}^{2}{b}^{3}d\\
    & \qquad +2\,{a}^{2}{b}^{2}cd+{a}^{2}{b}^{2}{d}^{2}-4\,{a}^{2}b{c}^{2}d-2\,{a}^{2}bc{d}^{2}-{a}^{2}b{d}^{3}+{a}^{2}{c}^{2}{d}^{2}
    +2\,{a}^{2}c{d}^{3}+3\,a{b}^{4}c-a{b}^{3}{c}^{2}\\
    & \qquad +2\,a{b}^{3}cd+5\,a{b}^{2}{c}^{3}+3\,a{b}^{2}{c}^{2}d-a{b}^{2}c{d}^{2}
    +2\,ab{c}^{3}d+2\,ab{c}^{2}{d}^{2}-3\,{b}^{4}{c}^{2}+2\,{b}^{3}{c}^{3}+2\,{b}^{3}{c}^{2}d\\
    & \qquad -3\,{b}^{2}{c}^{4}-4\,{b}^{2}{c}^{3}d+{a}^{3}{d}^{
      2}+5\,{a}^{2}{b}^{2}d+8\,{a}^{2}bcd+2\,{a}^{2}b{d}^{2}+2\,{a}^{2}{c}^{
      2}d+2\,{a}^{2}c{d}^{2}+{a}^{2}{d}^{3}\\
    & \qquad -2\,a{b}^{4}-11\,a{b}^{3}c-15\,a{
      b}^{2}{c}^{2}-3\,a{b}^{2}cd+2\,a{b}^{2}{d}^{2}-4\,ab{c}^{3}-3\,ab{c}^{
      2}d+2\,abc{d}^{2}-a{c}^{2}{d}^{2}\\
    & \qquad +2\,{b}^{4}c+7\,{b}^{3}{c}^{2}-2\,{b}^{3}cd+7\,{b}^{2}{c}^{3}-3\,{b}^{2}{c}^{2}d+2\,b{c}^{4}+b{c}^{3}d-4\,{a}^{2}bd
    -4\,{a}^{2}cd-3\,{a}^{2}{d}^{2}\\
    & \qquad +6\,a{b}^{3}+16\,a{b}^{2}c-2\,a{b}^{2}d+10\,ab{c}^{2}-3\,abcd-4\,ab{d}^{2}-2\,a{c}^{2}d-ac{d}^{2}-4\,{b}^{3}c\\
    & \qquad -6\,{b}^{2}{c}^{2}+4\,{b}^{2}cd-4\,b{c}^{3}+b{c}^{2}d+2\,{a}^{2}d-6\,a{b}^{2}-8\,abc+4\,abd+4\,acd+2\,a{d}^{2}+2\,{b}^{2}c\\
    & \qquad +2\,b{c}^{2}-2\,bcd+2\,ab-2\,ad ) {s}^{2}+3\, ( 2\,{a}^{5}{b}^{2}{d}^
    {3}-4\,{a}^{5}b{d}^{4}+2\,{a}^{5}{d}^{5}-6\,{a}^{4}{b}^{3}c{d}^{2}+8\,{a}^{4}{b}^{2}c{d}^{3}\\
    & \qquad -2\,{a}^{4}bc{d}^{4}+6\,{a}^{3}{b}^{4}{c}^{2}d-4\,{a}^{3}{b}^{2}{c}^{2}{d}^{3}-2\,{a}^{2}{b}^{5}{c}^{3}
    -8\,{a}^{2}{b}^{4}{c}^{3}d+4\,{a}^{2}{b}^{3}{c}^{3}{d}^{2}+4\,a{b}^{5}{c}^{4}\\
    & \qquad +2\,a{b}^{4}{c}^{4}d-2\,{b}^{5}{c}^{5}+4\,{a}^{5}b{d}^{3}-4\,{a}^{5}{d}^{4}+3\,{a}^{4}{b}^{3}{d}^{2}
    -3\,{a}^{4}{b}^{2}c{d}^{2}-4\,{a}^{4}{b}^{2}{d}^{3}+{a}^{4}bc{d}^{3}\\
    & \qquad +5\,{a}^{4}b{d}^{4}+{a}^{4}c{d}^{4}-4\,{a}^{4}{d}
    ^{5}-6\,{a}^{3}{b}^{4}cd-6\,{a}^{3}{b}^{3}{c}^{2}d+2\,{a}^{3}{b}^{3}c{
      d}^{2}+2\,{a}^{3}{b}^{2}{c}^{2}{d}^{2}\\
    & \qquad -4\,{a}^{3}{b}^{2}c{d}^{3}+2\,{a}^{3}b{c}^{2}{d}^{3}+5\,{a}^{3}bc{d}^{4}+3\,{a}^{2}{b}^{5}{c}^{2}+5\,{
      a}^{2}{b}^{4}{c}^{3}+8\,{a}^{2}{b}^{4}{c}^{2}d+9\,{a}^{2}{b}^{3}{c}^{3}d\\
    & \qquad -4\,{a}^{2}{b}^{3}{c}^{2}{d}^{2}-4\,{a}^{2}{b}^{2}{c}^{3}{d}^{2}-6\,
    a{b}^{5}{c}^{3}-8\,a{b}^{4}{c}^{4}-2\,a{b}^{3}{c}^{4}d+a{b}^{3}{c}^{3}
    {d}^{2}+3\,{b}^{5}{c}^{4}+3\,{b}^{4}{c}^{5}\\
    & \qquad -2\,{b}^{4}{c}^{4}d+2\,{a}^{5}{d}^{3}+{a}^{4}{b}^{2}{d}^{2}-5\,{a}^{4}bc{d}^{2}
    -8\,{a}^{4}b{d}^{3}+{a}^{4}c{d}^{3}+10\,{a}^{4}{d}^{4}-3\,{a}^{3}{b}^{4}d\\
    & \qquad +4\,{a}^{3}{b}^{3}cd+2\,{a}^{3}{b}^{3}{d}^{2}-5\,{a}^{3}{b}^{2}{c}^{2}d-4\,{a}^{3}{b}^{2}c{d}^{2}
    -2\,{a}^{3}{b}^{2}{d}^{3}+2\,{a}^{3}b{c}^{2}{d}^{2}-6\,{a}^{3}bc{d}^{3}\\
    & \qquad +{a}^{3}b{d}^{4}-2\,{a}^{3}c{d}^{4}+2\,{a}^{3}{d}^{5}+3\,{a}^{2}{b}^{5}c
    -5\,{a}^{2}{b}^{4}{c}^{2}+4\,{a}^{2}{b}^{4}cd+8\,{a}^{2}{b}^{3}{c}^{3}-4\,{a}^{2}{b}^{3}c{d}^{2}\\
    & \qquad +7\,{a}^{2}{b}^{2}{c}^{3}d+13\,{a}^{2}{b}^{2}{c}^{2}{d}^{2}-2\,{a}^{2}b{c}^{3}{d}^{2}-3\,{a}^{2}bc{
      d}^{4}-6\,a{b}^{5}{c}^{2}+12\,a{b}^{4}{c}^{3}+3\,a{b}^{4}{c}^{2}d\\
    & \qquad -10\,a{b}^{3}{c}^{4}-10\,a{b}^{3}{c}^{3}d+3\,a{b}^{3}{c}^{2}{d}^{2}-a{b}^{2}{c}^{4}d
    -6\,a{b}^{2}{c}^{3}{d}^{2}+3\,{b}^{5}{c}^{3}-7\,{b}^{4}{c}^{4}\\
    & \qquad -4\,{b}^{4}{c}^{3}d+3\,{b}^{3}{c}^{5}+8\,{b}^{3}{c}^{4}d+{b}^{3}{c}^{3}{d}^{2}+3\,{a}^{4}b{d}^{2}-2\,{a}^{4}c{d}^{2}
    -8\,{a}^{4}{d}^{3}+4\,{a}^{3}{b}^{3}d\\
    & \qquad +8\,{a}^{3}{b}^{2}cd+4\,{a}^{3}{b}^{2}{d}^{2}+8\,{a}^{3}b{c}^{2}d+13\,{a}^{3}bc{d}^{2}-2\,{a}^{3}c{d}^{3}-8\,{a}^{3}{d}^{4}
    -2\,{a}^{2}{b}^{5}-10\,{a}^{2}{b}^{4}c\\
    & \qquad -17\,{a}^{2}{b}^{3}{c}^{2}-4\,{a}^{2}{b}^{3}cd+4\,{a}^{2}{b}^{3}{d}^{2}-8\,{a}^{2}{b}^{2}{c}^{3}-8\,{a}^{2}{b}^{2}{c}^{2}d
    +4\,{a}^{2}{b}^{2}c{d}^{2}-4\,{a}^{2}b{c}^{3}d\\
    & \qquad -2\,{a}^{2}b{c}^{2}{d}^{2}+9\,{a}^{2}bc{d}^{3}-2\,{a}^{2}b{d}^{4}
    +{a}^{2}c{d}^{4}+4\,a{b}^{5}c+12\,a{b}^{4}{c}^{2}-4\,a{b}^{4}cd+15\,a{b}^{3}{c}^{3}\\
    & \qquad -8\,a{b}^{3}{c}^{2}d-4\,a{b}^{3}c{d}^{2}+8\,a{b}^{2}{c}^{4}-2\,a{b}^{2}{c}^{2}{d}^{2}+4\,a{b}^{2}c{d}^{3}+2\,ab{c}^{3}{d}^{2}
    -2\,ab{c}^{2}{d}^{3}\\
    & \qquad -2\,{b}^{5}{c}^{2}-2\,{b}^{4}{c}^{3}+4\,{b}^{4}{c}^{2}d-2\,{b}^{3}{c}^{4}+{b}^{3}{c}^{3}d-2\,{b}^{3}{c}^{2}{d}^{2}
    -2\,{b}^{2}{c}^{5}-2\,{b}^{2}{c}^{4}d\\
    & \qquad +{b}^{2}{c}^{3}{d}^{2}+2\,{a}^{4}{d}^{2}-5\,{a}^{3}{b}^{2}d-12\,{a}^{3}bcd-5\,{a}^{3}b{d}^{2}+4\,{a}^{3}c{d}^{2}
    +10\,{a}^{3}{d}^{3}+6\,{a}^{2}{b}^{4}\\
    & \qquad +17\,{a}^{2}{b}^{3}c-4\,{a}^{2}{b}^{3}d+16\,{a}^{2}{b}^{2}{c}^{2}-6\,{a}^{2}{b}^{2}cd-8\,{a}^{2}{b}^{2}{d}^{2}
    -4\,{a}^{2}b{c}^{2}d-12\,{a}^{2}bc{d}^{2}\\
    & \qquad +4\,{a}^{2}b{d}^{3}+{a}^{2}c{d}^{3}+2\,{a}^{2}{d}^{4}-8\,a{b}^{4}c-13\,a{b}^{3}{c}^{2}+12\,a{b}^{3}cd
    -12\,a{b}^{2}{c}^{3}+13\,a{b}^{2}{c}^{2}d\\
    & \qquad +4\,ab{c}^{3}d-4\,abc{d}^{3}+2\,{b}^{4}{c}^{2}+{b}^{3}{c}^{3}-4\,{b}^{3}{c}^{2}d+2\,{b}^{2}{c}^{4}
    -{b}^{2}{c}^{3}d+2\,{b}^{2}{c}^{2}{d}^{2}+4\,{a}^{3}bd\\
    & \qquad -4\,{a}^{3}{d}^{2}-6\,{a}^{2}{b}^{3}-10\,{a}^{2}{b}^{2}c+8\,{a}^{2}{b}^{2}d+12\,{a}^{2}bcd
    +2\,{a}^{2}b{d}^{2}-2\,{a}^{2}c{d}^{2}-4\,{a}^{2}{d}^{3}\\
    & \qquad +4\,a{b}^{3}c+4\,a{b}^{2}{c}^{2}-8\,a{b}^{2}cd-4\,ab{c}^{2}d+4\,abc{d}^{2}+2\,{a}^{2}{b}^{2}-4\,{a}^{2}bd+2\,{a}^{2}{d}^{2}) s\\
    & \qquad -2\,{a}^{6}{b}^{3}{d}^{3}+6\,{a}^{6}{b}^{2}{d}^{4}-6\,{a}^{6}b{d}^{5}+2\,{a}^{6}{d}^{6}+6\,{a}^{5}{b}^{4}c{d}^{2}
    -12\,{a}^{5}{b}^{3}c{d}^{3}+6\,{a}^{5}{b}^{2}c{d}^{4}\\
    & \qquad -6\,{a}^{4}{b}^{5}{c}^{2}d+12\,{a}^{4}{b}^{3}{c}^{2}{d}^{3}-6\,{a}^{4}{b}^{2}{c}^{2}{d}^{4}+2\,{a}^{3}{b}^{6}{c}^{3}
    +12\,{a}^{3}{b}^{5}{c}^{3}d-12\,{a}^{3}{b}^{4}{c}^{3}{d}^{2}\\
    & \qquad -6\,{a}^{2}{b}^{6}{c}^{4}-6\,{a}^{2}{b}^{5}{c}^{4}d+6\,{a}^{2}{b}^{4}{c}^{4}{d}^{2}+6\,a{b}^{6}{c}^{5}-2\,{b}^{6}{c}^{6}
    -6\,{a}^{6}{b}^{2}{d}^{3}+12\,{a}^{6}b{d}^{4}\\
    & \qquad -6\,{a}^{6}{d}^{5}-3\,{a}^{5}{b}^{4}{d}^{2}+6\,{a}^{5}{b}^{3}c{d}^{2}+6\,{a}^{5}{b}^{3}{d}^{3}-9\,{a}^{5}{b}^{2}c{d}^{3}
    -9 \,{a}^{5}{b}^{2}{d}^{4}+3\,{a}^{5}bc{d}^{4}\\
    & \qquad +12\,{a}^{5}b{d}^{5}-6\,{a} ^{5}{d}^{6}+6\,{a}^{4}{b}^{5}cd+6\,{a}^{4}{b}^{4}{c}^{2}d-3\,{a}^{4}{b}^{4}c{d}^{2}
    -3\,{a}^{4}{b}^{3}{c}^{2}{d}^{2}+3\,{a}^{4}{b}^{3}c{d}^{3}\\
    & \qquad -3\,{a}^{4}{b}^{2}{c}^{2}{d}^{3}-12\,{a}^{4}{b}^{2}c{d}^{4}
    +3\,{a}^{4 }b{c}^{2}{d}^{4}+6\,{a}^{4}bc{d}^{5}-3\,{a}^{3}{b}^{6}{c}^{2}-6\,{a}^{3}{b}^{5}{c}^{3}\\
    & \qquad -12\,{a}^{3}{b}^{5}{c}^{2}d-15\,{a}^{3}{b}^{4}{c}^{3}d+12\,{a}^{3}{b}^{4}{c}^{2}{d}^{2}
    +9\,{a}^{3}{b}^{3}{c}^{3}{d}^{2}+3\,{a}^{3}{b}^{2}{c}^{2}{d}^{4}+9\,{a}^{2}{b}^{6}{c}^{3}\\
    & \qquad +15\,{a}^{2}{b}^{5}{c}^{4}+3\,{a}^{2}{b}^{5}{c}^{3}d+9\,{a}^{2}{b}^{4}{c}^{4}d-9\,{a}^{2}{b}^{4}{c}^{3}{d}^{2}
    -6\,{a}^{2}{b}^{3}{c}^{4}{d}^{2}-3\,{a}^{2}{b}^{3}{c}^{3}{d}^{3}\\
    & \qquad -9\,a{b}^{6}{c}^{4}-12\,a{b}^{5}{c}^{5}+6\,a{b}^{5}{c}^{4}d
    +3\,a{b}^{4}{c}^{4}{d}^{2}+3\,{b}^{6}{c}^{5}+3\,{b}^{5}{c}^{6}-3\,{b}^{5}{c}^{5}d\\
    & \qquad -6\,{a}^{6}b{d}^{3}+6\,{a}^{6}{d}^{4}-3\,{a}^{5}{b}^{3}{d}^{2}+12\,{a}^{5}{b}^{2}c{d}^{2}
    +15\,{a}^{5}{b}^{2}{d}^{3}-15\,{a}^{5}bc{d}^{3}-30\,{a}^{5}b{d}^{4}\\
    & \qquad +18\,{a}^{5}{d}^{5}+3\,{a}^{4}{b}^{5}d-6\,{a}^{4}{b}^{4}cd-3\,{a}^{4}{b}^{4}{d}^{2}+6\,{a}^{4}{b}^{3}{c}^{2}d
    +6\,{a}^{4}{b}^{3}c{d}^{2}+3\,{a}^{4}{b}^{3}{d}^{3}\\
    & \qquad +6\,{a}^{4}{b}^{2}{c}^{2}{d}^{2}+18\,{a}^{4}{b}^{2}c{d}^{3}-3\,{a}^{4}{b}^{2}{d}^{4}+3\,{a}^{4}b{c}^{2}{d}^{3}
    -18\,{a}^{4}bc{d}^{4}-6\,{a}^{4}b{d}^{5}+6\,{a }^{4}{d}^{6}\\
    & \qquad -3\,{a}^{3}{b}^{6}c+9\,{a}^{3}{b}^{5}{c}^{2}-6\,{a}^{3}{b}^{5}cd
    -12\,{a}^{3}{b}^{4}{c}^{3}+9\,{a}^{3}{b}^{4}{c}^{2}d+9\,{a}^{3}{b}^{4}c{d}^{2}-21\,{a}^{3}{b}^{3}{c}^{3}d\\
    & \qquad -27\,{a}^{3}{b}^{3}{c}^{2}{d}^{2}+3\,{a}^{3}{b}^{2}{c}^{3}{d}^{2}+15\,{a}^{3}{b}^{2}{c}^{2}{d}^{3}
    +12\,{a}^{3}{b}^{2}c{d}^{4}-6\,{a}^{3}b{c}^{2}{d}^{4}-12\,{a}^{3}bc{d}^{5}\\
    & \qquad +9\,{a}^{2}{b}^{6}{c}^{2}-30\,{a}^{2}{b}^{5}{c}^{3}-3\,{a}^{2}{b}^{5}{c}^{2}d
    +24\,{a}^{2}{b}^{4}{c}^{4}+6\,{a}^{2}{b}^{4}{c}^{3}d-12\,{a}^{2}{b}^{4}{c}^{2}{d}^{2}\\
    & \qquad +9\,{a}^{2}{b}^{3}{c}^{4}d+18\,{a}^{2}{b}^{3}
    {c}^{3}{d}^{2}-3\,{a}^{2}{b}^{3}{c}^{2}{d}^{3}-3\,{a}^{2}{b}^{2}{c}^{4}{d}^{2}
    +6\,{a}^{2}{b}^{2}{c}^{3}{d}^{3}+9\,{a}^{2}{b}^{2}{c}^{2}{d}^{4}\\
    & \qquad -9\,a{b}^{6}{c}^{3}+33\,a{b}^{5}{c}^{4}+12\,a{b}^{5}{c}^{3}d-15\,a{b}^{4}{c}^{5}
    -21\,a{b}^{4}{c}^{4}d+3\,a{b}^{4}{c}^{3}{d}^{2}-12\,a{b}^{3}{c}^{4}{d}^{2}\\
    & \qquad -6\,a{b}^{3}{c}^{3}{d}^{3}+3\,{b}^{6}{c}^{4}-12\,{b}^{5}{c}^{5}
    -6\,{b}^{5}{c}^{4}d+3\,{b}^{4}{c}^{6}+12\,{b}^{4}{c}^{5}d+3\,{b}^{4}{c}^{4}{d}^{2}-2\,{a}^{6}{d}^{3}\\
    & \qquad -6\,{a}^{5}{b}^{2}{d}^{2}+12\,{a}^{5}bc{d}^{2}+24\,{a}^{5}b{d}^{3}-18\,{a}^{5}{d}^{4}-3\,{a}^{4}{b}^{4}d
    -6\,{a}^{4}{b}^{3}cd-6\,{a}^{4}{b}^{3}{d}^{2}\\
    & \qquad -24\,{a}^{4}{b}^{2}{c}^{2}d-30\,{a}^{4}{b}^{2}c{d}^{2}+3\,{a}^{4}{b}^{2}{d}^{3}-6\,{a}^{4}b{c}^{2}{d}^{2}
    +36\,{a}^{4}bc{d}^{3}+24\,{a}^{4}b{d}^{4}\\
    & \qquad -18\,{a}^{4}{d}^ {5}+2\,{a}^{3}{b}^{6}+9\,{a}^{3}{b}^{5}c+18\,{a}^{3}{b}^{4}{c}^{2}+3\,{a}^{3}{b}^{4}cd
    -6\,{a}^{3}{b}^{4}{d}^{2}+16\,{a}^{3}{b}^{3}{c}^{3}\\
    & \qquad +9\,{a}^{3}{b}^{3}{c}^{2}d-6\,{a}^{3}{b}^{3}c{d}^{2}+24\,{a}^{3}{b}^{2}{c}^{3}d
    -21\,{a}^{3}{b}^{2}{c}^{2}{d}^{2}-33\,{a}^{3}{b}^{2}c{d}^{3}+6\,{a}^{3}{b}^{2}{d}^{4}\\
    & \qquad -6\,{a}^{3}b{c}^{2}{d}^{3}+27\,{a}^{3}bc{d}^{4}-2\,{a}^{3}{d}^{6}-6\,{a}^{2}{b}^{6}c-15\,{a}^{2}{b}^{5}{c}^{2}
    +6\,{a}^{2}{b}^{5}cd-21\,{a}^{2}{b}^{4}{c}^{3}\\
    & \qquad +18\,{a}^{2}{b}^{4}{c}^{2}d+12\,{a}^{2}{b}^{4}c{d}^{2}-24\,{a}^{2}{b}^{3}{c}^{4}+21\,{a}^{2}{b}^{3}{c}^{3}d
    +12\,{a}^{2}{b}^{3}{c}^{2}{d}^{2}-12\,{a}^{2}{b}^{3}c{d}^{3}\\
    & \qquad -6\,{a}^{2}{b}^{2}{c}^{4}d-18\,{a}^{2}{b}^{2}{c}^{2}{d}^{3}-6\,{a}^{2}{b}^{2}c{d}^{4}
    +3\,{a}^{2}b{c}^{2}{d}^{4}+6\,{a}^{2}bc{d}^{5}+6\,a{b}^{6}{c}^{2}+3\,a{b}^{5}{c}^{3}\\
    & \qquad -12\,a{b}^{5}{c}^{2}d-12\,a{b}^{4}{c}^{3}d+
    12\,a{b}^{3}{c}^{5}-3\,a{b}^{3}{c}^{4}d+15\,a{b}^{3}{c}^{3}{d}^{2}
    +12\,a{b}^{3}{c}^{2}{d}^{3}+3\,a{b}^{2}{c}^{4}{d}^{2}\\
    & \qquad -6\,a{b}^{2}{c}^{3}{d}^{3}-6\,a{b}^{2}{c}^{2}{d}^{4}-2\,{b}^{6}{c}^{3}+3\,{b}^{5}{c}^{4}
    +6\,{b}^{5}{c}^{3}d+3\,{b}^{4}{c}^{5}-6\,{b}^{4}{c}^{4}d-6\,{b}^{4}{c}^{3}{d}^{2}\\
    & \qquad -2\,{b}^{3}{c}^{6}-3\,{b}^{3}{c}^{5}d+3\,{b}^{3}{c}^{4}{d}^{2}
    +2\,{b}^{3}{c}^{3}{d}^{3}-6\,{a}^{5}b{d}^{2}+6\,{a}^{5}{d}^{3}+6\,{a}^{4}{b}^{3}d+24\,{a}^{4}{b}^{2}cd\\
    & \qquad +6\,{a}^{4}{b}^{2}{d}^{2}-24\,{a}^{4}bc{d}^{2}-30\,{a}^{4}b{d}^{3}+18\,{a}^{4}{d}^{4}-6\,{a}^{3}{b}^{5}
    -18\,{a}^{3}{b}^{4}c+6\,{a}^{3}{b}^{4}d\\
    & \qquad -24\,{a}^{3}{b}^{3}{c}^{2}+9\,{a}^{3}{b}^{3}cd+12\,{a}^{3}{b}^{3}{d}^{2}+12\,{a}^{3}{b}^{2}{c}^{2}d
    +36\,{a}^{3}{b}^{2}c{d}^{2}-12\,{a}^{3}{b}^{2}{d}^{3}\\
    & \qquad +12\,{a}^{3}b{c}^{2}{d}^{2}-27\,{a}^{3}bc{d}^{3}-6\,{a}^{3}b{d}^{4}+6\,{a}^{3}{d}^{5}
    +12\,{a}^{2}{b}^{5}c+21\,{a}^{2}{b}^{4}{c}^{2}-24\,{a}^{2}{b}^{4}cd\\
    & \qquad +24\,{a}^{2}{b}^{3}{c}^{3}-39\,{a}^{2}{b}^{3}{c}^{2}d-24\,{a}^{2}{b}^{2}{c}^{3}d
    +15\,{a}^{2}{b}^{2}{c}^{2}{d}^{2}+24\,{a}^{2}{b}^{2}c{d}^{3}+3\,{a}^{2}b{c}^{2}{d}^{3}\\
    & \qquad -12\,{a}^{2}bc{d}^{4}-6\,a{b}^{5}{c}^{2}-3\,a{b}^{4}{c}^{3}
    +18\,a{b}^{4}{c}^{2}d-6\,a{b}^{3}{c}^{4}+6\,a{b}^{3}{c}^{3}d-18\,a{b}^{3}{c}^{2}{d}^{2}\\
    & \qquad +6\,a{b}^{2}{c}^{4}d-3\,a{b}^{2}{c}^{3}{d}^{2}+6\,a{b}^{2}{c}^{2}{d}^{3}-6\,{a}^{4}{b}^{2}d+12\,{a}^{4}b{d}^{2}
    -6\,{a}^{4}{d}^{3}+6\,{a}^{3}{b}^{4}\\
    & \qquad +12\,{a}^{3}{b}^{3}c-12\,{a}^{3}{b}^{3}d-24\,{a}^{3}{b}^{2}cd+12\,{a}^{3}bc{d}^{2}+12\,{a}^{3}b{d}^{3}
    -6\,{a}^{3}{d}^{4}-6\,{a}^{2}{b}^{4}c\\
    & \qquad -6\,{a}^{2}{b}^{3}{c}^{2}+18\,{a}^{2}{b}^{3}cd+12\,{a}^{2}{b}^{2}{c}^{2}d-18\,{a}^{2}{b}^{2}c{d}^{2}
    -6\,{a}^{2}b{c}^{2}{d}^{2}+6\,{a}^{2}bc{d}^{3}-2\,{a}^{3}{b}^{3}\\
& \qquad +6\,{a}^{3}{b}^{2}d-6\,{a}^{3}b{d}^{2}+2\,{a}^{3}{d}^{3} \Big).
\end{align*}

\subsection{Class 2.11 : $2\I_2^* + 6\I_2$ with MWG = $\mathbb{Z}/2\mathbb{Z}$}
\begin{equation*}
        y^2 =  \Big( x- dt(t-1) (at-1) \Big)\Big( x- bt(t-1) (ct-1) \Big) \Big( x-t(t-1) (at-1) (ct-1) \Big).
\end{equation*}

\subsection{Class 2.12 : $3\I_0^* + 2\I_2+2\I_1$ with MWG = $\mathbb{Z}/2\mathbb{Z}$}
\begin{equation*}
    \begin{aligned}
        & y^2 = x^3 +9t(t-1) \Big( (ad-bc-2a+2c) t -ad+2a+b+d-2 \Big)x^2\\
        & -81 t^2 (t-1)^2 (at-ct-a+1) \Big( (ad-bc-a+c)  t -(d-1) (a+b-1) \Big) x.
    \end{aligned}
\end{equation*}

\end{document}